\documentclass[11pt]{article}
\usepackage[letterpaper,margin=1in]{geometry}
\usepackage{graphicx} 
\usepackage{amssymb}
\usepackage{amsmath}
\usepackage{amsthm}
\usepackage{enumerate}
\usepackage{epsfig}
\usepackage{tikz}
\usetikzlibrary{patterns.meta}
\usepackage{hyperref}
\hypersetup{hidelinks}
\usepackage{caption}
\usetikzlibrary{arrows.meta,calc,decorations.markings}
\usetikzlibrary{decorations.pathreplacing}

\newtheorem{thm}{Theorem}[section]
\newtheorem{lem}[thm]{Lemma}

\newtheorem{fact}[thm]{Fact}
\newtheorem{cor}[thm]{Corollary}
\newtheorem{df}[thm]{Definition}

\newtheorem{rem}[thm]{Remark}
\newtheorem*{fact*}{Fact}
\newtheorem{claim}[thm]{Claim}
\newtheorem{obs}[thm]{Observation}
\newtheorem{example}[thm]{Example}

\title{Major progress on hamiltonian saturation in hypergraphs}
\author{Andrzej Ruci\'nski\thanks{Research supported by Narodowe Centrum Nauki,  grant 2024/53/B/ST1/00164.}\\
\small Adam Mickiewicz University\\[-0.8ex]
\small Department of Discrete Mathematics\\[-0.8ex]
\small Pozna\'n, Poland\\
\small\tt rucinski@amu.edu.pl\\
\and
Andrzej {\.Z}ak \thanks{Research partially supported
by the Polish Ministry of Science and Higher Education.} \\
\small AGH University \\[-0.8ex]
\small Faculty of Applied Mathematics\\[-0.8ex]
\small Krak\'ow, Poland\\
\small\tt zakandrz@agh.edu.pl
}

\date{\today}

\newsavebox{\reviewfigurebox}

\begin{document}

\maketitle

\begin{abstract}
For $1\leq \ell< k$,  an \emph{$\ell$-overlapping $k$-cycle} is a $k$-uniform hypergraph in which,
for some cyclic vertex ordering, every edge consists of $k$ consecutive vertices and every two
consecutive edges share exactly $\ell$ vertices.

A $k$-uniform hypergraph $H$ is $\ell$-\emph{hamiltonian saturated} if $H$ does not contain an
$\ell$-overlapping Hamilton $k$-cycle but every hypergraph obtained from $H$ by adding a new edge
does contain such a cycle. For $n$ divisible by $k-\ell$, let sat$(n,k,\ell)$ be the smallest number of edges in an
$\ell$-hamiltonian saturated $k$-uniform hypergraph on $n$ vertices. In the case of graphs Clark
and Entringer showed in 1983 that, for sufficiently large $n$, sat$(n,2,1)=\lceil \frac{3n}2\rceil$. By ad hoc methods, the present authors proved previously
that sat$(n,k,\ell)=\Theta(n^{\ell})$ for some rather narrow and scattered ranges of  $\ell$. In this paper we present a unified proof of this equality for all $k$ and $\ell>k/2$.

\end{abstract}

\section{Introduction}

Saturation is a central notion of extremal graph theory, somewhat complementary to Tur\'an numbers.  Given a
fixed graph or hypergraph \(F\), a (hyper)graph $H$ is  \(F\)-saturated if it does not
contain a copy of \(F\), but the addition of any missing edge to $H$ creates
such a copy. The smallest number of edges in an $F$-saturated (hyper)graph on $n$ vertices is called the \emph{saturation number} of $F$. 

A systematic study of saturation goes back to the theorem of
Erd\H{o}s, Hajnal and Moon \cite{EHM} on complete graphs.  Since then, saturation problems
have been studied for many families of graphs and hypergraphs; see, for example,
the survey of Currie, Faudree, Faudree and Schmitt \cite{CFFS}.  In the hypergraph setting
saturation problems are substantially more complex.  A general upper bound of $O(n^{k-1})$ for
the saturation number of a fixed \(k\)-uniform hypergraph was obtained by
Pikhurko \cite{Pikhurko}, but for specific hypergraphs the correct order of magnitude
may be much smaller and usually requires more involved methods.
For instance, for the loose triangle $T$, that is, the $(1,3)$-cycle defined below, with exactly three edges,
English, Kostochka and Zirlin~\cite{EKZ} proved that the saturation number equals $\Theta(n)$.

The present paper is concerned with hamiltonian saturation.  In the graph case,
the smallest maximally non-hamiltonian graphs were studied by Clark and
Entringer~\cite{CE}, and later by Clark, Crane, Entringer and Shapiro~\cite{CCES}
and by Clark, Entringer and Shapiro~\cite{CES}.  Related hamiltonian extremal
questions were also considered by Bondy~\cite{Bondy} and by Wong and Wong
\cite{ww}.

 We will call  a $k$-uniform hypergraph briefly a \emph{$k$-graph}.
For hypergraphs, several different  notions of cycles are
possible, most notably, the Berge cycle (see, e.g.,~\cite{EGGMS} where the saturation number for the family of $k$-uniform Berge
triangles has been determined).

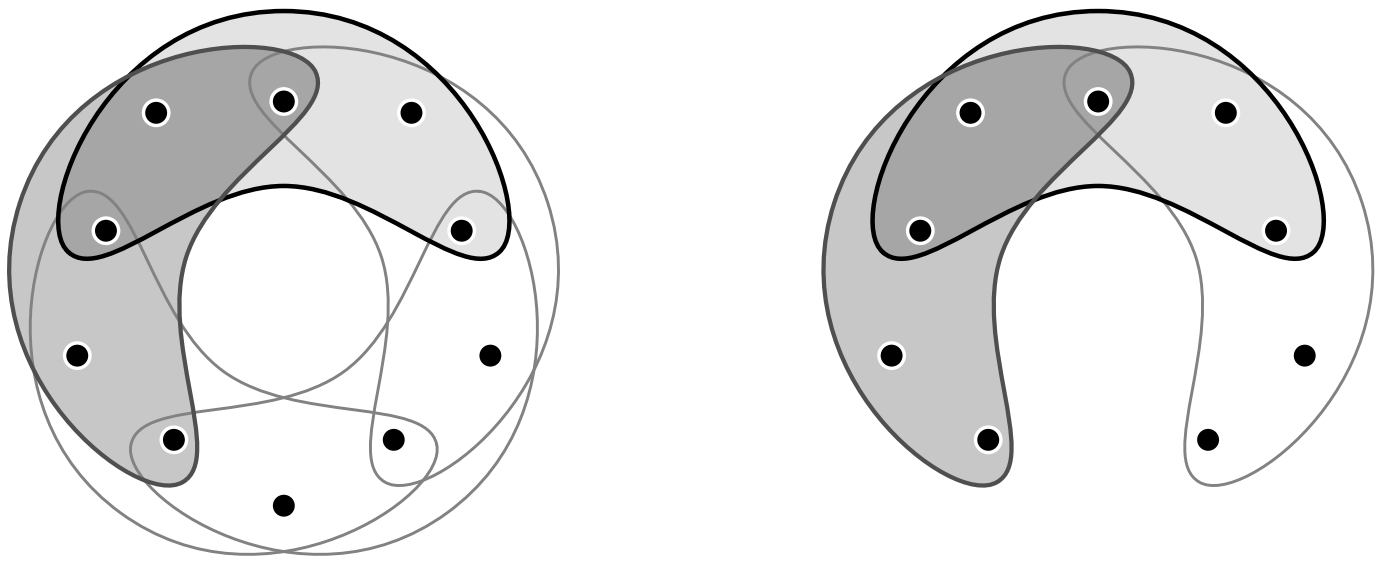
\begin{figure}[htbp]
\centering
\begingroup
\newcommand{\rysKontur}{%
  (-1.07,.18)
  .. controls (-1.27,.35) and (-.89,1.38) .. (0,1.38)
  .. controls (.89,1.38) and (1.27,.35) .. (1.07,.18)
  .. controls (.89,.02) and (.44,.51) .. (0,.51)
  .. controls (-.44,.51) and (-.89,.02) .. cycle}

\newcommand{\rysPanel}[3]{%
  \ifnum#1=1\relax
    \fill[black!11] \rysKontur;
    \fill[rotate=72,black!22] \rysKontur;
    \begin{scope}
      \clip \rysKontur;
      \fill[rotate=72,black!35] \rysKontur;
    \end{scope}
    \tikzset{druga krawedz/.style={draw=black!69,line width=.8pt}}
  \else
    \begin{scope}
      \clip \rysKontur;
      \path[rotate=72,pattern color=black!62,
        pattern={Lines[angle=-27,distance=2.9pt,line width=.32pt]}]
        \rysKontur;
    \end{scope}
    \tikzset{druga krawedz/.style={draw=black!85,line width=.75pt,
      dash pattern=on 3.1pt off 1.6pt}}
  \fi
  \foreach \a in {#2}
    \draw[rotate=\a,black!49,line width=.55pt] \rysKontur;
  \draw[black,line width=.85pt] \rysKontur;
  \draw[rotate=72,druga krawedz] \rysKontur;
  \foreach \i in {#3}{
    \coordinate (v) at
      ({90+36*\i}:{ifthenelse(mod(\i,2)==0,.93,1.08)});
    \fill[white] (v) circle[radius=2.04pt];
    \fill[black] (v) circle[radius=1.55pt];
  }%
}

\newcommand{\rysCyklSciezka}[1]{%
  \begin{tikzpicture}[scale=1.33,line cap=round,line join=round]
    \rysPanel{#1}{144,216,288}{0,...,9}
    \begin{scope}[xshift=4.05cm]
      \rysPanel{#1}{288}{0,1,2,3,4,6,7,8,9}
    \end{scope}
    \pgfresetboundingbox
    \path[use as bounding box] (-1.40,-1.39) rectangle (5.46,1.40);
  \end{tikzpicture}%
}

\rysCyklSciezka{1}\par 
\endgroup

\caption{A $(3,5)$-cycle on ten vertices and a $(3,5)$-path on nine vertices. Two consecutive edges are shaded to exhibit their intersection.}
\label{fig:cycle-and-path}
\end{figure}

In this paper, however, we study a more symmetric notion, introduced by G.Y. Katona and Kierstead in~\cite{kk} and popularized in a series of papers by R\"odl, Ruci\'nski and Szemer\' edi, in the context of Dirac-type questions (see, e.g., \cite{RR}). For \(1\le \ell<k\), an
\emph{\(\ell\)-overlapping \(k\)-cycle}, or briefly an \emph{\((\ell,k)\)-cycle} is a \(k\)-graph whose vertices can
be cyclically ordered so that every edge consists of \(k\) consecutive vertices,
and every two consecutive edges overlap in exactly \(\ell\) vertices.  Note that the number of vertices in an \((\ell,k)\)-cycle is divisible by \(k-\ell\).
An \emph{\(\ell\)-overlapping \(k\)-path}, or briefly an
\emph{\((\ell,k)\)-path}, is defined analogously using a linear ordering of its
vertices (see Figure~\ref{fig:cycle-and-path} for an example of each).

An \((\ell,k)\)-cycle or an \((\ell,k)\)-path contained in a $k$-graph $H$ is called \emph{Hamilton}  if it is a spanning sub-hypergraph of $H$.
A \(k\)-graph is \emph{\(\ell\)-hamiltonian} if it contains a Hamilton \((\ell,k)\)-cycle.
This
notion, for \(\ell=k-1\), reduces to hamiltonian chains in hypergraphs studied by
G. Y. Katona and Kierstead~\cite{kk}; see also G. Y. Katona's survey~\cite{k-sur}.

\medskip

For \(n\) divisible by \(k-\ell\), let \(\mathrm{sat}(n,k,\ell)\) denote the saturation number of a Hamilton $(\ell,k)$-cycle, that is, the minimum number of edges in an
\(\ell\)-hamiltonian saturated \(k\)-graph on \(n\) vertices. Clark and Entringer~\cite{CE} proved that   sat$(n,2,1)=\lceil \frac{3n}2\rceil$ for $n\ge52$.

 \.Zak  in~\cite{zak} (for \(\ell=k-1\)) and  Ruci\'nski and \.Zak in \cite{RZ,RZ2,RZ3,RZ4} initiated a systematic study of this function for hypergraphs.  They proved that
\[
        \mathrm{sat}(n,k,\ell)=\Theta(n^\ell),
\]
where the hidden constant depends on $k$ and $\ell$ only,
for several ranges of the parameters \(k,\ell\), including $0.8k\le \ell<k$ and $(k-1)/3\le \ell<(k-1)/2$, and conjectured that the same
formula should hold for all \(1\le \ell<k\).   In this paper we complete the determination of the order of magnitude of
\(\mathrm{sat}(n,k,\ell)\) in the upper half of the range of parameter \(\ell\). In this range even non-consecutive edges of a $(\ell,k)$-path or cycle may intersect -- a source of potential difficulty.
\begin{thm}\label{main_result}
For all \(k\geq 3\) and \(k/2\leq \ell \leq k-1\), assuming \(n\) is divisible by \(k-\ell\),
\[
        \mathrm{sat}(n,k,\ell)=\Theta(n^\ell).
\]
\end{thm}

The cases $\ell=k/2$ and $\ell=k-1$ were proved in \cite{RZ3} and \cite{zak}, respectively, 
thus, for $k\in\{3,4\}$, all cases for $\ell$ are already covered. Consequently,  throughout the paper we assume  that $k\ge5$ and $\ell>k/2$.

We anticipate that the other half of the range, $2\le \ell<k/2$, will require a different, possibly simpler approach, and we intend to treat it
in a separate  paper~\cite{LRZ}.

The case $\ell=0$, not studied in this paper, corresponds to a perfect $k$-matching (provided $k$ divides $n$) and may be viewed as a degenerate case of a Hamilton $(\ell,k)$-cycle. Thus, one may define
$\mathrm{sat}(n,k,0)$ as the minimum number of edges in an $n$-vertex
$k$-graph with no perfect $k$-matching, in which the addition of any missing
edge creates such a matching. For graphs, it was observed  in~\cite{HKNS} that, as a consequence of a structural result of Mader,  $\mathrm{sat}(n,2,0)=\Theta(n^{3/2})$. The question of what is the order of magnitude of
$\mathrm{sat}(n,k,0)$ for a fixed $k\ge3$ is widely open and, most likely, will require a completely different approach than the one presented here.

In the entire range $1\le\ell\le k-1$, the lower bound follows  by a  standard double counting argument (see \cite[Proposition 2.1]{RZ}).
Therefore, the proof of Theorem~\ref{main_result} is entirely concentrated on
 the upper bound, for which all we need is just a construction of an \(\ell\)-hamiltonian saturated \(k\)-graph $H$ on \(n\) vertices and with $O(n^\ell)$ edges.
For notational reasons, we now switch from \(n\) to \(N\) in the statement of Theorem~\ref{main_result} (to free letter \(n\) to denote the number of vertices in an auxiliary graph).

The construction, indeed, starts from a bounded-degree maximally non-hamiltonian graph
\(G_1\) on \(n\) vertices, extends it to another graph $G$ on $2n$ vertices, upon which one builds a $k$-graph $H$ on $N$ vertices with the required property.  The edge set of
\(H\) is divided into three parts.  Roughly, part \(H_1\) follows the graph
structure of \(G_1\), while part \(H_2\) contains edges which ``prefer''  particular classes of a prescribed vertex set partition.
The third part of \(H\) is non-constructive, but contained in a well-defined $k$-graph $H_3$, that is, $H_1\cup H_2\subseteq H\subseteq H_3$ -- see the details in Section~\ref{main_con}.

Another difficulty in the large-overlap range \(\ell>k/2\), besides the intersection of non-consecutive edges, is the
loss of sufficiently precise local control. Function~\(\nu\) introduced
below (and also used in previous articles on this topic) bounds the length of an individual \(H_2\)-path in terms of the number of
available vertices in its assigned class. For \(\ell>k/2\), however,
determining \(\nu\) with the required precision becomes unexpectedly difficult.
Indeed, even a small error in a local estimate
may  produce an
unbounded accumulated error when applied to many paths at once.

The main new idea here is therefore to replace this local ``bookkeeping'' by a global
one. To this end, we introduce  a function \(\Phi\) which controls the total contribution
of all relevant \(H_2\)-paths simultaneously, optimizing over admissible
distributions of the available vertices. 
In this way, errors that would accumulate linearly when
applying  \(\nu\) separately to each path, are replaced by a constant global error.

The paper is organized as follows.  In Section~\ref{fniu} we define  functions \(\nu\) and \(\Phi\), the latter based on admissible sequences, and  collect all their relevant properties, the proofs of which are deferred to the Appendix.
 In Section~\ref{main_con}, for every sufficiently large $N$, we   construct, via auxiliary $k$-graphs $H_1$, $H_2$, and $H_3$, a non-$\ell$-hamiltonian  $k$-graph \(H\)  on $N$ vertices. In  Section~\ref{notHam}  we prove that $H_1\cup H_2$ is itself non-hamiltonian (which, by construction, is required for $H$ to be well-defined).  Then, in Section~\ref{eHam}, we complete the proof of Theorem~\ref{main_result}
by showing that adding any new edge $e$ to $H$ creates a Hamilton
\((\ell,k)\)-cycle (again, by construction, it suffices to prove it for $(H_1\cup H_2)+e$ and $e\not\in H_3$). The final Section~\ref{sec:correct-loads} is devoted to proving a quite technical, but crucial lemma which, in turn, calls upon two further lemmas. One of these lemmas, Lemma~\ref{lem:matching-stability},  establishes stability of star matchings in bipartite graphs under local perturbations, and may be of independent interest.


\section{Two pivotal functions}\label{fniu}
In this section, we define two functions related to $(\ell,k)$-paths which turn out to be instrumental in the forthcoming proofs. 

\medskip

\subsection{Function $\nu$}\label{2.1}
Our proofs will heavily depend on $(\ell,k)$-paths, where each edge includes at least $k-\ell+1$ vertices from a predetermined, relatively small set, and the rest from a significantly larger set.
 Here is a formal definition.

 \begin{df}[dominant path]\label{dom}\rm Let $U$ and $W$ be two disjoint sets. An $(\ell,k)$-path $P$ with $V(P)\subseteq U\cup W$  is called \emph{$U$-dominant} if
 \begin{itemize}
 \item[(i)] $U\subseteq V(P)$,
\item[(ii)] $|e\cap U|\ge k-\ell+1\qquad\mbox{for all}\qquad e\in P$,
\item[(iii)] $P$ begins and ends with $\ell$ vertices from $U$.
\end{itemize}
\end{df}
\noindent The set $U$ will be referred to as \emph{mandatory}, while the set $W$ as \emph{the reservoir}.
\begin{rem}\label{remark}\rm
Note that if $P$ satisfies assumptions (ii) and (iii) only,
then $P$ is $U'$-dominant, where $U'=U\cap V(P)$.
\end{rem}
We now state a simple upper bound on the length of an $(\ell,k)$-path which satisfies just the single assumption (ii) above. Its proof can be found in Appendix A1.

\begin{obs}\label{obs1}
Let $U$ be a finite set and $P$ be an $(\ell,k)$-path satisfying just assumption (ii) in Definition~\ref{dom}. Then
$$|V(P)|\le \ell|U|.$$
\end{obs}

 A pivotal role in our proofs is played by a function determining the greatest possible length of a $U$-dominant path with access to an infinite reservoir set. Throughout, we adopt a convention of indicating infinite sets by the asterisk $^*$. \space

\begin{df}[function $\nu$]\label{defniu}\rm  Given an integer $\alpha\ge2\ell$ and disjoint sets $U$ and $W^*$ with $|U|=\alpha$ and $|W^*|=\infty$, 
$$\nu(\alpha)=\max\{|V(P)|:\; V(P)\subseteq U\cup W^* \quad\mbox{and $P$ is $U$-dominant}\}.$$
\end{df}

In the lemma below we collect all the properties of function $\nu$ we will need throughout the paper.
 We defer their proofs to Appendix A1.

\begin{lem}\label{nu1}

For parts (a)-(e), let $\alpha\ge2\ell$. Then \space

\begin{enumerate}
\item[(a)] $\nu(\alpha)\ge \alpha+1$,

\item[(b)] $\nu(\alpha)\le\ell\alpha$,

\item[(c)] $\nu(\alpha+1)\ge\nu(\alpha)$,

\item[(d)] \( \nu(\alpha+k-\ell)\ge \nu(\alpha)+k-\ell. \),

\item[(e)] 
$\nu(\alpha)\ge \nu(\alpha+1)-(k-1)$.

\item[(f)] Moreover, $\nu(\alpha)-\alpha\to \infty$, as $\alpha\to\infty$.
\end{enumerate}

\end{lem}

Properties (a), (b), and (f) compare $\nu(\alpha)$ with $\alpha$, properties (c) and (d) say that $\nu$ is a non-decreasing function of $\alpha$ which, in the long run, increases steadily, while property (e) limits the rate of increase of $\nu$ from above.

\subsection{Function $\Phi$}\label{secPhi}
In the proof of Theorem~\ref{main_result} we shall also need a multi-path variant of \(\nu\), allowing us to
control several \(H_2\)-paths simultaneously when a prescribed total number of
vertices from the mandatory sets is unavailable. To define it, we  need two other definitions.

First, we introduce a fundamental notion of a gap.

\begin{df}\label{defgap}\rm Given $k\ge5$ and $\ell>k/2$, let
the \emph{gap} $g:=g(k,\ell)$
be the number of vertices between any  two
\emph{consecutive} \emph{disjoint} edges belonging to an $(\ell,k)$-path (or $(\ell,k)$-cycle).
\end{df}

\begin{figure}[htbp]
\centering
\begin{tikzpicture}[x=.60cm,y=1cm,every node/.style={font=\small},
  end edge/.style={draw,rounded corners=4pt,fill=black!5,line width=.75pt},
  bracket/.style={line width=.75pt,line cap=round}]
\draw[end edge] (.68,-.18) rectangle (7.32,.18);
\draw[end edge] (9.68,-.18) rectangle (16.32,.18);
\draw[bracket] (4,.22) -- (4,.56) -- (10,.56) -- (10,.22);
\draw[bracket] (7,-.22) -- (7,-.91) -- (13,-.91) -- (13,-.22);
\node[above] at (2.45,.19) {$e_1$};
\node[above] at (7,.56) {$e_2$};
\node[below] at (10,-.91) {$e_3$};
\node[above] at (13.85,.19) {$e_4$};
\foreach \i in {1,...,16}
  \fill (\i,0) circle[radius=1.75pt];
\foreach \i in {8,9}
  \draw[fill=white,line width=.7pt] (\i,0) circle[radius=2.1pt];
\draw[decorate,decoration={brace,mirror,amplitude=3pt}]
  (7.68,-.15) -- (9.32,-.15);
\node[below] at (8.5,-.28) {$g=2$};
\end{tikzpicture}
\caption{Four consecutive edges of a $(4,7)$-path. The first edge disjoint
from $e_1$ is $e_4$, so $g=2$.}
\label{fig:gap}
\end{figure}
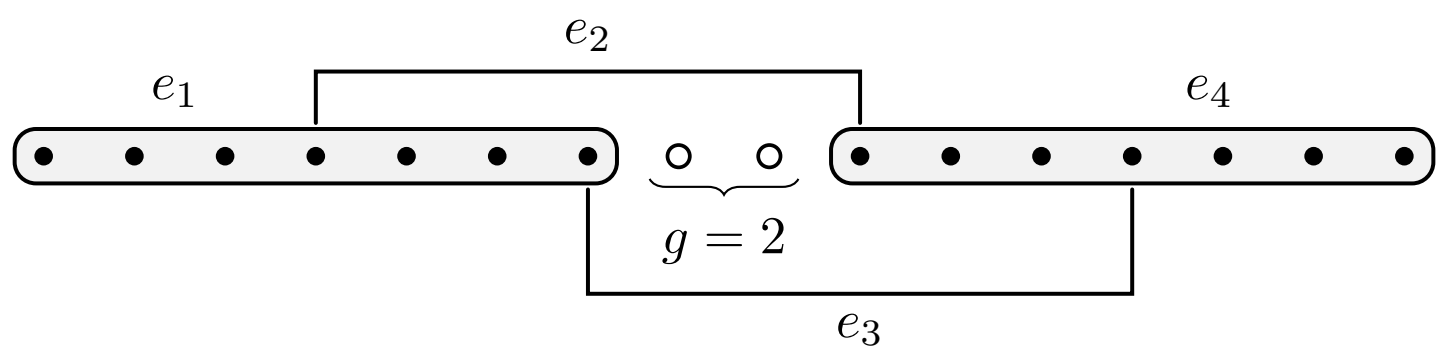

Notice
that gap is well defined, as it does not depend on a particular choice of two consecutive disjoint edges. Moreover, there is an explicit formula for it,
\begin{align}\label{gap}
0\le g=\left\lceil \frac{k}{k-\ell}\right\rceil(k-\ell)-k<k-\ell<k,
\end{align}
so that $g=0$ if and only if $k-\ell$ divides $k$ (see Figure~\ref{fig:gap}).

Next, we define admissible sequences of integers.

\begin{df}\label{defadm}\rm
Fix $k$ and $\ell$ as in Theorem~\ref{main_result} and let $m\ge2$ be an integer.
 We say that an integer sequence
\((r_i)_{i=1}^m\) is \emph{admissible} if
$$0\le r_i\le k-\ell,\;\;\;i=1, \dots,m$$
and there exists another integer sequence,
\((s_i)_{i=1}^{m-1}\), with all 
$$ s_i\in\{g,g+k-\ell,g+2(k-\ell)\},$$ called \emph{witness of admissibility}, or just \emph{witnessing},
such that, for every \(t=1,\ldots,m-1\),
\begin{equation}\label{adm}
        \sum_{i=1}^{t-1} s_i
        \le
        \sum_{i=1}^t r_i
        \le
        \sum_{i=1}^{t} s_i,
\quad\mbox{
and}\quad
\sum_{i=1}^{m-1}s_i=\sum_{i=1}^m r_i
\end{equation}
We denote by
$\mathcal A_m$
the family of all admissible sequences of length \(m\).
The family $\mathcal A_m$ is nonempty: $(g,\ldots,g,0)$ is witnessed by $s_i=g$. We also set $\mathcal A_1=\{(0)\}$ and let $\mathcal A_0$ consist of the empty sequence.

\end{df}

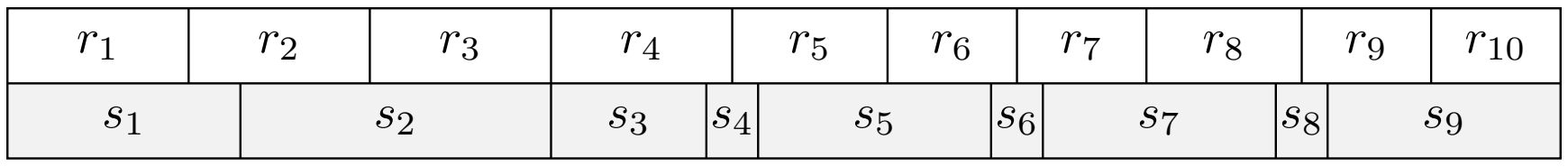
\begin{figure}[htbp]
\centering
\begin{tikzpicture}[x=.19cm,y=1cm,font=\small,line width=.45pt]
  \fill[black!5] (0,-.55) rectangle (60,0);
  \fill[black!5] (27,-.55) rectangle (29,0);
  \draw (0,-.55) rectangle (60,.55);
  \draw (0,0) -- (60,0);
  \foreach \a/\b/\i in {0/7/1,7/14/2,14/21/3,21/28/4,28/34/5,
      34/39/6,39/44/7,44/50/8,50/55/9,55/60/10}{
    \node[inner sep=0pt] at ({(\a+\b)/2},.275) {$r_{\i}$};
    \ifnum\i<10 \draw (\b,0) -- (\b,.55);\fi
  }
  \foreach \a/\b/\i in {0/9/1,9/21/2,21/27/3,27/29/4,29/38/5,
      38/40/6,40/49/7,49/51/8,51/60/9}{
    \node[inner sep=0pt] at ({(\a+\b)/2},-.275) {$s_{\i}$};
    \ifnum\i<9 \draw (\b,-.55) -- (\b,0);\fi
  }
\end{tikzpicture}
\caption{An admissible sequence $(r_i)$ and its witness $(s_i)$ ($m=10$). The numbers $r_i$ and $s_i$ are represented symbolically by the widths of the bricks (rectangles).}
\label{fig:admissible-bricks}
\end{figure}

The intuition behind conditions~\eqref{adm} is that it allows to line up two layers of ``bricks'' of lengths $r_i$ and $s_i$, the latter beneath the former, so that each brick in one row touches only at most two bricks in the other; the equality of the total lengths leaves no unused part of any brick (see Figure~\ref{fig:admissible-bricks}).

Note that condition~\eqref{adm} implies that, for all $1\le a\le b\le m$,
\begin{equation}\label{follows}
\sum_{i=a}^{b-1}s_i\le\sum_{i=a}^b r_i\le(b-a+1)(k-\ell).
\end{equation}
Indeed,
$$\sum_{i=a}^{b-1}s_i=\sum_{i=1}^{b-1}s_i-\sum_{i=1}^{a-1}s_i\le\sum_{i=1}^{b} r_i-\sum_{i=1}^{a-1}r_i=\sum_{i=a}^b r_i.$$

Now, we are ready to define a multi-path variant of function \(\nu\).

\begin{df}[function \(\Phi\)]\label{defphi}\rm
Fix $k$ and $\ell$. For integers \(\alpha\ge k+\ell\) and \(m\ge0\), define
\begin{align*}
\Phi(\alpha,m)
=
\max\left\{
\sum_{i=1}^m\bigl(\nu(\alpha-r_i)+r_i\bigr):\;(r_i)_{i=1}^m\in \mathcal A_m
\right\}.
\end{align*}
\end{df}
\noindent
Thus $\Phi(\alpha,0)=0$ and $\Phi(\alpha,1)=\nu(\alpha)$.

Although, in general, Definition~\ref{defadm} allows three distinct values of $s_i$, for extremal sequences the largest of them can be avoided. The proof of the next claim can be found in Appendix A2.

\begin{claim}\label{just2}
The maximum defining $\Phi(\alpha,m)$ can always be attained
with a witnessing sequence satisfying
\begin{equation}\label{eq:phi-two-value-witness}
s_i\in\{g,g+k-\ell\},\qquad 1\le i<m.
\end{equation}
\end{claim}

It turns out that the growth of $\Phi(\alpha,m)$, as a function of $m$, oscillates around $\nu(\alpha)$.
\begin{lem}\label{phi-step}
There exists a constant $c$ such that for all $\alpha\ge k+\ell$
and \(m\ge1\)
\begin{equation}\label{phi-dif-bounds}
   \nu(\alpha)-c \le     \Phi(\alpha,m)-\Phi(\alpha,m-1)\le \nu(\alpha)+c.
\end{equation}
It follows that 
\begin{equation}\label{phi-bounds}
m(\nu(\alpha)-c)\overset{(i)}{\le} \Phi(\alpha,m)\overset{(ii)}{\le} m(\nu(\alpha)+c).
\end{equation}
\end{lem}


The next lemma  establishes super-additivity of $\Phi$ with respect to the second variable, up to an additive constant.

\begin{lem}\label{phi-glue}
There exists a constant $c$ such that, for all $\alpha\ge k+\ell$ and all
positive integers $a,b$,
\[
        \Phi(\alpha,a)+\Phi(\alpha,b)
        \le
        \Phi(\alpha,a+b)+c .
\]
Consequently, for every sequence of positive integers $m_1,\ldots,m_q$,
\[
        \sum_{s=1}^q\Phi(\alpha,m_s)
        \le
        \Phi\left(\alpha,\sum_{s=1}^q m_s\right)+(q-1)c .
\]
\end{lem}

Proofs of Lemmas~\ref{phi-step} and~\ref{phi-glue} are deferred to Appendix A2.

\section{Main construction and outline of proof}\label{main_con}

Fix $k\ge5$ and $\ell>k/2$. Throughout, we identify a $k$-graph, $k\ge3$, with its edge set when no confusion arises; for ordinary graphs $G$, however, we still write $E(G)$. Using symbols $O(...)$ and  $\Theta(...)$ we will typically mean that the constant hidden there depends on $k$ and $\ell$ only. If, in addition, the constant will depend on another parameter, say $\alpha$, then we will reveal it by writing $O_\alpha(...)$, etc.
 
 In this section, after setting up several constants, we  construct, for a sufficiently large $N$, a $k$-graph $H$ with $|V(H)|=N$ and  $|H|=O(N^{\ell})$ and give a proof that $H$ is $\ell$-hamiltonian saturated, based on two main lemmas which themselves will be proved in the remainder of the paper.

\subsection{Setting up the vertex set}
Let $c_3$ and $c_4$ be the constants appearing in, respectively, Lemma~\ref{phi-step} and Lemma~\ref{phi-glue}. It can be deduced from their proofs (see Appendix A2) that $c_3<k^2$ and $c_4<3k^4$.

  Set $c_-=3c_4+46k$ and $c_+=2k^4+2k^3+5k^2$,  and put
\begin{equation}\label{defh}
h:=1+\left\lceil\frac{c_-+c_+}{k-\ell}\right\rceil.
\end{equation}
Next, choose $\alpha\ge7k$ so large that  
\begin{align}
(i)\;\; \nu(\alpha)>c_3+c_-+c_++(k-\ell)(kh+1)
\quad\&\quad (ii)\;\;\nu(\alpha)-\alpha>c_3+2k, \label{nu_i_alpha}
\end{align}
where function $\nu$ has been defined in  Section~\ref{2.1} (see Definition~\ref{defniu}).
Such a choice of $\alpha$ is possible by Lemma~\ref{nu1}(f).

 Recall Definition~\ref{defphi} of function $\Phi$ and, for every sufficiently large integer $N$ divisible by $k-\ell$, define $n$ to be the smallest integer such that
$N\le\Phi(\alpha,n)-c_+$. Note that, by~\eqref{phi-bounds}(ii),  
\begin{equation}\label{Nn}
n\ge\frac{\Phi(\alpha,n)}{\nu(\alpha)+c_3}\ge\frac{N+c_+}{\nu(\alpha)+c_3},
\end{equation}
thus, given $k,\ell$, and $\alpha$, the larger  $N$, the larger  $n$. In particular, we take $N$ so large that $h\le n$. Also, it follows from the definition of $n$ that
\begin{equation}\label{Nbetween}
\Phi(\alpha,n-1)-c_+<N\le \Phi(\alpha,n)-c_+.
\end{equation}
 

The following pair of inequalities (see Claim~\ref{NN} below), which adjust the ones just stated, are of the greatest importance for us. They will serve as punch lines in the proofs of Lemmas~\ref{main1} (see~\eqref{VCN}) and~\ref{main2} (see~\eqref{NprimN}).
For their succinct formulation, we need another  parameter, $h_N$, whose role  is explained below. Let
\[
h_N=\begin{cases}
0,&N>\Phi(\alpha,n-1)+c_-,\\
h,&N\le\Phi(\alpha,n-1)+c_-.
\end{cases}
\]

\begin{claim}\label{NN} For every sufficiently large $N$ divisible by $k-\ell$,
\begin{equation}\label{nn} 
\Phi(\alpha,n-1)+c_--(k-\ell)\max\{h_N-1,0\}
<N\le\Phi(\alpha,n)-c_+-(k-\ell)kh_N.
\end{equation}
\end{claim}
\noindent
With $h_N=0$, inequalities~\eqref{nn} leave the interval
\[
  (\Phi(\alpha,n-1)-c_+,\,\Phi(\alpha,n-1)+c_-]
\]
uncovered.
Setting $h_N=h$ shifts the lower endpoint in~\eqref{nn} to the left by
$(h-1)(k-\ell)\ge c_-+c_+$, the width of that interval, by~\eqref{defh}.
The upper endpoint moves farther to the left, by $kh(k-\ell)$.
Condition~\eqref{nu_i_alpha}(i) keeps it above $\Phi(\alpha,n-1)+c_-$.
We provide a short proof for completeness.

\begin{proof}[Proof of Claim~\ref{NN}]
This is immediate if $h_N=0$ (the right-hand side by~\eqref{Nbetween}, the left-hand side by the definition of $h_N$). Otherwise, by the left side of~\eqref{Nbetween},
$N>\Phi(\alpha,n-1)-c_+$, so the inequality $(h-1)(k-\ell)\ge c_-+c_+$ (see~\eqref{defh}) proves the lower bound. For the upper bound note that, by the definition of $h_N$, Lemma~\ref{phi-step}, and~\eqref{nu_i_alpha}(i),
\[
N\le\Phi(\alpha,n-1)+c_-\le \Phi(\alpha,n)-\nu(\alpha)+c_3+c_-
<\Phi(\alpha,n)-c_+-(k-\ell)k h_N.
\]
 \end{proof}

\medskip

The construction of a $k$-graph $H$ will be laid upon the following partition of the vertex set $U=V(H)$ of size $N$.
To this end, put 
\begin{equation}\label{alphai}\alpha_i=\begin{cases}\alpha-k+\ell\quad\mbox{for}\quad i\le h_N,\\\alpha\quad\mbox{for}\quad i>h_N.\end{cases}
\end{equation}

\begin{figure}[htbp]
\centering
\begin{tikzpicture}[x=0.98cm,y=0.88cm,font=\small,
  part/.style={draw=black!70,line width=.45pt},
  measure/.style={decorate,decoration={brace,amplitude=3pt},thin}]
  \foreach \x/\idx/\bottom in
    {0/1/-1.07,2/{h_N}/-1.07,3/{h_N+1}/-1.63,5/n/-1.63,
     6.16/{n+1}/-1.27,8.16/{2n}/-1.17}{
    \ifdim\x pt<6pt \def\split{0.00}\else\def\split{0.42}\fi
    \path[part,fill=black!12] (\x,1.15) rectangle ++(1,{\split-1.15});
    \path[part,fill=white] (\x,\split) rectangle ++(1,{\bottom-\split});
    \node at (\x+.5,1.44) {$U_{\idx}$};
    \node at (\x+.5,{(1.15+\split)/2}) {$A_{\idx}$};
    \node at (\x+.5,{(\split+\bottom)/2}) {$B_{\idx}$};
  }
  \foreach \x/\split/\bottom in {1/0/-1.07,4/0/-1.63,7.16/.42/-1.22}{
    \node at (\x+.5,1.44) {$\cdots$};
    \node at (\x+.5,{(1.15+\split)/2}) {$\cdots$};
    \node at (\x+.5,{(\split+\bottom)/2}) {$\cdots$};
  }
  \draw[measure] (-.10,0) -- (-.10,1.15)
    node[midway,left=4pt] {$k+\ell$};
  \draw[measure] (9.26,1.15) -- (9.26,.42)
    node[midway,right=4pt] {$k-1$};
  \draw[measure] (9.26,.42) -- (9.26,-1.17)
    node[midway,right=4pt] {$b_{2n}$};
  \draw[measure] (-.10,-1.07) -- (-.10,0)
    node[midway,left=4pt] {$\alpha-2k$};
  \draw[measure] (6,-1.76) -- (3,-1.76)
    node[midway,below=4pt] {$|U_i|=\alpha$};
  \draw[densely dashed,black!50] (0,-1.63) -- (3,-1.63);
  \draw[measure] (-.10,-1.63) -- (-.10,-1.07)
    node[midway,left=4pt] {$k-\ell$};
\end{tikzpicture}
\caption{The partition of vertex set $U=\bigcup U_i$.}
\label{fig:vertex-partition}
\end{figure}
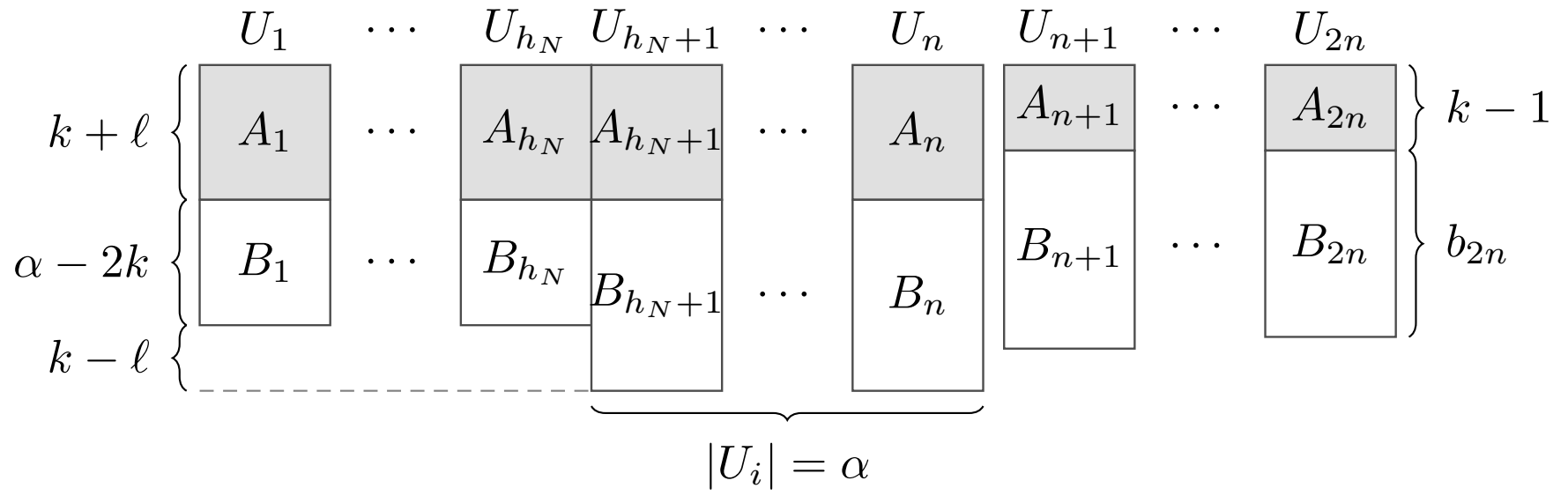

Let $A_i$, $B_i$, $i=1,\dots,2n$, be a family of $4n$ pairwise disjoint sets of sizes
\begin{align}\label{mocai2026}
|A_i|=\begin{cases} k+\ell\mbox{ for } i=1,...,n \\
k-1 \mbox{ for }i=n+1,...,2n, \end{cases}
\end{align}
and
\begin{align}\label{mocbi2026}
|B_i|=\begin{cases} \alpha_i-k-\ell\mbox{ for } i=1,...,n \\
b_i\quad\qquad\;\;\mbox{ for }i=n+1,...,2n, \end{cases}
\end{align}
where the $b_i$'s differ from each other by at most one and are chosen in such a way that
\begin{align}\label{total}
\sum_{i=1}^{2n}(|A_i|+|B_i|)=N.
\end{align}
For each $i=1,\dots,2n$, set $U_i=A_i\cup B_i$ and note that for $i\le n$ we have $|U_i|=\alpha_i$ (see Figure~\ref{fig:vertex-partition}). Our original intention was to have $|U_i|=\alpha$ for every $i\le n$. Diminishing the first $h_N$ sets $B_i$ by $k-\ell$ elements each allows us to obtain the lower threshold in~\eqref{nn} in the proof of Lemma~\ref{main1}.

Next, we estimate the sizes $b_i=|B_i|$, $i=n+1,\dots, 2n$.
By the minimality of $n$ and Lemma~\ref{phi-step},
\[
N>\Phi(\alpha,n-1)-c_+
 \ge(n-1)(\nu(\alpha)-c_3)-c_+.
\]
Consequently, the total size prescribed for
$B_{n+1},\ldots,B_{2n}$ satisfies
\[
\begin{aligned}
\sum_{i=n+1}^{2n}b_i
&=N-n\alpha+(k-\ell)h_N-n(k-1)\\
&>n\bigl(\nu(\alpha)-\alpha-c_3-k+1\bigr)-(\nu(\alpha)-c_3+c_+)\\
&>(k+1)n-(\nu(\alpha)-c_3+c_+),
\end{aligned}
\]
where  the last inequality follows from our choice of $\alpha$.
It follows that, for sufficiently large $N$ (and thus $n$), the balanced choice of
the integers $b_{n+1},\ldots,b_{2n}$ gives
\begin{align}\label{moc_bi}
\min_{n<i\le2n}b_i
=\left\lfloor\frac1n\sum_{i=n+1}^{2n}b_i\right\rfloor
\ge k.
\end{align}
In particular, $b_i\ge k$, for all $i>n$.

On the other hand, 
since the $b_i$'s differ by at most one and their sum is at
most $N$, by~\eqref{Nn}, we obtain
\begin{equation}\label{biO}
b_i\le\left\lceil\frac{N}{n}\right\rceil
=O_\alpha(1)
\qquad\text{for each }n<i\le2n.
\end{equation}
Together with the previous estimates for $i\le n$, we conclude that
$|U_i|=O_\alpha(1)$ for every $i\in[2n]$.


\subsection{Construction}
 The desired $k$-graph $H$ will be defined on the $N$-vertex set
$$U=\bigcup_{i=1}^{2n} U_i.$$

Our construction will be also  based on a graph $G$ obtained as follows. Take any maximally non-hamiltonian graph $G_1$ on the vertex set $\{1,\dots,n\}$ with bounded maximum degree $\Delta(G_1)$. Combining results from~\cite{CES} and~\cite{Bondy}, for each $n\ge52$, one may assume that $\Delta(G_1)\le 5$ (see~\cite[Corollary 2.6]{RZ}).
To form $G$, add to $G_1$  vertices $\{n+1,\dots,2n\}$ and join $i$ to $n+i$, $i = 1,\dots,n$. Thus, 
$$V(G)=[2n],\quad E(G)=E(G_1)\cup\{\{i,n+i\}:\;i=1,\dots,n\},\quad\mbox{and}\quad
\Delta(G)\le6.$$ Note also that, quite trivially, $G_1$ does not have isolated vertices, so the only pendant vertices of $G$ are $n+1,\dots,2n$.

Prior to specifying the edge set of $H$, we introduce additional terminology and notation, which we illustrate by means of an example.
For a subset $S \subset U$, define its \emph{trace} as
\[
\operatorname{tr}(S) = \{ i : S \cap U_i \neq \emptyset \} \quad \text{and} \quad \min(S) = \min\{ i \in \operatorname{tr}(S)\}.
\]
In particular, for $u\in U$, we have $\operatorname{tr}(\{u\})=\{i\}$, where $u\in U_i$.

We  construct  the ultimate $k$-graph $H$  by means of three other hypergraphs: $H_1,H_2$, and $H_3$. Set
$$A:=\bigcup_{i=1}^{2n} A_i$$
for convenience,
 and let

$$H_1=\left\{e\in \binom{A}{k}:\operatorname{tr}(e)\in E(G)   \right\}.$$

Defining $H_2$ requires some preparations. Set
\[ \mathcal{L} = \left\{ L\in\binom{A}{\ell}: \operatorname{tr}(L)\subseteq \{x,y\} \text{ for some }xy\in E(G) \right\} \]
and 
$$\mathcal K=\left\{f\in\binom Uk: \forall L\in\binom{f\cap A}{\ell}\cap\mathcal{L}:\; L\subseteq A_{\min(f)} \right\}. $$
Now,
\[ H_2 := \left\{ f\in\mathcal K:  \left|f\cap U_{\min(f)}\right|\ge k-\ell+1\right\}.
\]

Before defining $H_3$, we make some observations about $H_2$.

\begin{fact}\label{minH2}
If $Q$ is an $(\ell,k)$-path in $H_2$, then there exists an index $i$ such that every edge of $Q$ contains at least $k-\ell+1$ vertices from $U_i$.
\end{fact}
\proof Let $e,e'\in Q$ be two edges with $|e\cap e'|=\ell$, and set $i=\min(e)$. Since $|e\cap U_i|\ge k-\ell+1$, it follows that $|e'\cap U_i|\ge 1$. Consequently, $i\in \operatorname{tr}(e')$, and therefore $\min(e')\le \min(e)$. By symmetry, we also obtain $\min(e)\le \min(e')$. Hence, $\min(e')=\min(e)=i$. By transitivity, this implies that $\min(f)=i$ for every $f\in Q$. \qed

\medskip

Thus, for any $H_2$-path $Q$, we may define
\begin{align}\label{minP}
    \min(Q) = \min(e) \quad \text{for some/every } e \in Q.
\end{align}

\begin{cor}\label{boundH2}
For every $(\ell,k)$-path $Q$ in $H_2$ we have $|V(Q)|=O_\alpha(1)$.
\end{cor}
\proof Note that, in view of~\eqref{minP}, $Q$ satisfies assumption (ii) in Definition~\ref{dom} (with $U_i$, where $i=\min(Q)$), and, thus,  by Observation~\ref{obs1}, $|V(Q)|\le\ell\max_{i\in[2n]}|U_i|$. Since $|U_i|\le \alpha_i$ for $i\le n$ and $|U_i|=k-1+b_i$ for $i>n$, the desired bound follows by~\eqref{biO}\qed

\medskip

For the final ingredient of the construction, $H_3$, we need yet another piece of notation. First, for a graph $F$, set $\operatorname{comp}(F)$ for the number of connected components in $F$ and, for a vertex subset $W \subset V(F)$, let $F[W]$ denote the subgraph of $F$ induced by $W$. Returning to our construction, for $S\subseteq U$, let
$$c(S)=\operatorname{comp}(G^3[\operatorname{tr}(S)]),$$
 where $G^3$ is the third power of $G$, that is, the graph on the same vertex set as $G$ in which two distinct vertices are adjacent whenever their distance in $G$ is at most three (the use of the third power will be explained shortly).

Finally, define
\[
        H_3:=
        \left\{
        e\in\binom{U}{k}:
        c(e)\le \ell
        \right\}.
\]
It is easy to see that $H_1\cup H_2\subseteq H_3$.

The ultimate $\ell$-hamiltonian saturated $k$-graph is defined via the standard ``sandwiching'' device used in our earlier papers on this topic:
take \(H\) as a maximal non-\(\ell\)-hamiltonian \(k\)-graph satisfying
\[
        H_1\cup H_2\subseteq H\subseteq  H_3.
\]
\noindent (Such an \(H\) exists by Lemma~\ref{main1} below.)

\subsection{Outline of proof of Theorem \ref{main_result}}

Theorem \ref{main_result} will be derived quickly from the following two lemmas
\begin{lem}\label{main1}
The \(k\)-graph \(H_1\cup H_2\) is not \(\ell\)-hamiltonian.
\end{lem}
Given a $k$-graph $H$ and a $k$-element set $e\subset V(H)$ such that $e\not\in H$, we denote by $H+e$ the $k$-graph $(V(H),H\cup\{e\})$.
\begin{lem}\label{main2}
For every
\[
        e\in \binom{U}{k}\setminus H_3,
\]
the \(k\)-graph
\[
        H_1\cup H_2+e
\]
contains a Hamilton \((\ell,k)\)-cycle.
\end{lem}

The use of the third power in the definition of $H_3$ can be now justified as follows. In order to construct a Hamilton $(\ell,k)$-cycle in $H_1\cup H_2+e$, $e\notin H_3$, we will seek a Hamilton path in $G_1$ between two non-adjacent vertices chosen from among the vertices of $\operatorname{tr}(e)$ or their neighbors. In the worst-case scenario, we have $\operatorname{tr}(e)\subseteq \{n+1,\dots,2n\}$, and we are thus compelled to find a Hamilton path between neighbors $u$ and $v$ of vertices $n+u$ and $n+v$, respectively.  Our construction ensures that $c(e)\ge\ell+1\ge2$ (with some margin), which enables us to select $n+u$ and $n+v$ from $\operatorname{tr}(e)$ so that they are non-adjacent in $G^3$. Consequently, by construction of $G$, $u$ and $v$ are non-adjacent in $G_1$, and, by the hamiltonian saturation of $G_1$, this non-adjacency guarantees the existence of a Hamilton path in $G_1$ with endpoints $u$ and $v$.

In bounding the size of $H$, we shall use the following elementary estimate whose proof is omitted (for a one line proof, see \cite[Proof of Fact 2.2]{RZ}).

\begin{claim}\label{claim:bounded-subsets} Fix positive integers $\Delta,r$ and $p$, and let \(F\) be a graph on \(O(n)\) vertices with maximum degree at most $\Delta$. Then the number of subsets
$ T\subseteq V(F)$
such that
$|T|\le r$ and $\operatorname{comp}(F[T])\le p$
is \(O_{\Delta,r,p}(n^p)\). \qed
\end{claim}

We are now ready for a short proof of our main result.

\begin{proof}[Proof of Theorem~\ref{main_result}]
By its definition, the hypergraph \(H\) is non-\(\ell\)-hamiltonian.
We first show the saturation of $H$. Let
\[
        e\in\binom{U}{k}\setminus H.
\]
If \(e\in H_3\), then since \(H\) is maximal non-\(\ell\)-hamiltonian among all \(k\)-graphs containing  \(H_1\cup H_2\) and
contained in \(H_3\), the hypergraph \(H+e\) is \(\ell\)-hamiltonian.

If \(e\notin H_3\), then, by Lemma~\ref{main2},
        $H_1\cup H_2+e$
is \(\ell\)-hamiltonian, and so is $H+e$.
Thus \(H\) is \(\ell\)-hamiltonian saturated.

 It remains to
estimate its number of edges. Since
        $H\subseteq H_3$,
it is enough to bound $|H_3|$.
 We may write
\[
        H_3=
        \bigcup_{\substack{
        T\subseteq V(G)\\
        |T|\le k,\ \operatorname{comp}(G^3[T])\le \ell}}
        \left\{
        e\in\binom{U}{k}:\operatorname{tr}(e)=T
        \right\}.
\]
Since \(G\) has bounded maximum degree, so does \(G^3\). By
Claim~\ref{claim:bounded-subsets}, the number of sets \(T\) in the above union is
        $O(n^\ell)$.
For a fixed \(T\), the number of \(k\)-sets \(e\) with \(\operatorname{tr}(e)=T\)
is bounded by a constant depending only on \(k,\ell\), because each \(U_i\) has
size \(O_{k,\ell}(1)\). Consequently,
\[
       |H|\le  |H_3|=O(n^\ell)
=O(N^\ell).
\]
This yields
$
        \operatorname{sat}(N,k,\ell)=O(N^\ell)
$
and, together with the lower bound, ultimately
$
        \operatorname{sat}(N,k,\ell)=\Theta(N^\ell).$

\end{proof}

The remainder of the paper is devoted to the proofs of Lemmas~\ref{main1} and~\ref{main2}.

\section{$H_1 \cup H_2$ is not hamiltonian }\label{notHam}

 In this section we prove Lemma \ref{main1}. The non-hamiltonicity of
\(H_1\cup H_2\) will be proved by showing that any potential Hamilton \((\ell,k)\)-cycle therein would
either force the presence of a Hamilton cycle in \(G_1\), a contradiction with the choice of $G_1$, or otherwise have fewer than \(N\) vertices, a contradiction with the notion of hamiltonicity.

 Suppose that $C=(e_1,\dots, e_{N/(k-\ell)})$ is a Hamilton $(\ell,k)$-cycle in $H_1 \cup H_2$. Recall that, on the one hand, $V(H_1)\subseteq A$, while on the other hand, by Corollary~\ref{boundH2},  every $H_2$-path has a bounded length. Therefore, $C$ cannot be entirely contained in either $H_1$ or $H_2$.  
 
  Let $\beta=(\beta_1,\dots,\beta_{N/(k-\ell)})$ be the binary (cyclic) sequence with $\beta_i=1$ if $e_i\in H_1$ and $\beta_i=2$ if $e_i\in H_2$.
We say that a (proper) sub-hypergraph of $C$ is an $H_i$-path in $C$, $i=1,2$, if the corresponding subsequence of $\beta$ forms an $i$-run, that is, a maximal subsequence of consecutive elements, all equal to $i$.

Let $m$ be the number of $H_1$-paths in $C$, which, clearly, equals the number of $H_2$-paths in $C$. Above we have argued that $m\ge1$, but, in fact,
since  
$$\sum_{i=1}^n|B_i|\ge n(\alpha-k-\ell)-h_N(k-\ell)\ge n(\alpha-2k)$$
 and each $H_2$-path has length $O_\alpha(1)$, we must have $m=\Theta(n)$. 

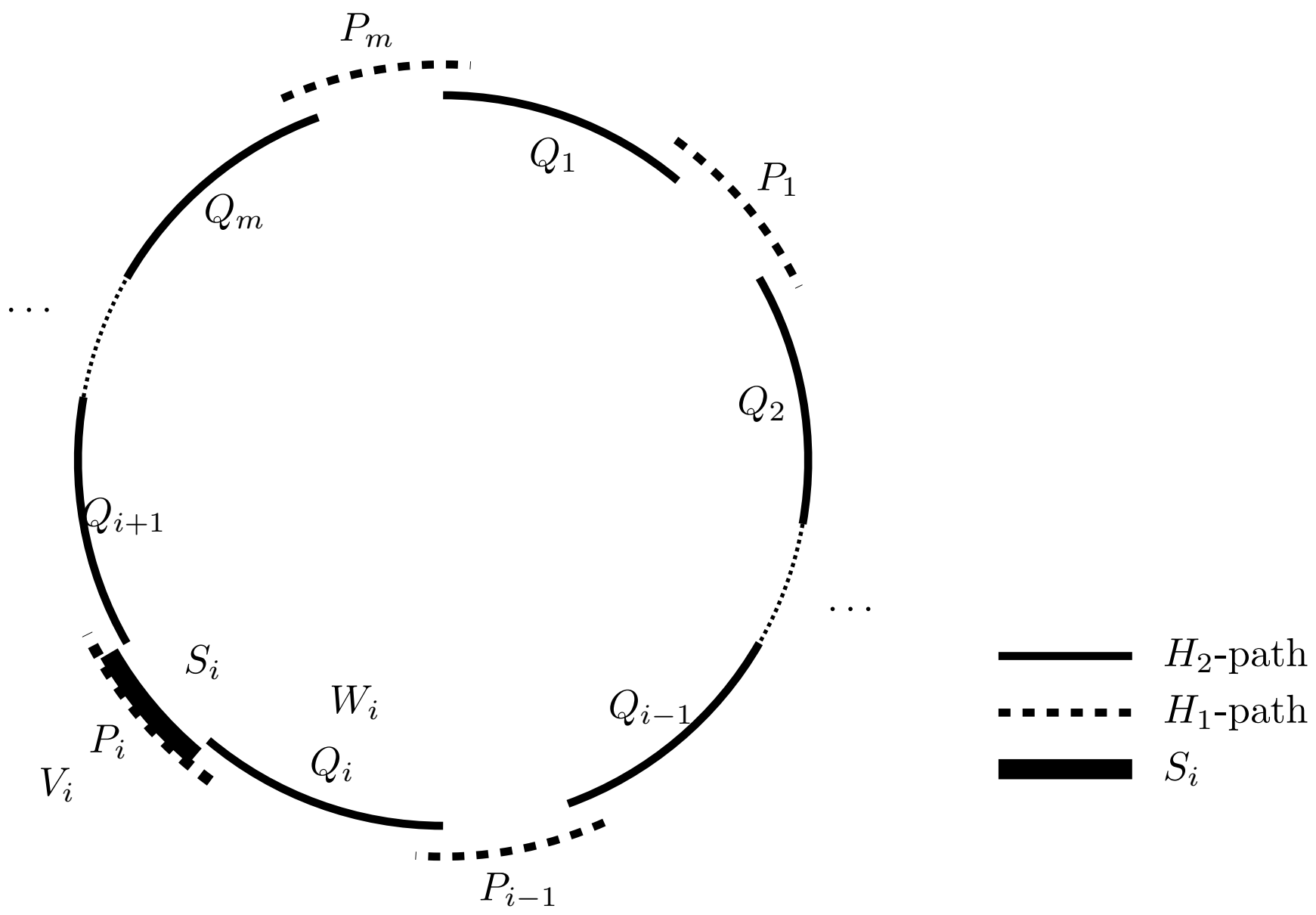
\begin{figure}[ht]
\centering
\begin{lrbox}{\reviewfigurebox}\begin{tikzpicture}[
  scale=0.9,every node/.style={font=\small}]
\def\RQ{3.55}
\def\RP{3.85}
\def\RLQ{3.15}
\def\RLP{4.25}
\def\Qlen{40}
\def\Overlap{4}
\def\A{90}
\def\B{30}
\def\C{-30}
\def\D{-90}
\def\E{-150}
\def\F{-210}

\foreach \ang/\lab in {90/{Q_1},30/{Q_2},-30/{Q_{i-1}},
                      -90/{Q_i},-150/{Q_{i+1}},-210/{Q_m}}{
  \draw[line width=2pt] (\ang:\RQ)
    arc[start angle=\ang,end angle={\ang-\Qlen},radius=\RQ];
  \node at ({\ang-0.5*\Qlen}:\RLQ) {$\lab$};
}

\draw[line width=2pt,dashed]
  ({\A-\Qlen+\Overlap}:\RP)
  arc[start angle={\A-\Qlen+\Overlap},end angle={\B-\Overlap},radius=\RP];
\node at ({((\A-\Qlen+\Overlap)+(\B-\Overlap))/2}:\RLP) {$P_1$};
\draw[line width=2pt,dashed]
  ({\F-\Qlen+\Overlap}:\RP)
  arc[start angle={\F-\Qlen+\Overlap},end angle={\A-360-\Overlap},radius=\RP];
\node at ({((\F-\Qlen+\Overlap)+(\A-360-\Overlap))/2}:\RLP)
  {$P_m$};

\draw[line width=1pt,densely dotted]
  ({\B-\Qlen}:\RQ)
  arc[start angle={\B-\Qlen},end angle=\C,radius=\RQ];
\node at ({((\B-\Qlen)+\C)/2}:\RLP) {$\cdots$};
\draw[line width=1pt,densely dotted]
  ({\E-\Qlen}:\RQ)
  arc[start angle={\E-\Qlen},end angle=\F,radius=\RQ];
\node at ({((\E-\Qlen)+\F)/2}:\RLP) {$\cdots$};

\draw[line width=2pt,dashed]
  ({\C-\Qlen+\Overlap}:\RP)
  arc[start angle={\C-\Qlen+\Overlap},end angle={\D-\Overlap},radius=\RP];
\node at ({((\C-\Qlen+\Overlap)+(\D-\Overlap))/2}:\RLP)
  {$P_{i-1}$};
\draw[line width=2.8pt,dashed]
  ({\D-\Qlen+\Overlap}:\RP)
  arc[start angle={\D-\Qlen+\Overlap},end angle={\E-\Overlap},radius=\RP];
\node at ({((\D-\Qlen+\Overlap)+(\E-\Overlap))/2}:\RLP)
  {$P_i$};

\draw[line width=5pt]
  ({\D-\Qlen}:{\RP-0.1})
  arc[start angle={\D-\Qlen},end angle=\E,radius={\RP-0.1}];
\node at ({((\D-\Qlen)+\E)/2}:3.05) {$S_i$};
\node at ({\D-0.5*\Qlen}:2.50) {$W_i$};
\node at ({((\D-\Qlen+\Overlap)+(\E-\Overlap))/2}:4.90)
  {$V_i$};

\begin{scope}[shift={(5.4,-1.9)}]
  \draw[line width=2pt] (0,0)--(1.3,0);
  \node[right] at (1.45,0) {$H_2$-path};
  \draw[line width=2pt,dashed] (0,-0.55)--(1.3,-0.55);
  \node[right] at (1.45,-0.55) {$H_1$-path};
  \draw[line width=5pt] (0,-1.1)--(1.3,-1.1);
  \node[right] at (1.45,-1.1) {$S_i$};
\end{scope}
\end{tikzpicture}\end{lrbox}
\ifdim\wd\reviewfigurebox>\linewidth\resizebox{\linewidth}{!}{\usebox{\reviewfigurebox}}\else\usebox{\reviewfigurebox}\fi
\caption{The $H_1$-paths $P_1,\dots,P_m$ and the $H_2$-paths $Q_1,\dots,Q_m$ alternating around the cycle $C$, where $P_i$ lies between $Q_i$ and $Q_{i+1}$.}
\label{fig:alternating_paths}
\end{figure}

Naturally, the $H_1$-paths and $H_2$-paths alternate around $C$. Let us denote the consecutive $H_1$-paths by $P_1,\dots,P_m$ and the $H_2$-paths by $Q_1,\dots,Q_m$, so that
$P_i$ is between $Q_i$ and $Q_{i+1}$, $i=1,\dots, m$, where $m+1:=1$ (see Figure~\ref{fig:alternating_paths}).
Further, let
$$V_i=V(P_i),\quad W_i=V(Q_i),\quad S_i=V_i\setminus(W_i \cup W_{i+1}),\quad  v_i=\min(W_i).$$ 

\subsection{Useful facts}\label{usfacts}

We will first prove some useful facts about the $H_1$-paths $P_i$ and the $H_2$-paths $Q_i$ and their mutual interference. Let us call the set of, respectively, the first and the last $\ell$ vertices of an $(\ell,k)$-path  the
\emph{initial} and \emph{terminal \(\ell\)-set}. (If the direction of  a path is irrelevant, we call both these $\ell$-sets \emph{terminal}.)
\begin{fact}\label{fakt1}
    For each $i\in[m]$, both terminal $\ell$-sets of $Q_i$ are subsets of~$A_{v_i}$.
\end{fact}
 \proof

Let $e$ be the last edge of $P_{i-1}$ and let $f$ be the first edge of $Q_i$. Since $P_{i-1}$ and $Q_i$ are consecutive $(\ell,k)$-paths on  $C$, we have \( |e\cap f|=\ell. \) Moreover, the set \( L:=e\cap f \) is precisely the initial block of $\ell$ vertices of $Q_i$. Since $e\in H_1$, we have $e\subseteq A$. Consequently, \[ L\in\binom{f\cap A}{\ell}. \] Furthermore, since $\operatorname{tr}(L)\subseteq \operatorname{tr}(e)$, by the definition of $\mathcal{L}$ one has \( L\in\mathcal{L}. \) In turn, since $f\in H_2\subseteq\mathcal K$, by the definition of $\mathcal K$, it follows that \( L\subseteq A_{\min(f)}. \) Since $f$ is an edge of $Q_i$, by~\eqref{minP} we have \( \min(f)=\min(Q_i)=v_i. \) Hence \[ e\cap f=L\subseteq A_{v_i}, \] that is, the initial block of $Q_i$ consists of $\ell$ vertices from $A_{v_i}$. The argument for the terminal block is analogous.
\qed

\begin{fact}\label{fakt2}
  For each $i\in[m]$,  $v_i\in\operatorname{tr} (V_{i-1})\cap\operatorname{tr} (V_{i})$. It follows that \newline $\operatorname{tr}(V_{i})\supseteq\{v_i,v_{i+1}\}$.
\end{fact}
\proof Follows trivially from Fact~\ref{fakt1}.  \qed

\medskip

Although the paths  $Q_i$ are, by definition, pairwise edge-disjoint, in principle, they do not need to be pairwise vertex-disjoint.  Below we show that, in fact, they are.

\begin{fact}\label{fakt3}
    The paths $Q_1,\dots, Q_m$ are pairwise vertex disjoint.
\end{fact}
\proof Suppose $W_i\cap W_{i+1}\neq\emptyset$, for some $i\in[m]$. Then,  
$|V_{i}|\le 2\ell-1$, because $V_{i}\subseteq W_i\cup W_{i+1}$, $|V_{i}\cap W_i|=\ell$ and $|V_{i}\cap W_{i+1}|=\ell$ . But, by Fact~\ref{fakt1}, $v_i=v_{i+1}$, so $V_{i}\subseteq A_{v_i}$. This means that $\operatorname{tr}(V_{i})=\{v_i\}$ is a singleton, a contradiction with the definition of $H_1$, where each edge has the trace of size two.
 \qed

\begin{fact}\label{fakt3_5}
For each $i\in[m]$,    $v_i\leq n$.
\end{fact}
\proof Let $L$ and $R$ be, respectively, the sets of the first and the last $\ell$ vertices of $Q_i$. If $L\cap R\neq\emptyset$, then $W_i\subseteq L\cup R$ and so $|L\cup R|\ge k$ (as $Q_i$ has at least one edge). If $L\cap R=\emptyset$, then, clearly, $|L\cup R|=2\ell\ge k$. By Fact~\ref{fakt1}, $L\cup R\subseteq A_{v_i}$. Thus, $|A_{v_i}|\ge k$ and, by~\eqref{mocai2026}, $v_i\le n$.\qed

\begin{fact}\label{fakt4}
    For all $i,j\in[m]$, $i \neq j$, we have $v_i\neq v_j$. It follows that $m\le n$.
\end{fact}
\proof By Fact~\ref{fakt3}, $W_i\cap W_j=\emptyset$.  By an argument similar to that in the proof of Fact~\ref{fakt3_5}, $|W_i \cap A_{v_i}|\ge k$ and $|W_j \cap A_{v_j}|\ge k$. If $v_i=v_j$, then

 $$|A_{v_i}|\ge |W_i \cap A_{v_i}|+|W_j \cap A_{v_i}| \geq 2k > k+\ell,$$ a contradiction to \eqref{mocai2026}. \qed

\medskip

For each $i=1,\dots,m$, set  $\tau_i=|\operatorname{tr}(V_i)|$ and $R_i=\left\{ j\le n:\; |S_i \cap A_j|\ge \ell \right\}$.

\begin{fact}\label{fact:11_new}
For each $i\in[m]$,
\[
\left| R_i \right| \ge \frac{\tau_i-2}{7}.
\]
In particular, if $\tau_i\ge 3$, then $R_i\neq\emptyset$.
\end{fact}
\proof
Let $P_{i}=e_1,e_2,\dots,e_r$,
where the edges are listed in the order from
$Q_{i}$ to $Q_{i+1}$, so that $|e_1\cap W_{i}|=|e_r\cap W_{i+1}|=\ell$. Recall that, by Fact~\ref{fakt1}, $e_1\cap W_{i}\subseteq A_{v_{i}}$ and $e_r\cap W_{i+1}\subseteq A_{v_{i+1}}$. In particular, $v_{i}\in\operatorname{tr} (e_1)$ and $v_{i+1}\in \operatorname{tr}(e_r)$.

 Since for every $h=1,\dots,r$, $e_h\in H_1$, we have $e_h\subset A$ and $|\operatorname{tr}(e_h)|=2$.
 Thus, the only possibility for $\operatorname{tr}(e_{h+1})$ to differ from $\operatorname{tr}(e_h)$ is when, for some $j$, $\operatorname{tr}(e_h)\cap \operatorname{tr}(e_{h+1})=\{j\} $, that is, $L_h:=e_h\cap e_{h+1}\subseteq A_j$. Call such a vertex $j\in V(G)$ an \emph{$h$-swapper} (in $P_i$), or just a \emph{swapper}. By the definition of $H_1$, for a swapper~$j$ we have $deg_G(j)\ge2$. It follows that no vertex $j>n$, as having only one neighbor in $G$, can be a swapper.

By Fact~\ref{fakt2} we know that $\operatorname{tr}(V_i)\supseteq\{v_{i},v_{i+1}\}$.
 If $\operatorname{tr}(V_i)=\{v_{i},v_{i+1}\}$, then there is nothing to prove.
 Otherwise, let $T_i = \{j: j \text{ is a swapper in } P_i \}$ and $t:=|T_i|$. We are going to relate $t$ to $\tau_i$, and then to $|R_i|$. Let
$x\in tr(V_i)\setminus\{v_i,v_{i+1}\}$ but $x\not\in T_i$.  Some change of trace must introduce or remove $x$: otherwise all traces would contain $x$, and, since $x$ is not a swapper, they would all be the same pair throughout, a contradiction with $v_{i}\in \operatorname{tr}(e_1)$ and $v_{i+1}\in tr(e_r)$.  At such a change the common vertex is a swapper $y$ adjacent to $x$. So, for every $x\not\in T_i$ there is a $y\in T_i$ with $xy\in E(G)$, which implies that the  non-swappers (other than $v_{i}$ and $v_{i+1}$) are exclusively among the neighbors of the swappers. Therefore,
\begin{equation}\label{taut}
        \tau_i-2
        \le t +6t=
        7t.
\end{equation}

On the other hand, if $j$ is an $h$-swapper, then  $L_h\subset A_j$. But, $L_h \subseteq S_i$ too. Indeed, if $L_h \cap W_{i} \neq \emptyset$ or $L_h \cap W_{i+1} \neq \emptyset$, then, respectively, $|\operatorname{tr}(e_h)|=1$ or $|\operatorname{tr}(e_{h+1})|=1$, each a contradiction with the definition of $H_1$.
Hence,  $|S_i\cap A_j|\ge|L_h|=\ell$ and it follows that $j\in R_i$. Consequently, $t\le |R_i|$ which, together with~\eqref{taut}, completes the proof.
\qed

\medskip

Set $W=\bigcup_{i=1}^m W_i$ and, for $j=1,\dots,2n$, 
$$r_j=|A_j\setminus W|.$$
Further, let
\[
J:=\{j \in [n]: r_j \ge \ell\},
\]
Below we prove two simple statements about the set $J$.

\begin{fact}\label{JJ} 
\begin{itemize}
\item[(i)] If for some $i\in[m]$\; $\tau_i\ge 3$, then $|J|\geq 1$.
\item[(ii)] For every $i\in[m]$, if  $v_i\in J$, then $|W_i|\leq k$.
\end{itemize}
\end{fact}
\proof (i) Let $i\in[m]$ be such that $\tau_i\ge 3$. Then, by Fact~\ref{fact:11_new},
$|S_{i}\cap A_{j}|\ge \ell$
for some $j\le n$. Since $A_j\setminus W\supseteq S_i\cap A_j$, it follows that $|A_j\setminus W|\ge\ell$,
which means that $j\in J$.

(ii) Suppose that  $v_i\in J$ and the $H_2$-path $Q_i$ has at least two edges. Then, by Fact~\ref{fakt1}, the first block of $\ell$ vertices of $Q_i$ ($L$) and the last block of $\ell$ vertices ($R$) both are subsets of $A_{v_i}$. Thus, if $L\cap R=\emptyset$, then $|L\cup R|=2\ell$. Otherwise, $W_i=L\cup R\subseteq A_{v_i}$ and $|W_i|\ge 2k-\ell$, where the bound is achieved when $Q_i$ has exactly two edges. In any case, $|L\cup R|\ge\min\{2\ell,2k-\ell\}$. But, on the other hand, by~\eqref{mocai2026},
$$|A_{v_i}\cap W_i|\le |A_{v_i}|-r_{v_i}\le(k+\ell)-\ell=k<\min\{2\ell,2k-\ell\},$$
a contradiction.

\qed




\subsection{Proof of Lemma~\ref{main1}}

We are now in a position to finish the proof of Lemma~\ref{main1}. Recall that, by means of contradiction, we have assumed that there is a Hamilton $(\ell,k)$-cycle $C$ which consists of $H_1$-paths $P_1,\dots,P_m$ and $H_2$-paths $Q_i$ appearing on $C$ alternately. We have also set up notation 
$$V_i=V(P_i),\quad W_i=V(Q_i),\quad S_i=V_i\setminus(W_i \cup W_{i+1}),\quad  v_i=\min(W_i)$$
 and proved, among others,  that the sets $W_i$ are pairwise disjoint (Fact~\ref{fakt3}), $v_i\le n$ (Fact~\ref{fakt3_5}), and $m\le n$ (Fact~\ref{fakt4}).
 Our ultimate goal is to get a contradiction by showing that $|V(C)|<N$.

 Recall that $W=\bigcup_{i=1}^m W_i$, and, for $j=1,\dots,2n$, set
$$A_j'=A_j\setminus W\;,\quad\mbox{and (recall)}\quad r_j=|A_j'|.$$
Further, set $S=\bigcup_{i=1}^m S_i$ and recall that $A=\bigcup_{j=1}^{2n}A_j$.

By the definition of $H_1$, we have $S\subseteq A$, and, since $S\cup W=U$ and $S\cap W=\emptyset$,
\begin{equation}\label{S=sum_r_j}
S=A\setminus W=\bigcup_{j=1}^{2n}A_j'\quad\mbox{and thus}\quad |S|= \sum_{j=1}^{2n} r_j.
\end{equation}


In bounding $|V(C)|$ we will be using the auxiliary function $\Phi(\alpha,m)$.
We consider two cases, but in each of them we will get the same intermediate bound
\begin{equation}\label{eq:common-nonham}
\begin{gathered}
|V(C)|\le\Phi(\alpha,L)+Dc_--(k-\ell)t_N,\\
\end{gathered}
\end{equation}
where $L$ and $t_N$ will be determined later and be different in each case, while $D:=n-L\ge\max\{ h_N-t_N,1\}$.

Using $L=n-D$ and Lemma~\ref{phi-step}, it follows from~\eqref{eq:common-nonham} that
\begin{equation}
\begin{aligned}\label{VCN}
|V(C)|
&\le\Phi(\alpha,n-1)+c_-
 -(D-1)(\nu(\alpha)-c_3-c_-)-(k-\ell)t_N\\
&\le\Phi(\alpha,n-1)+c_--(k-\ell)(D-1+t_N)\\
&\le\Phi(\alpha,n-1)+c_--(k-\ell)\max\{h_N-1,0\}<N.
\end{aligned}
\end{equation}
Here the second inequality  uses $D\ge1$ and
$\nu(\alpha)-c_3-c_-\ge k-\ell$, the latter following
from~\eqref{nu_i_alpha}(i).
The third one follows from
$D\ge h_N-t_N$, and the last one from~\eqref{nn}.
This contradicts the assumption that
$C$ is a Hamilton $(\ell,k)$-cycle on $N$ vertices.
Thus, it remains to establish~\eqref{eq:common-nonham}. Recall that, by Fact~\ref{fakt2}, we have $\tau_i=|\operatorname{tr}(V_i)|\ge2$, for all $i=1,\dots,m$.

\medskip

{\bf Case 1.} Assume first that $\tau_i=2$ and $|W_i| \geq 2\ell$ for every $i=1,\dots,m$.
Then, by Facts~\ref{fakt2} and~\ref{fakt4}, the set of pairs
$\{\operatorname{tr}(V_i): i=1,\dots,m\}$ forms a connected $2$-regular subgraph of $G_1$, that is, a cycle, on the vertex set $\{v_1,\dots,v_m\}$ (in this order).
Hence, by the non-hamiltonicity of $G_1$,
$m\le n-1$. 
Put $L=m$ and $t_N=|[h_N]\cap\{v_1,\dots,v_m\}|$. Then, quite trivially, $D=n-m\ge \max\{h_N-t_N,1\}$.

By the definition of $H_2$ and by  Facts~\ref{minH2} and~\ref{fakt1}, each path $Q_i$ satisfies assumptions (ii) and (iii) of Definition~\ref{dom} with respect to $U_{v_i}$. Thus, $|U_{v_i}\cap W_i|\ge2\ell$ and, by the remark after this definition, $Q_i$ is $(U_{v_i}\cap W_{i})$-dominant. Consequently, 
by Definition~\ref{defniu} and the monotonicity of $\nu$ (see Lemma~\ref{nu1}(c)), for every $i\in[m]$,
\begin{equation}\label{Wunu}
|V(Q_i)|=|W_i|\le \nu(|U_{v_i}\cap W_{i}|)\le \nu(\alpha_{v_i}-r_{v_i}),
\end{equation}
because $U_{v_i}\cap W_{i}\subseteq U_{v_i}\setminus A'_{v_i}$.

       \smallskip

In turn, by~\eqref{S=sum_r_j}, $S=\bigcup_{i=1}^m A_{v_i}'$; in fact, since $\tau_i=2$,  $S_i\subseteq A_{v_{i}}'\cup A_{v_{i+1}}'$, $i=1,\dots m-1$,
and $S_m\subseteq A_{v_m}'\cup A_{v_1}'$.

\medskip

 It follows that 
 \begin{equation}\label{S3}
 |S_i|\in\{g,g+k-\ell,g+2(k-\ell)\}.
 \end{equation}
 
 Indeed, since \(P_i\) is an \((\ell,k)\)-path, by  Definition \ref{defgap},
\[
 |S_i|\equiv g\pmod {k-\ell}.
\]
On the other hand, Fact~\ref{fakt1} and \(|W_j|\ge2\ell\) imply that
\begin{align}\label{a_priv_i}
 |A'_{v_j}|\le k+\ell-2\ell=k-\ell, \quad j\in[m].
\end{align}
Thus, \(|S_i|\le2(k-\ell)\), and we get~\eqref{S3}
(the third option only when $g=0$.)

\medskip

\begin{figure}[htbp]
\centering
\begin{lrbox}{\reviewfigurebox}\begin{tikzpicture}[
  x=0.48cm,y=0.70cm,
  every node/.style={font=\small},
  brace/.style={decorate,decoration={brace,amplitude=4pt}}
]
  \foreach \a/\b/\lab in {
    0/3/{Q_1},
    4.5/7.7/{Q_{i-1}},
    9.1/13.7/{Q_i},
    15.1/18.3/{Q_{i+1}},
    19.8/22.8/{Q_m}}{
    \fill[black!6] (\a,0) rectangle ({\a+0.7},0.8);
    \fill[black!6] ({\b-0.7},0) rectangle (\b,0.8);
    \draw (\a,0) rectangle (\b,0.8);
    \draw ({\a+0.7},0)--({\a+0.7},0.8);
    \draw ({\b-0.7},0)--({\b-0.7},0.8);
    \node at ({\a+0.35},0.4) {$\ell$};
    \node at ({(\a+\b)/2},0.4) {$\cdots$};
    \node at ({\b-0.35},0.4) {$\ell$};
    \draw[brace] (\a,1.10)--(\b,1.10)
      node[midway,above=6pt] {$\lab$};
  }
  \draw[thick] (9.1,0) rectangle (13.7,0.8);

  \draw (3,0.4)--(4.5,0.4);
  \node[fill=white,inner sep=2pt] at (3.75,0.4) {$\cdots$};
  \draw (18.3,0.4)--(19.8,0.4);
  \node[fill=white,inner sep=2pt] at (19.05,0.4) {$\cdots$};

  \fill[black!35] (8.35,0) rectangle (9.1,0.8);
  \draw (7.7,0) rectangle (9.1,0.8);
  \draw (8.35,0)--(8.35,0.8);
  \node[font=\scriptsize] at (8.025,0.4) {$\cdots$};
  \draw[brace] (7.7,1.10)--(9.1,1.10)
    node[midway,above=6pt] {$S_{i-1}$};

  \fill[black!35] (13.7,0) rectangle (14.45,0.8);
  \draw (13.7,0) rectangle (15.1,0.8);
  \draw (14.45,0)--(14.45,0.8);
  \node[font=\scriptsize] at (14.775,0.4) {$\cdots$};
  \draw[brace] (13.7,1.10)--(15.1,1.10)
    node[midway,above=6pt] {$S_i$};

  \draw[brace] (7,2.10)--(9.8,2.10)
    node[midway,above=6pt] {$P_{i-1}$};
  \draw[brace] (13,2.10)--(15.8,2.10)
    node[midway,above=6pt] {$P_i$};

  \draw[brace] (9.1,-0.28)--(8.35,-0.28)
    node[midway,below=8pt] {$A'_{v_i}\cap S_{i-1}$};
  \draw[brace] (14.45,-0.28)--(13.7,-0.28)
    node[midway,below=8pt] {$A'_{v_i}\cap S_i$};

\end{tikzpicture}\end{lrbox}
\ifdim\wd\reviewfigurebox>\linewidth\resizebox{\linewidth}{!}{\usebox{\reviewfigurebox}}\else\usebox{\reviewfigurebox}\fi
\caption{Case~1, with the closing part $S_m$ omitted.}
\label{fig:admissible-sequence}
\end{figure}
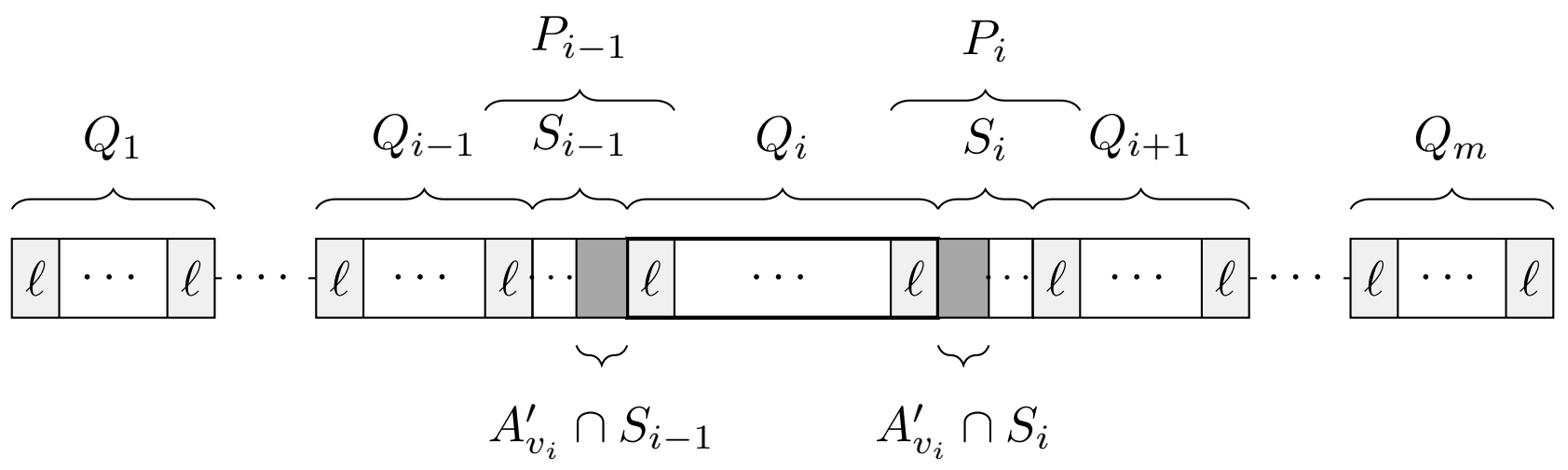

Note that, by \eqref{a_priv_i}, $r_{v_i}\le k-\ell$, $i=1,\dots,m$.
Let
\[
r^-_{i} = |A'_{v_i}\setminus S_m| \leq r_{v_i},\;i=1,\dots,m,\quad\mbox{and}\quad s_i=|S_{i}|,\; i=1,\dots,m-1.
\]
Then, it is easy to see that the sequence \(r^-_{1},\ldots,r^-_{m}\) is admissible in the sense of Definition~\ref{defadm}, with $s_1,\dots,s_{m-1}$ as a witnessing sequence.
Indeed, for every $1\le t\le m-1$, the set
$\bigcup_{i=1}^{t}(A'_{v_i}\setminus S_m)$ consists precisely
of all vertices of $S_1,\ldots,S_{t-1}$ and those vertices
of $S_{t}$ that belong to $A'_{v_t}$. Hence
\[
\sum_{i=1}^{t}r^-_i
=\sum_{j=1}^{t-1}s_j+|S_{t}\cap A'_{v_t}|.
\]
Since $0\le |S_{t}\cap A'_{v_t}|\le |S_{t}|=s_t$,
this yields both inequalities in~\eqref{adm}.
Moreover, the sets $A'_{v_i}\setminus S_m$, $i=1,\ldots,m$,
partition $\bigcup_{i=1}^{m-1}S_i$, so
\[
\sum_{i=1}^{m}r^-_i
=\sum_{i=1}^{m-1}|S_i|
=\sum_{j=1}^{m-1}s_j.
\]
Thus, the equality in~\eqref{adm} also holds 
(see Figure \ref{fig:admissible-sequence}).

This allows us to use function \(\Phi(\alpha,m)\) below. 

Finally, we bound the number of vertices in $C$ as follows.
\begin{align*}
        |V(C)|&=|W|+|S|=\sum_{i=1}^m |W_i|+\sum_{i=1}^m r^-_{i}+|S_m|\\
        &\overset{\eqref{Wunu}}{\le} \sum_{i=1}^m\bigl(\nu(\alpha_{v_i}-r_{v_i})+r^-_{i}\bigr)+|S_m|\\
        &\overset{Lem.\ref{nu1}(d)}{\le} \sum_{i=1}^m\bigl(\nu(\alpha-r_{v_i})+r^-_{i}\bigr)+|S_m|-(k-\ell)t_N\\
        &\overset{Lem.\ref{nu1}(c)}{\le} \sum_{i=1}^m\bigl(\nu(\alpha-r^-_{i})+r^-_{i}\bigr)+|S_m|-(k-\ell)t_N\\
        &\overset{Def.\ref{defphi}}{\le}
        \Phi(\alpha,m)+|S_m|-(k-\ell)t_N,
\end{align*}
where the second inequality uses Lemma~\ref{nu1}(d)  exactly $t_N$ times. Since $|S_m|\le2(k-\ell)<c_-$ and $D=n-m\ge1$, this implies~\eqref{eq:common-nonham} with $L=m$.

\medskip

{\bf Case 2.} Suppose now that either $\tau_i\ge 3$ or $|W_i|<2\ell$, for at least one $i\in[m]$.  Recall  that, by Fact~\ref{fakt4}, $m\le n$.

 For the forthcoming argument recall also that for each $i=1,\dots,m$ we put $ R_i:=\{j\le n:\ |S_i\cap A_j|\ge \ell\}$, and that we defined $J:=\{j \in [n]: r_j \ge \ell\}$, where $r_j=|A_j'|$.
  Trivially,  $R_i\subseteq J$ for each $i$.

\medskip

 Call an
$H_2$-path $Q_i$ \emph{short} if $|W_i| < 2\ell$. By Fact~\ref{JJ}(ii), since $k<2\ell$,  if
$v_i\in J$, then $Q_i$ is short. Let $s$ denote the number of short $H_2$-paths. It follows that $|J\cap\{v_1,\dots,v_m\}|\le s$.   Note also that by Fact~\ref{JJ}(i) and the assumptions of Case 2 here, if $J=\emptyset$, then $s\ge1$. Since $J\subseteq[n]$, we infer that
\begin{align}\label{D_ge1}
\max\{1,|J|\}\le s+n-m.
\end{align}

In turn, an $H_1$-path $P_i$ is called \emph{next-to-short} if it appears on $C$ next to a short $H_2$-path (either before, or after).
Clearly, the number of next-to-short $P_i$'s is at most $2s$.

Finally, call an $H_1$-path $P_i$ \emph{flawed}  if $\tau_i\ge 3$. Flawed $P_i$'s are related to the set $J$ as follows.

\begin{figure}[htbp]
\centering
\begin{tikzpicture}[
  x=1cm,y=1cm,
  every node/.style={font=\small},
  brace/.style={decorate,decoration={brace,amplitude=4pt}}
]
  \fill[black!20] (0,0) rectangle (3,0.8);
  \fill[black!10] (3,0) rectangle (6,0.8);
  \draw (0,0) rectangle (8.8,0.8);
  \draw (3,0)--(3,0.8) (6,0)--(6,0.8);
  \node at (1.5,0.4) {$S_{i_1}\cap A_j$};
  \node at (4.5,0.4) {$S_{i_2}\cap A_j$};
  \node at (7.4,0.4) {remainder};

  \draw[brace] (0,1.02)--(8.8,1.02)
    node[midway,above=6pt] {$|A_j|=k+\ell$};
  \draw[brace] (3,-0.17)--(0,-0.17)
    node[midway,below=6pt] {$\ge\ell$ vertices};
  \draw[brace] (6,-0.17)--(3,-0.17)
    node[midway,below=6pt] {$\ge\ell$ vertices};
  \draw[brace] (8.8,-0.17)--(6,-0.17)
    node[midway,below=6pt] {$\le k-\ell<\ell$ vertices};
\end{tikzpicture}
\caption{For distinct indices $i_1,i_2$ with $j\in R_{i_1}\cap R_{i_2}$,
the two displayed subsets of $A_j$ are disjoint. Any other subset $S_i\cap A_j$
lies in the remainder part, so $j\notin R_i$ for $i\notin\{i_1,i_2\}$.}
\label{fig:dirty-incidences}
\end{figure}
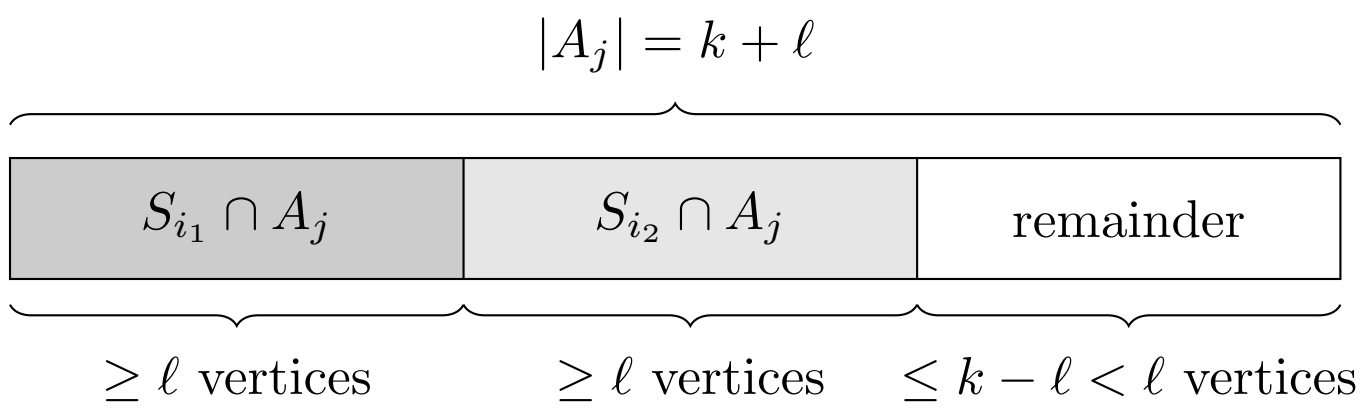

For every $i$ with $\tau_i\ge3$, by Fact~\ref{fact:11_new}, $|R_i|\ge1$.
On the other hand, for a fixed $j\in J$, the sets $S_i\cap A_j$ are pairwise
disjoint as $i$ varies, and of size at least $\ell$ whenever $j\in R_i$ (see Figure~\ref{fig:dirty-incidences}). Thus, since $|A_j|=k+\ell$ and $\ell>k/2$, a fixed $j$ belongs to at
most two sets $R_i$. Therefore,
\begin{equation}\label{dirt}
\sum_{P_i\text{ flawed}}|R_i|\le 2|J|.
\end{equation}
In particular, the number of flawed paths is itself bounded by $2|J|$.

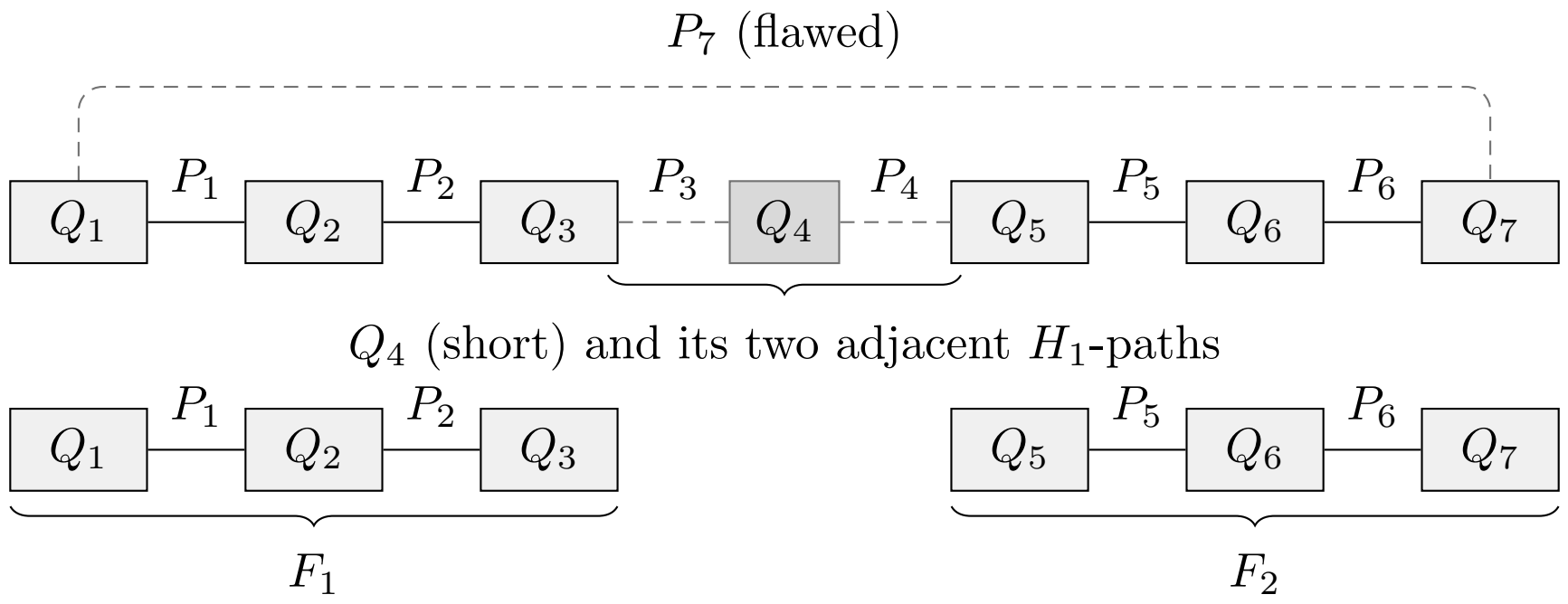
\begin{figure}[htbp]
\centering
\begin{tikzpicture}[
  x=0.86cm,y=0.90cm,
  every node/.style={font=\small},
  qpath/.style={draw,fill=black!6,minimum width=10mm,minimum height=6mm},
  removed/.style={black!55,densely dashed},
  brace/.style={decorate,decoration={brace,amplitude=4pt}}
]
  \draw[removed,rounded corners=5pt] (0,0.34)--(0,1.10)--(12,1.10)--(12,0.34);
  \node[above=2pt] at (6,1.10) {$P_7$ (flawed)};

  \foreach \i/\x in {1/0,2/2,3/4,5/8,6/10,7/12}
    \node[qpath] (q\i) at (\x,0) {$Q_{\i}$};
  \node[qpath,draw=black!55,fill=black!15,minimum width=8mm] (q4) at (6,0) {$Q_4$};
  \foreach \i/\j in {1/2,2/3,5/6,6/7}
    \draw (q\i.east)--(q\j.west) node[midway,above=1pt] {$P_{\i}$};
  \draw[removed] (q3.east)--(q4.west) node[midway,above=1pt,text=black] {$P_3$};
  \draw[removed] (q4.east)--(q5.west) node[midway,above=1pt,text=black] {$P_4$};
  \draw[brace] (7.50,-0.43)--(4.50,-0.43)
    node[midway,below=6pt] {$Q_4$ (short) and its two adjacent $H_1$-paths};

  \foreach \i/\x in {1/0,2/2,3/4,5/8,6/10,7/12}
    \node[qpath] (f\i) at (\x,-1.85) {$Q_{\i}$};
  \foreach \i/\j in {1/2,2/3,5/6,6/7}
    \draw (f\i.east)--(f\j.west) node[midway,above=1pt] {$P_{\i}$};
  \draw[brace] (4.58,-2.31)--(-0.58,-2.31)
    node[midway,below=6pt] {$F_1$};
  \draw[brace] (12.58,-2.31)--(7.42,-2.31)
    node[midway,below=6pt] {$F_2$};
\end{tikzpicture}
\caption{Schematic cyclic order of the paths in $C$ (top), with vertex overlaps
suppressed. Removing the flawed path $P_7$, the short path $Q_4$, and its two
neighbors $P_3,P_4$ leaves the clean strips $F_1$ and $F_2$ (bottom).}
\label{fig:clean-strips}
\end{figure}

Remove (virtually) all flawed, short, and next-to-short paths from $C$. What remains is a set of vertex disjoint $(\ell,k)$-paths, called \emph{clean strips}, all beginning and ending with an $H_2$-path (see Figure~\ref{fig:clean-strips}).
Let $q$ be the number of clean strips. Then, by the above considerations,
$$0\le q\le 2|J|+s.$$

We first bound the total number of vertices in all removed paths. Trivially, short $Q_i$'s contribute
at most $2\ell s$. The next-to-short $P_i$'s which are not flawed have each a trace of size two, so their sets
$S_i$ satisfy $|S_i|\leq2(k+\ell)$.
Finally, again by Fact~\ref{fact:11_new}, for every $i$ with $\tau_i\ge3$ we have $|R_i|\ge (\tau_i-2)/ 7$,
and so,
$$|S_i|\le (k+\ell)\tau_i\le7(k+\ell)|R_i|+2(k+\ell)\le9(k+\ell)|R_i|\le18k|R_i|.$$
Summing this up and using~\eqref{dirt}, 
\[
        \sum_{P_i\text{ flawed}} |S_i|
        \le
        18k\sum_{P_i\text{ flawed}} |R_i|
        \le
        36k|J|\le 36k(s+n-m).
\]

Thus, the whole  contribution of the removed paths to $V(C)$ is less than
\begin{equation}\label{whole_dirty}
        2\ell s+2s\cdot 2(k+\ell)+36k(s+n-m)
        \le 46k(s+n-m).
\end{equation}

We now turn to estimating the contribution of the $q$ clean strips $F_1,\dots, F_q$.
For each $t=1,\dots,q$,  the \emph{weight} $a_t$ of $F_t$ is defined as the number of $H_2$-paths contained in it. Set  $L:=a_1+\dots a_q$ and note that $L= m-s$.

 Fix $t$, set $F:=F_t$ and $a:=a_t$, and  let  $Q_{j+1},\dots, Q_{j+a}$ be the $H_2$-paths contained in $F$ and  $P_{j+1},\dots,P_{j+a-1}$ the $H_1$-paths between them, where the indices are read cyclically in $[m]$.
 For the ease of notation assume that $j=0$. Recall our notation  $v_{i}:=\min(Q_{i})$, $i=1,\dots,a$.
Recall also that, since the paths $Q_{i}$ are not short, while the paths $P_{i}$ are not flawed, we have $|W_{i}|\ge 2\ell$ for all $i=1,\dots a$, and  
 $tr(S_{i})\subseteq\{v_{i},v_{i+1}\}$ for all $i=1,\dots,a-1$.

Let $T=\bigcup_{i=1}^{a-1}S_{i}$ and put
\[
\rho_i=|A_{v_{i}}\cap T|,\qquad i=1,\ldots,a.
\]
Then \(0\le\rho_i\le r_{v_{i}}\le k-\ell\) (cf. \eqref{a_priv_i} and recall that, since $F$ is clean, we have $|W_{i}|\ge 2\ell$).

For \(i=1,\ldots,a-1\), again put \(s_i=|S_{i}|\). We have
\[
S_{i}\subseteq A_{v_{i}}\cup A_{v_{i+1}},
\]
and so for every \(h=1,\ldots,a-1\) we have
\[
\sum_{i=1}^{h-1}s_i
   \le \sum_{i=1}^h\rho_i
   \le \sum_{i=1}^{h}s_i,
\qquad
\sum_{i=1}^{a}\rho_i
   =\sum_{i=1}^{a-1}s_i.
\]
As in Case~1, \(s_i\in\{g,g+k-\ell,g+2(k-\ell)\}\), and hence
\((\rho_1,\ldots,\rho_{a})\in\mathcal A_{a}\).
If \(a=1\), then \(T=\emptyset\) and \(\rho_1=0\), so the same conclusion follows from
\(\mathcal A_1=\{(0)\}\).

 By the definition of \(\nu\) and monotonicity (see Lemma~\ref{nu1}(c)), since  \(\rho_i\le r_{v_{i}}\), we have
\[
|W_{i}|\le \nu(\alpha_{v_{i}}-r_{v_{i}})\le
\nu(\alpha_{v_{i}}-\rho_i).
\]

Consequently,
\[
\begin{aligned}
|V(F)|
&=\sum_{i=1}^{a}\bigl(|W_{i}|+\rho_i\bigr)
\le
\sum_{i=1}^{a}
\bigl(\nu(\alpha_{v_{i}}-\rho_i)+\rho_i\bigr).
\end{aligned}
\]
For \(v_{i}>h_N\) we have \(\alpha_{v_{i}}=\alpha\), whereas for \(v_{i}\le h_N\),
\(\alpha_{v_{i}}=\alpha-k+\ell\). In the latter case, by Lemma~\ref{nu1}(d) we bound
\[
\nu(\alpha_{v_{i}}-\rho_i)
\le
\nu(\alpha-\rho_i)-(k-\ell).
\]
Since \((\rho_1,\ldots,\rho_a)\in\mathcal A_{a}\), the definition of \(\Phi\) now yields
\[
|V(F)|
\le
\Phi(\alpha,a)
-(k-\ell)\cdot|[h_N]\cap\{v_{1},\ldots,v_{a}\}|.
\]
By \eqref{D_ge1},
\[
D:=s+n-m=n-L\ge1.
\]
Let \(t_N\) be the total number of minima $v_i\le h_N$ that occur in the union of \emph{all} clean strips. Altogether, clean strips contain \(L=m-s\) \(H_2\)-paths with distinct minima in \([n]\). Hence, at most \(n-L=D\) elements of \([h_N]\) are missing, and therefore
\[
D+t_N\ge h_N.
\]
Summing the above estimate over the \(q\) clean strips $F_t$ and applying Lemma~\ref{phi-glue}, we obtain
\[
\sum_{t=1}^q |V(F_t)|
<
\Phi(\alpha,L)+c_4q-(k-\ell)t_N.
\]
 Recall that the removed paths contribute at most \(46kD\) vertices and
\[
q\le2|J|+s\le3D,
\]
cf. \eqref{D_ge1}.
It follows, recalling~\eqref{whole_dirty}, that
\begin{equation*}
|V(C)|
<
46kD+\Phi(\alpha,L)+3c_4D-(k-\ell)t_N=\Phi(\alpha,L)+Dc_--(k-\ell)t_N.
\end{equation*}
 (If no clean strip remains, then \(L=t_N=0\), and the bound on the removed paths gives the same inequality, with \(\Phi(\alpha,0)=0\).)
 Together with \(D+t_N\ge h_N\), this is~\eqref{eq:common-nonham}, completing Case~2 and thus the whole proof of Lemma~\ref{main1}.

\section{$(H_1\cup H_2)+e$ is hamiltonian}\label{eHam}
In this section we prove Lemma~\ref{main2}, that is, for every $e\in\binom{U}{k}\setminus H_3$ we build a Hamilton $(\ell,k)$-cycle $C$ in $(H_1\cup H_2)+e$. Naturally, it will be seeded at $e$, more specifically, at two carefully selected vertices $v_1$ and $v_n$ of $G_1$. Their choice will be made in Section~\ref{v1vn}.

While  constructing $C$, we will profit from our experience gained when analyzing the  structure of a ``presumed-to-exist'' Hamilton $(\ell,k)$-cycle in the proof of Lemma~\ref{main1}. Thus, $C$ will consist of $H_2$-paths $Q_1,\dots,Q_n$ and 
connecting paths $P_1,\dots, P_{n-1}$ in $H_1 \cup H_2$ alternating around it, but also of a special path $\Pi$, playing the role of $P_n$, and containing $e$ as well as some edges from both, $H_1$ and $H_2$. To emphasize its hybrid nature we denote it by $\Pi$. Path $\Pi$ will be constructed in Section~\ref{Qzero}.

The main difficulty of this proof lies, however, in guaranteeing that all vertices  have been used. To this end, we need first to `tune up'' all the parameters involved in the construction. This is done in Section~\ref{numpre} which ends with a core lemma, Lemma~\ref{88}, the proof of which is postponed to Section~\ref{sec:correct-loads}. Then the final construction comes in two stages: the preliminary one where, using some ``fake'' vertices, one builds a cycle $C'$ with more than $N$ vertices (Section~\ref{toolong}), and the ultimate one, when the excess fake vertices are trimmed out, while the others are replaced by still missing elements of the vertex set $U$ (Section~\ref{getham}).

In this section, we will often use the notion of a terminal $\ell$-set of an $(\ell,k)$-path defined in Section~\ref{usfacts}. For a sequence of vertices, by a \emph{$k$-window} we mean any block of $k$ consecutive vertices in that sequence.
A $k$-window is \emph{relevant} if its initial position is congruent to 1 (mod $k-\ell$), so that it corresponds to a potential edge of an $(\ell,k)$-path or cycle.
\subsection{The choice of $v_1$ and $v_n$}\label{v1vn}
Let $e\in\binom{U}{k}\setminus H_3$.
Then, by the definition of $H_3$, we have \(c(e)\ge \ell+1\) where, recall,  $c(e)=\operatorname{comp}(G^3[\operatorname{tr}(e)])$. For $Z\subseteq\operatorname{tr}(e)$, write
\[
    e(Z):=\{u\in e:u\in\bigcup_{i\in Z}U_i\}=\{u\in e:\; \operatorname{tr}(\{u\})\subseteq  Z\}.
\]

In other words, $e(Z)$ is the largest subset $e'\subseteq e$ with $\operatorname{tr}(e')\subseteq Z$.
We shall later use the following observations:
\begin{equation}\label{eq:subset-bound-main2}
    |e'|\le k - \ell-1+|\operatorname{tr}(e')|\
    \qquad\text{for every}\quad e'\subseteq e
\end{equation}
and
\begin{equation}\label{eq:component-bound-main2}
    |e(Z)|\le k - \ell
    \qquad\text{for every connected component \(Z\) of }
    G^3[\operatorname{tr}(e)].
\end{equation}

 To prove them, set $e'':=e(\operatorname{tr}(e'))$ and note that $\operatorname{tr}(e'')=\operatorname{tr}(e')$, $|e''|\ge |e'|$, \newline and $\operatorname{tr}(e'')\cap\operatorname{tr}(e\setminus e'')=\emptyset$. Then, $c(e)\le c(e'')+c(e\setminus e'')$ and, since $c(e)\ge \ell+1$, we get
$$ k-|e''|\ge c(e\setminus e'')\ge \ell+1-c(e''),$$
equivalently,
\begin{equation}\label{ebis}
|e''|\le k-\ell-1+c(e'').
\end{equation}

Inequality~\eqref{eq:subset-bound-main2} follows, because $|e'|\le|e''|$ and $c(e'')\le |\operatorname{tr}(e'')|=|\operatorname{tr}(e')|$.
To deduce \eqref{eq:component-bound-main2}, set $e'=e(Z)$  and notice that  $e''=e'$ and  $c(e')=1$. Then,~\eqref{ebis} becomes~\eqref{eq:component-bound-main2}.

\smallskip

In particular, since $k-\ell<\ell$, it follows from~\eqref{eq:component-bound-main2}   that
\begin{equation}\label{eq:no-L1-in-e}
    \binom{e\cap A}{\ell}\cap\mathcal L=\varnothing,
\end{equation}
where, recalling  from Section~\ref{main_con},
\[ \mathcal{L} = \left\{ L\in\binom{A}{\ell}: \operatorname{tr}(L)\subseteq \{x,y\} \text{ for some }xy\in E(G) \right\}. \]

Let \(X_0\) be the component of \(G^3[\operatorname{tr}(e)]\) containing
\[
    x:=\min(e),
\]
put \(Y_0:=\operatorname{tr}(e)\setminus X_0\), and set
\[
    y:=\min(e(Y_0)).
\]
Observe that $X_0\cup Y_0=\operatorname{tr}(e)$,  $e=e(X_0)\cup e(Y_0)$, and, by~\eqref{eq:component-bound-main2},
\[
    |e(X_0)|\le k-\ell,
    \qquad\mbox{so}\qquad
    |e(Y_0)|\ge\ell.
\]
For \(z\in\{x,y\}\), define
\[
    z':=
    \begin{cases}
        z,   & z\le n,\\
        z-n, & z>n.
    \end{cases}
\]
Since \(x\) and \(y\) belong to different components of
\(G^3[\operatorname{tr}(e)]\), the vertices \(x'\) and \(y'\)
are non-adjacent in \(G_1\). Hence, by the hamiltonian saturation  of
\(G_1\), there is a Hamilton path
    $v_1v_2\cdots v_n$
in \(G_1\) with $ v_1=x',$
    and $v_n=y'$. (Note that $\{v_1,\dots,v_n\}=[n]$.)

\subsection{Constructing path $\Pi$}\label{Qzero}
We first construct a bounded length path \(\Pi\) containing \(e\), whose left
terminal \(\ell\)-set belongs to \(U_{v_n}\) and whose right terminal
\(\ell\)-set belongs to \(U_{v_1}\). Throughout this section,  we adopt a convention that, given a set of vertices  $S$ and an integer $r\le|S|$,  symbol \(S^r\) represents \emph{any} sequence $s_1,\dots,s_r$ of  \emph{distinct} vertices of \(S\). 
Also, due to the shortage of letters in the alphabet, we decided to locally reuse some letters again. So, the $t$ and $h$ appearing here have nothing to do with their counterparts in earlier sections.

Order the vertices of \(e\) so that all vertices of \(e(Y_0)\) precede
all vertices of \(e(X_0)\), the first vertex belongs to \(U_y\), and the
last vertex belongs to \(U_x\).

Two cases with respect to $y$ are considered. Recall the gap $g$ defined in Section~\ref{secPhi} (see  Definition~\ref{defgap}) and note that
\begin{align}\label{rel_e}
    \ell+g=p(k-\ell)
\end{align}
for some integer \(p\ge2\).

\begin{figure}[htbp]
\centering
\begin{lrbox}{\reviewfigurebox}
\begin{tikzpicture}[x=0.82cm,y=0.85cm,
 every node/.style={font=\small},
 brace/.style={decorate,decoration={brace,amplitude=4pt}}]
\fill[black!6] (0,-0.35) rectangle (2.5,0.35);
\draw (0,-0.35) rectangle (2.5,0.35);
\node at (1.25,0) {$B_y^{\ell+g}$};
\fill[white] (2.5,-0.35) rectangle (5.1,0.35);
\draw (2.5,-0.35) rectangle (5.1,0.35);
\node at (3.8,0) {$e(Y_0)$};
\fill[white] (5.1,-0.35) rectangle (7.2,0.35);
\draw (5.1,-0.35) rectangle (7.2,0.35);
\node at (6.15,0) {$e(X_0)$};
\fill[black!6] (7.2,-0.35) rectangle (9.7,0.35);
\draw (7.2,-0.35) rectangle (9.7,0.35);
\node at (8.45,0) {$B_x^{\ell+g}$};

\foreach \x/\u in {2.63/y,7.07/x}{
  \fill (\x,0) circle[radius=1.5pt];
  \draw[thin] (\x,-.10) -- (\x,-.58);
  \node[below,inner sep=2pt] at (\x,-.58) {$u_{\u}$};
}

\draw[brace] (2.5,0.65)--(7.2,0.65) node[midway,above=6pt] {$e$};
\end{tikzpicture}
\end{lrbox}
\ifdim\wd\reviewfigurebox>\linewidth\resizebox{\linewidth}{!}{\usebox{\reviewfigurebox}}\else\usebox{\reviewfigurebox}\fi
\caption{Description of $\Pi$, case I.}
\label{eee}
\end{figure}
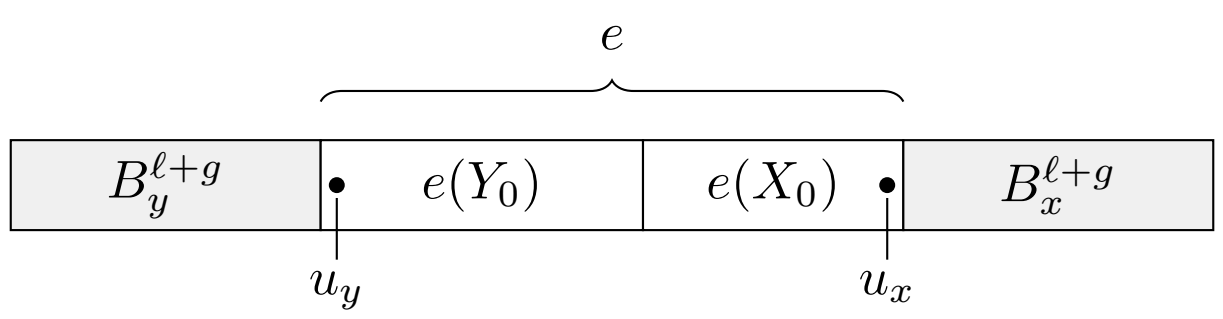

{\bf Case I:} If \(y\le n\) (and so $x\le n$), prepend to \(e\) a vertex sequence $B_y^{\ell+g}$ and append to $e$ a sequence $B_x^{\ell+g}$, that is, take
$$\Pi=B_y^{\ell+g}eB_x^{\ell+g},$$ see also Figure~\ref{eee}. Here $u_y$ and $u_x$ denote the first and last vertices of $e$, respectively.

By \eqref{rel_e}, $e$ is a relevant $k$-window of $\Pi$. Moreover,  since $\ell+g<k$, every  $k$-window $w$ of $\Pi$ containing an element to the right of~$e$ must contain $u_x$ too, and, by the choice of $x$, we have $\min(w)=x$. Similarly (but not quite symmetrically),  every \emph{relevant} $k$-window $w$ to the left of $e$ (that is, containing a multiple of $k-\ell$ elements to the left of $e$) must contain $u_y$ and, crucially, $w\cap e\subseteq e(Y_0)$. Thus, by the choice of $y$, we have $\min(w)=y$. This and  \eqref{eq:no-L1-in-e} imply that  $w\in H_2$   (membership in $\mathcal K$ holds vacuously), and thus \(\Pi\subset H_2+e\).

\medskip

{\bf Case II.} Suppose now that \(y>n\). By the definitions of $x$ and $y$, this assumption implies that $y'=y-n\not\in \operatorname{tr}(e)$, since $\{y,y'\}\in E(G)$. In what follows, we describe only the left part of \(\Pi\); the right part is obtained symmetrically by reversing the construction (and using $x$ instead of $y$)  when $x>n$, and by using the right hand side construction of Case I if $x\le n$. (Note that if $x>n$, then, trivially, $x'=x-n\not\in \operatorname{tr}(e)$.) Put
\[
    t:=t_y:=|e\cap A_y|
\]
 and notice that by~\eqref{eq:component-bound-main2} applied to the component of \(G^3[\operatorname{tr}(e)]\) containing $y$ we have $0\le t\le k-\ell$.

We find it convenient to introduce one more parameter related to $g$, namely
\begin{equation}\label{d}\delta:=k-\ell-g.\end{equation} Thus, \(1\le \delta\le k-\ell\), cf. \eqref{gap}.

Now, we reveal our construction which depends on the value of $t$.

\medskip

\begin{figure}[htbp]
\centering
\begin{lrbox}{\reviewfigurebox}
\begin{tikzpicture}[x=0.82cm,y=0.85cm,
 every node/.style={font=\small},
 brace/.style={decorate,decoration={brace,amplitude=4pt}}]
\fill[black!6] (0,-0.35) rectangle (2.2,0.35);
\draw (0,-0.35) rectangle (2.2,0.35);
\node at (1.1,0) {$A_{y'}^{\ell}$};
\fill[black!6] (2.2,-0.35) rectangle (3.4,0.35);
\draw (2.2,-0.35) rectangle (3.4,0.35);
\node at (2.8,0) {$A_y^g$};
\fill[black!6] (3.4,-0.35) rectangle (5.6,0.35);
\draw (3.4,-0.35) rectangle (5.6,0.35);
\node at (4.5,0) {$A_y^{\ell}$};
\fill[black!6] (5.6,-0.35) rectangle (6.9,0.35);
\draw (5.6,-0.35) rectangle (6.9,0.35);
\node at (6.25,0) {$B_y^g$};
\fill[white] (6.9,-0.35) rectangle (8.4,0.35);
\draw (6.9,-0.35) rectangle (8.4,0.35);
\node at (7.65,0) {$[e|\delta]$};
\fill[white] (8.4,-0.35) rectangle (10.9,0.35);
\draw (8.4,-0.35) rectangle (10.9,0.35);
\node at (9.65,0) {$e\setminus[e|\delta]$};
\draw[brace] (5.6,-0.65)--(2.2,-0.65) node[midway,below=6pt] {$\ell+g$};
\draw[brace] (6.9,0.65)--(10.9,0.65) node[midway,above=6pt] {$e$};
\draw[brace] (3.4,1.5)--(8.4,1.5) node[midway,above=6pt] {$\widehat e$};
\end{tikzpicture}
\end{lrbox}
\ifdim\wd\reviewfigurebox>\linewidth\resizebox{\linewidth}{!}{\usebox{\reviewfigurebox}}\else\usebox{\reviewfigurebox}\fi
\caption{Description of $\Pi$, case II.1.}
\label{fig:configuration2}
\end{figure}
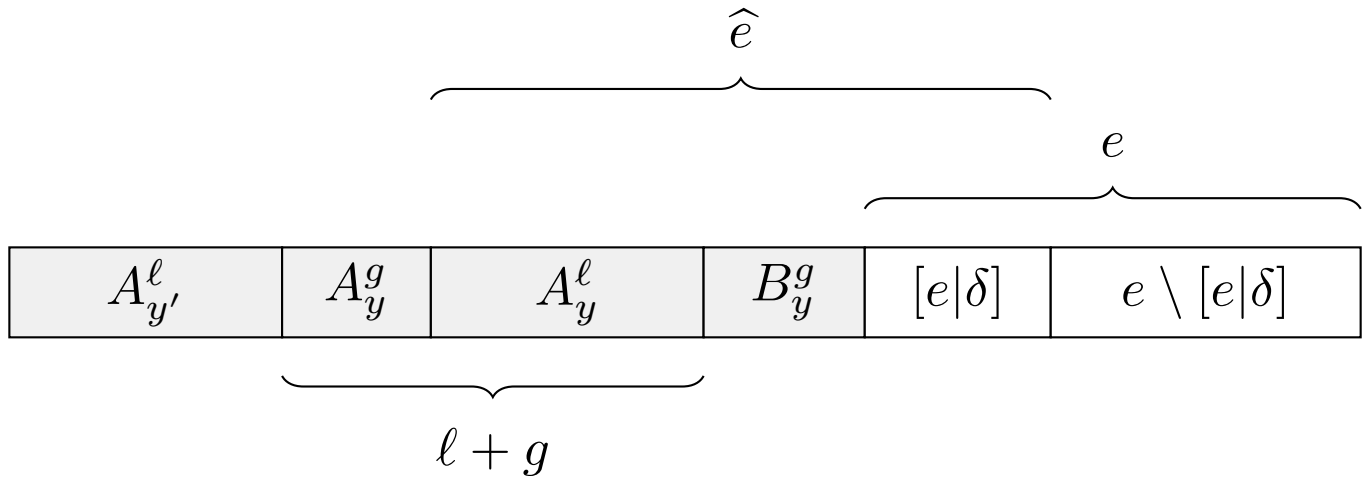

{\bf Case II.1.} If \(t\le \delta-1\), recall that $y'=y-n$ and prepend
\begin{equation}\label{eq:left-disjoint}
    A_{y'}^\ell A_y^{g+\ell}B_y^g
\end{equation}
to \(e\). This is possible since, by~\eqref{d} and~\eqref{mocai2026},
\[
    t+\ell+g
    \le (\delta-1)+\ell+g
    =k-1
    =|A_y|.
\]
\noindent
Let
\(
        \widehat e = A_y^\ell B_y^g [e|\delta],
\)
where \([e|\delta]\) denotes the set  of the first \(\delta\) vertices of \(e\), see Figure~\ref{fig:configuration2}. (Note that, by~\eqref{d}, $|\widehat e|=k$.)

\noindent

Recall that  we deem a $k$-window relevant if the distance between its initial position and the initial position  of $e$ is divisible by $k-\ell$.
We first verify that all relevant \(k\)-windows preceding \(e\) belong to
\(H_1\cup H_2\). Every such window $w$ preceding \(\widehat e\) is contained in
\(A_{y'}\cup A_y\) and meets both these sets (as  $g< k-\ell$); hence, $\operatorname{tr}(w)=\{y', y\}\in E(G)$ and, by definition, \(w\in H_1\).

Now, let \(w\) be any relevant $k$-window from \(\widehat e\) up to, but not including,
\(e\).  We are going to show that $w\in H_2$.
Note first that \(w\setminus e\)  is a suffix of
\[
        A_y^\ell B_y^g
\]
of size at least \(k-\ell\). Since \(g<k-\ell\), this suffix meets
\(A_y\). Moreover, since, $|e(Y_0)|\ge\ell$, we have \(w\cap e\subseteq e(Y_0)\). Thus
\[
        \min(w)=y
        \qquad\text{and}\qquad
        |w\cap U_y|\ge k-\ell+1,
\]
and it only remains to check that $w\in\mathcal K$. If
\[
        L\in\binom{w\cap A}{\ell}\cap\mathcal L,
\]
then, by~\eqref{eq:no-L1-in-e}, \(L\not\subseteq e\).
Hence,  \(L\cap (A_y\setminus e)\neq\emptyset\). Since the only neighbor of $y$ in $G$ is $y'$  and
\(y'=y-n\notin\operatorname{tr}(e)\), this, by the definition of $\mathcal L$, forces \(L\subseteq A_y\).
Therefore, \(w\in H_2\).

\medskip

{\bf Case II.2.} Let \(t\ge \delta\). Recall that the vertices of $e=e(X_0)\cup e(Y_0)$ have been ordered so that $e(Y_0)$ is to the left of $e(X_0)$. Now, in addition, if necessary, reorder the vertices within $e(Y_0)$, so that $e\cap A_y$ forms a prefix of $e$.
To construct (the left wing of)  $\Pi$, define
\begin{equation}\label{ha}
   h:= h_y:=\max\{0,t_y+2g-2(k-\ell)+1\}
\end{equation}
and note that, by  \eqref{gap},
\begin{equation}\label{eq:h-properties}
    0\le h\le g
    \qquad\mbox{and}\qquad
    h\le t-1,
\end{equation}
as well as, that $h=0$ when $t\le \delta-1$.

Now, prepend to $e$ the sequence
\begin{equation}\label{eq:left-overlap}
  W':=A_{y'}^{\ell+h}
    A_y^{\ell+g-h-\delta}.
\end{equation}
Together with the first \(\delta\) vertices of \(e\), this forms the sequence
\[
  W'':=A_{y'}^{\ell+h}A_y^{\ell+g-h}.
\]
Note that, by the definitions of \(h\) and $\delta$, the whole sequence $W'e$ uses
\[
    t+\ell+g-h-\delta\le k-1
\]
vertices of \(A_y\), and thus, in view of~\eqref{mocai2026}, is feasible.

We verify that the resulting sequence $W'e$ is an
\((\ell,k)\)-path in \((H_1\cup H_2)+e\). First note, recalling~\eqref{rel_e}, that  $\ell+g=p(k-\ell)$. This, in turn, implies, using~\eqref{d},
that the length of  prefix~\eqref{eq:left-overlap}  is\[
    2\ell+g-\delta=(2p-1)(k-\ell),
\]
so \(e\) itself is a relevant \(k\)-window.

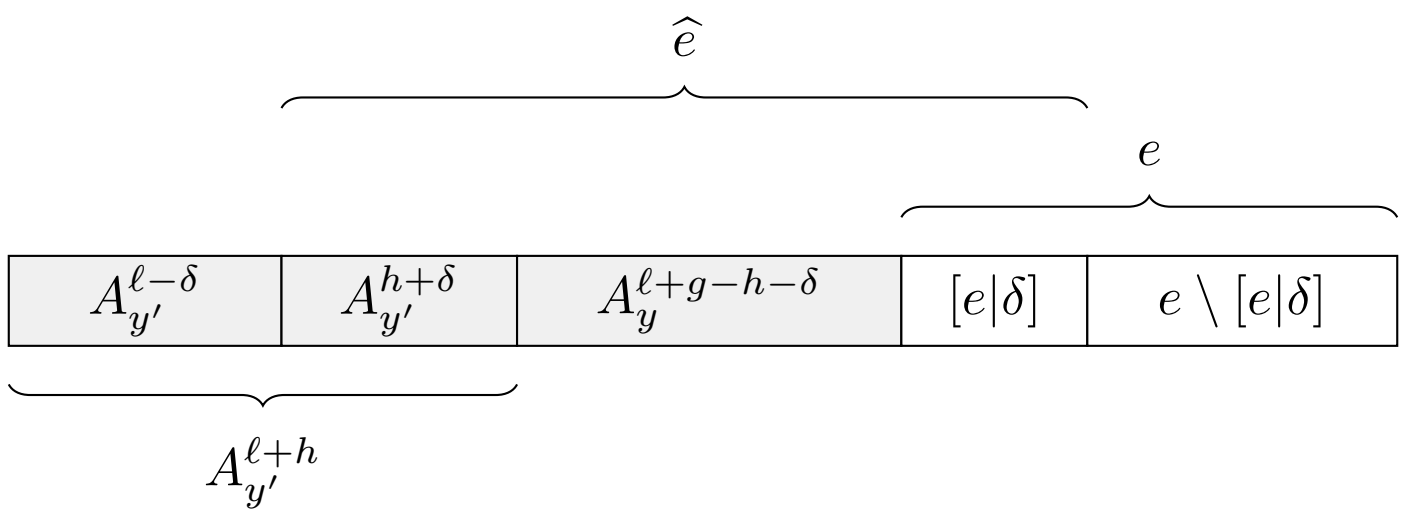
\begin{figure}[htbp]
\centering
\begin{lrbox}{\reviewfigurebox}
\begin{tikzpicture}[x=0.82cm,y=0.85cm,
 every node/.style={font=\small},
 brace/.style={decorate,decoration={brace,amplitude=4pt}}]
\fill[black!6] (0,-0.35) rectangle (2.2,0.35);
\draw (0,-0.35) rectangle (2.2,0.35);
\node at (1.1,0) {$A_{y'}^{\ell-\delta}$};
\fill[black!6] (2.2,-0.35) rectangle (4.1,0.35);
\draw (2.2,-0.35) rectangle (4.1,0.35);
\node at (3.15,0) {$A_{y'}^{h+\delta}$};
\fill[black!6] (4.1,-0.35) rectangle (7.2,0.35);
\draw (4.1,-0.35) rectangle (7.2,0.35);
\node at (5.65,0) {$A_y^{\ell+g-h-\delta}$};
\fill[white] (7.2,-0.35) rectangle (8.7,0.35);
\draw (7.2,-0.35) rectangle (8.7,0.35);
\node at (7.949999999999999,0) {$[e|\delta]$};
\fill[white] (8.7,-0.35) rectangle (11.2,0.35);
\draw (8.7,-0.35) rectangle (11.2,0.35);
\node at (9.95,0) {$e\setminus[e|\delta]$};
\draw[brace] (4.1,-0.65)--(0,-0.65) node[midway,below=6pt] {$A_{y'}^{\ell+h}$};
\draw[brace] (7.2,0.65)--(11.2,0.65) node[midway,above=6pt] {$e$};
\draw[brace] (2.2,1.5)--(8.7,1.5) node[midway,above=6pt] {$\widehat e$};
\end{tikzpicture}
\end{lrbox}
\ifdim\wd\reviewfigurebox>\linewidth\resizebox{\linewidth}{!}{\usebox{\reviewfigurebox}}\else\usebox{\reviewfigurebox}\fi
\caption{Description of $\Pi$, case II.2. The overlap $\widehat e\cap e=[e|\delta]$ lies in $A_y$.}
\label{fig:configuration3}
\end{figure}

Let \(\widehat e\) be the relevant $k$-window which overlaps with \(e\) in precisely its first
$\delta=k-\ell-g$ vertices -- Figure~\ref{fig:configuration3}.

\noindent

Then, since each block in $W''$, that is,
$ A_{y'}^{\ell+h}$ and $A_y^{\ell+g-h}$,
is shorter than \(k\) (cf. \eqref{eq:h-properties} and \eqref{gap}), every \(k\)-window
contained in $W''$ (that is,  to the left of $\widehat e$, including $\widehat e$) meets both \(A_y\) and \(A_{y'}\),
 and hence belongs to \(H_1\).

 On the other hand, every remaining relevant $k$-window \(w\ne e\) (that is, to the right of $\widehat e$)
 satisfies $|w\cap e| \ge \delta+k-\ell$ and thus
 $$ |w\setminus e| \le k - (k-\ell+\delta)=\ell-\delta\overset{\eqref{eq:h-properties}}\leq \ell+g-h-\delta.$$
 It follows that $w\setminus e\subset A_y$.
 But also, $|w\setminus e|\ge k-\ell$ which, together with $\delta\ge1$ and the previous observation, yields
  $|w\cap U_{y}|\ge|w\cap A_y|\ge k-\ell+1$. Moreover, since $|e(Y_0)|\ge\ell$ and $e(Y_0)$ is to the left of $e(X_0)$ in $e$,
 we have $y = \min(w)$.
 The argument that $w\in\mathcal K$ is identical to the one in Case II.1 and therefore is omitted.
 Thus,  \(w\in H_2\).

This completes the construction of \(\Pi\).  Denote the terminal
\(\ell\)-sets of $\Pi$  by $T_1$ (on the right) and $T_n$ (on the left) and note that
\[
    T_1\subseteq U_{v_1},
    \quad\mbox{and}\quad
    T_n\subseteq U_{v_n}
\]
 (recall that $v_1=x'$ and $v_n=y'$). According to the construction, if  $x>n$, then
 \(T_i\subseteq A_{v_i}\), $i=1,n$,  if $y\le n$, then  \(T_i\subseteq B_{v_i}\), $i=1,n$, while in the remaining case $T_1\subseteq B_{v_1}$ and $T_n\subseteq A_{v_n}$.

 Note that in all cases the constructed path $\Pi$ satisfies $|V(\Pi)|<5k$.

\subsection{Numerical preparations}\label{numpre}
In this section we first measure what portion of each set \(A_{v_i}\), $i\in[n]$, has been consumed by the path $\Pi$ constructed in the previous section (not counting the two terminals, $T_1$ and $T_n$, though).

Recall that, outside $e$, the only  sets $A_i$, $i\le n$, used by $\Pi$ are $A_{v_1}$ and $A_{v_n}$. More precisely, it happens only if $x>n$ or $y>n$. For instance, if $y>n$, then $|(V(\Pi)\setminus (e\cup T_n))\cap  A_{v_n}|=h_y$, where $h_y$ was defined in~\eqref{ha}. Moreover, then $y$ is a pendant vertex of $G$, and so \(v_n\notin\operatorname{tr}(e)\),  implying that $e\cap A_{v_n}=\emptyset$, where recall $v_n=y'$ (and, by symmetry, all of this is true for $x$ and $v_1=x'$ as well).

Hence, we define
 for all $i\in\{2,\dots,n-1\}$,
\[
    p_i:=|e\cap A_{v_i}|,
\]
while
$$p_1:=\begin{cases}|e\cap A_{v_1}|\quad\mbox{if}\quad x\le n\\ h_x \quad\mbox{if}\quad x>n \end{cases}$$
and
$$p_n:=\begin{cases}|e\cap A_{v_n}|\quad\mbox{if}\quad y\le n\\ h_y \quad\mbox{if}\quad y>n. \end{cases}$$
(Recall that if $t_y\le \delta-1$, which corresponds to Case II.1 of the construction,  then $h_y=0$, and so $p_n=0$ -- and the same is true for $x$ and $p_1$.)


\begin{fact}\label{sumpi}
For all $I\subseteq[n]$,
\begin{equation}\label{eq:p-interval-corrected}
    \sum_{i\in I}p_i\le k-\ell-1+|I|.
    \end{equation}
    In particular, for each $i\in[n]$, one has $0\le p_i\le k-\ell$.
     Moreover, \space
     \begin{equation}\label{npk}
 |\{i:p_i>0\}|\le k.
     \end{equation}
    
\end{fact}

\proof
Put $w_1=x$, $w_n=y$, and $w_i=v_i$ for
$i=2,\ldots,n-1$.
The sets $e\cap A_{w_i}$, $i\in[n]$, are
pairwise disjoint subsets of $e$.
Moreover, for $z\in\{x,y\}$, if $z>n$,
then $h_z\le t_z=|e\cap A_z|$. Thus,
$$
\sum_{i\in I}p_i\le |e'|
\quad\mbox{where}\quad
e'=\bigcup_{i\in I}(e\cap A_{w_i})
$$
and $|\operatorname{tr}(e')|\le |I|$.
In particular, $\sum_{i=1}^n p_i\le k$ which implies~\eqref{npk}.
In turn, \eqref{eq:p-interval-corrected} follows
from \eqref{eq:subset-bound-main2}.
\qed

\medskip

Guided by an analogy with the proof of Lemma~\ref{main1}  (Case 1), \space
let \(r'_1,\ldots,r'_n\) achieve the maximum defining \(\Phi(\alpha,n)\)  with a
witnessing sequence

\[
        s'_1,\ldots,s'_{n-1}\in\{g,g+k-\ell\}
\]
(cf. Definitions~\ref{defphi} and~\ref{defadm}). Such a choice is possible by Claim~\ref{just2}.
(We use primes here, because this is not our final sequence.) Since the sets $A_{v_i}$, while constructing path $\Pi$,  have suffered an additional loss  of $p_i$ vertices, we need to adjust the sequence $(r_i')_1^n$ to $(r_i)_1^n$, so that $r_i\le k-\ell-p_i$ (and not just $r_i'\le k-\ell$). This is achieved by the following lemma the proof of which is presented in Section~\ref{sec:correct-loads}. Note that, in principle, we do not know which of $r_i$ and $r_i'$ is larger.

\begin{lem}\label{88}
Let \(r'_1,\ldots,r'_n\) be an admissible sequence with a
witnessing sequence $s'_1,\ldots,s'_{n-1}$ with $s'_i\in\{g,g+k-\ell\}$ for every $i$, while $p_1,\dots,p_n$ be a nonnegative integer sequence satisfying the conclusions of Fact~\ref{sumpi}. Then there exists an admissible sequence \(r_1,\ldots,r_n\) with a
witnessing sequence $s_1,\ldots,s_{n-1}$ such that $s_i\in\{g,g+k-\ell\}$

\begin{equation}\label{eq:corrected-loads}
        0\le r_i\le k-\ell-p_i, \quad\mbox{for every $i$},
\end{equation}
and 
\begin{equation}\label{eq:corrected-stability}
        \sum_{i=1}^n|r_i-r'_i|
        \le 2k^3.
        \end{equation}

\end{lem}
\subsection{Constructing a (too long) cycle}\label{toolong}
In this section, we construct an auxiliary $(\ell,k)$-cycle $C'$ of length $N'>N$, containing all vertices of \space $\bigcup_{i=1}^nU_i\cup  V(\Pi)$, while all other vertices will be taken from an auxiliary, infinite set $W^*$. Later, we will shorten this cycle to length $N$ and replace all ``fake'' vertices by the missing vertices of $\bigcup_{i=n+1}^{2n}U_i$. This way, we will ultimately  obtain a Hamilton $(\ell,k)$-cycle $C$ with $V(C)=U=\bigcup_{i=1}^{2n}U_i$. Moreover, our construction will guarantee that $C\subseteq H_1\cup H_2+e$.

Cycle $C'$ will consist, besides path $\Pi$ already constructed in Section~\ref{Qzero}, of vertex-disjoint paths $Q'_1,\dots, Q'_n$, ``glued'' together by paths $P_i$, $i=1,\dots, n-1$. In turn, paths $Q'_1$ and $Q'_n$ will be connected with $\Pi$, so, in a sense, $\Pi$ plays the role of $P_n$ (see Section~\ref{notHam}.)

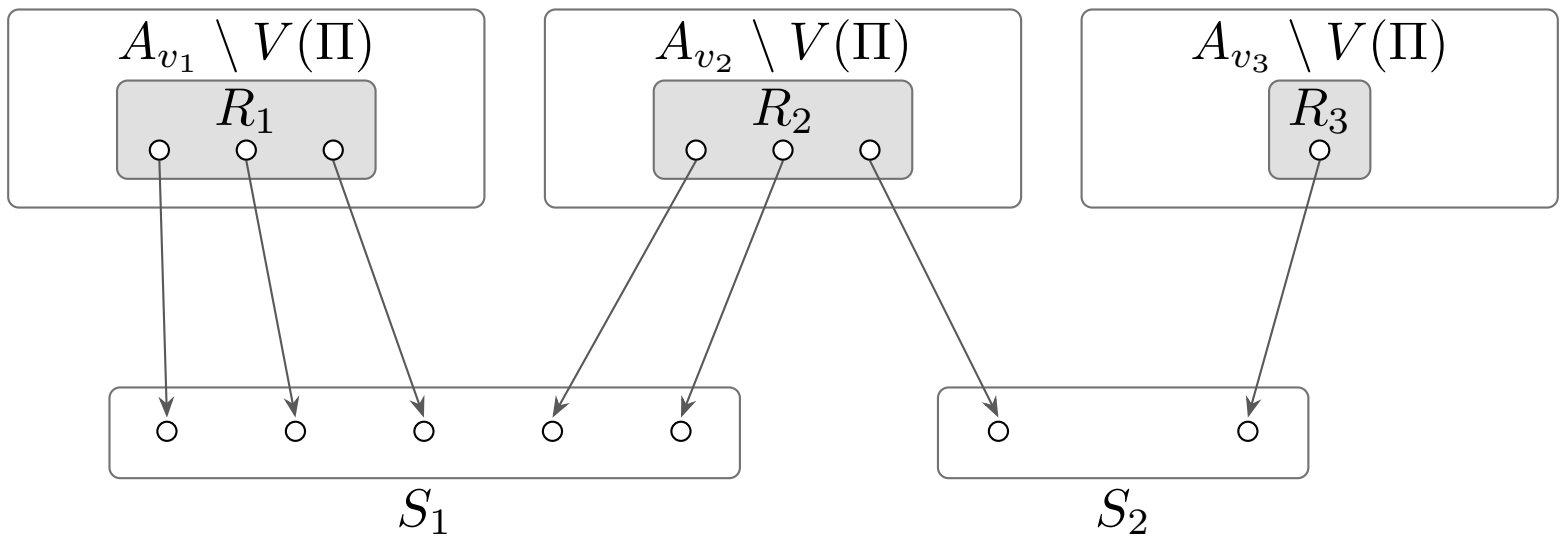
\begin{figure}[htbp]
\centering
\begin{lrbox}{\reviewfigurebox}
\begin{tikzpicture}[
  x=1cm,y=1cm,
  every node/.style={font=\small,inner sep=2pt},
  available/.style={draw=black!55,rounded corners=2pt},
  reserved/.style={draw=black!55,fill=black!12,rounded corners=2pt},
  assigned/.style={-{Stealth[length=4pt]},draw=black!65,thin}
]
  \foreach \a/\b/\m/\i in {0/3.15/1.575/1,3.55/6.70/5.125/2,7.10/10.25/8.675/3}{
    \draw[available] (\a,1.79) rectangle (\b,3.10);
    \node at (\m,2.86) {$A_{v_\i}\setminus V(\Pi)$};
  }
  \foreach \a/\b/\m/\i in {0.72/2.43/1.575/1,4.27/5.98/5.125/2,8.34/9.01/8.675/3}{
    \draw[reserved] (\a,1.98) rectangle (\b,2.63);
    \node at (\m,2.43) {$R_\i$};
  }
  \draw[available] (0.67,0.00) rectangle (4.84,0.60);
  \draw[available] (6.15,0.00) rectangle (8.60,0.60);
  \node[below] at (2.755,0.00) {$S_1$};
  \node[below] at (7.375,0.00) {$S_2$};
  \foreach \u/\v in {1.00/1.05,1.575/1.90,2.15/2.75,4.55/3.60,5.125/4.45,5.70/6.55,8.675/8.20}{
    \draw[assigned] (\u,2.10) -- (\v,0.40);
    \draw[fill=white] (\u,2.17) circle[radius=1.8pt];
    \draw[fill=white] (\v,0.31) circle[radius=1.8pt];
  }
\end{tikzpicture}
\end{lrbox}
\ifdim\wd\reviewfigurebox>\linewidth\resizebox{\linewidth}{!}{\usebox{\reviewfigurebox}}\else\usebox{\reviewfigurebox}\fi
\caption{Assigning reserved vertices to the sets $S_i$.
  In this example, $k=7$, $\ell=4$, $(r_1,r_2,r_3)=(3,3,1)$ and
  $(s_1,s_2)=(5,2)$. Two vertices of $R_2$ go to $S_1$ and one to $S_2$.
  Each reserved vertex is assigned exactly once.}
\label{fig:reserved-sets}
\end{figure}

Given the admissible sequence $(r_1,\dots,r_n)$ and its witness $(s_1,\dots,s_{n-1})$, both delivered by Lemma~\ref{88}, choose sets 
\[
    R_i\subseteq A_{v_i}\setminus V(\Pi),
    \qquad
    |R_i|=r_i,\quad i=1,\dots,n,
\]
and  pairwise disjoint sets
\[
    S_i\subseteq R_i\cup R_{i+1},
    \qquad
    |S_i|=s_i,\quad i=1,\dots,n-1
\]
such that
\[\bigcup_{i=1}^nR_i=\bigcup_{i=1}^{n-1}S_i.\]
This is possible by~\eqref{follows} and the right-hand side of~\eqref{adm} (see Figure~\ref{fig:reserved-sets} for the allocation).

It is crucial that, after constructing $\Pi$ and ``reserving'' the sets $R_i$, there is still enough room in the sets $A_{v_i}$ to form terminal $\ell$-sets of the paths $Q_1,\dots,Q_n$ we are going to build.
Indeed, if $i\in\{2,\dots,n-1\}$, or $i=1$ and $x=v_1$, or $i=n$ and $y=v_n$,  then, by~\eqref{mocai2026} and~\eqref{eq:corrected-loads},
\begin{equation}\label{easy_case}
    |A_{v_i}\setminus(V(\Pi)\cup R_i)|
    \ge k+\ell-p_i-r_i
    \ge 2\ell
\end{equation}
and two terminals can be formed. In the remaining two cases, $i=1$ and $x>n$, or $i=n$ and $y>n$, since $p_i$ does not count the vertices of $T_i$, we only have
\begin{equation}\label{hard_case}
    |A_{v_i}\setminus(V(\Pi)\cup R_i)|
    \ge k+\ell-(\ell+p_i)-r_i
    \ge\ell,
\end{equation}
but this is fine, because we just need one more terminal $\ell$-set.

We first build paths $Q'_i$, $i=1,\dots,n$.
To this end, for $i=1,\dots,n$, let $W^*_i$  be infinite disjoint subsets of $W^*$,
while
\[
    D_i:=
    U_{v_i}\setminus(V(\Pi)\cup R_i).
\]
Further, let
\[
    \sigma_i:=
    |U_{v_i}\cap V(\Pi)|,
\]
so that
\[
    |D_i|=\alpha_{v_i}-\sigma_i-r_i
    \quad\mbox{and}\quad
    \sum_{i=1}^{n}\sigma_i\le |V(\Pi)|< 5k.
\]

For $i=2,\ldots,n-1$, let $Q'_i$ be a
$(D_i,W^*_i)$-dominant path of maximum
length, with both terminal $\ell$-sets
contained in $A_{v_i}$. This is possible by~\eqref{easy_case} and, if necessary, by swapping around some vertices of $A_{v_i}$ and $B_{v_i}$.

\begin{figure}[htbp]
\centering
\begin{lrbox}{\reviewfigurebox}
\begin{tikzpicture}[
  x=1cm,y=1cm,
  every node/.style={font=\small,inner sep=2pt},
  shared/.style={fill=black!20},
  brace/.style={decorate,decoration={brace,amplitude=4pt}}
]
  \fill[black!8] (0,-0.30) rectangle (0.95,0.30);
  \fill[shared] (2.30,-0.30) rectangle (3.50,0.30);
  \fill[shared] (6.28,-0.30) rectangle (7.48,0.30);
  \fill[black!8] (8.83,-0.30) rectangle (9.78,0.30);
  \draw (0,-0.30) rectangle (9.78,0.30);
  \foreach \x in {0.95,2.30,3.50,6.28,7.48,8.83}{
    \draw (\x,-0.30) -- (\x,0.30);
  }
  \foreach \x/\t in {0.475/\ell,1.625/\cdots,2.90/T_n,4.89/\cdots,6.88/T_1,8.155/\cdots,9.305/\ell}{
    \node at (\x,0) {$\t$};
  }
  \draw[brace] (0,0.45) -- (3.50,0.45)
    node[midway,above=5pt] {$Q'_n$};
  \draw[brace] (6.28,0.45) -- (9.78,0.45)
    node[midway,above=5pt] {$Q'_1$};
  \draw[brace] (7.48,-0.45) -- (2.30,-0.45)
    node[midway,below=5pt] {$\Pi$};
  \node[below=4pt] at (0.475,-0.30) {in $A_{v_n}$};
  \node[below=4pt] at (9.305,-0.30) {in $A_{v_1}$};
\end{tikzpicture}
\end{lrbox}
\ifdim\wd\reviewfigurebox>\linewidth\resizebox{\linewidth}{!}{\usebox{\reviewfigurebox}}\else\usebox{\reviewfigurebox}\fi
\caption{Joining $Q'_n$, $\Pi$, and $Q'_1$ in the orientation of $\Pi$. The darker blocks are the shared terminal $\ell$-sets $T_n$ and $T_1$.}
\label{fig:pi-join}
\end{figure}
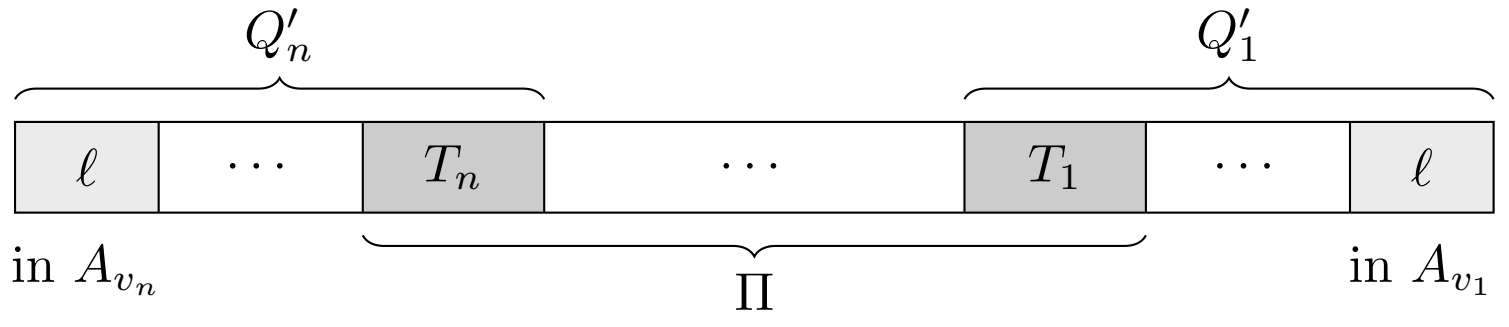

For $i\in\{1,n\}$, the argument needs to be more subtle. Let $T_i'$ be an abstract
set of $\ell$ vertices, disjoint from $D_i\cup W^*_i$,
and put $D_i'=D_i\cup T_i'$.
Choose a $(D_i',W^*_i)$-dominant path $Q'_i$
of maximum length. By moving around the vertices
of $D_i'$, we may arrange that one terminal
$\ell$-set is $T_i'$ and the other is
contained in $A_{v_i}\setminus V(\Pi)$.
This is possible by~\eqref{hard_case} and the fact that the two terminal $\ell$-sets are
disjoint, since $\alpha$ is sufficiently large.
Finally, replace $T_i'$ by $T_i$ obtaining the ultimate path $Q'_i$ for $i\in\{1,n\}$. Consequently,
\begin{equation}\label{QPiT} 
V(Q'_i)\cap V(\Pi)=T_i,
\qquad i=1,n,
\end{equation}
so the paths $Q'_1,\Pi,$ and $Q'_n$ now form one $(\ell,k)$-path, while all paths $Q'_i$, $i=1,\dots,n$,  are pairwise vertex disjoint and disjoint from $\Pi$, except for $Q'_1$ and  $Q'_n$ (see Figure~\ref{fig:pi-join}).

As for their lengths, by the definition of function $\nu$, for $i=2,\dots,n-1$,
$$ |V(Q'_i)|=\nu(\alpha_{v_i}-\sigma_i-r_i),$$
while, for $i\in\{1,n\}$, by Lemma~\ref{nu1}(c),

\begin{equation*}
    |V(Q'_i)|= \nu(\alpha_{v_i}-\sigma_i-r_i+\ell)\ge\nu(\alpha_{v_i}-\sigma_i-r_i).
\end{equation*}
In summary, for each $i\in[n]$, we have
\begin{equation}\label{eq:Qi-length}
    |V(Q'_i)|\ge\nu(\alpha_{v_i}-\sigma_i-r_i).
\end{equation}

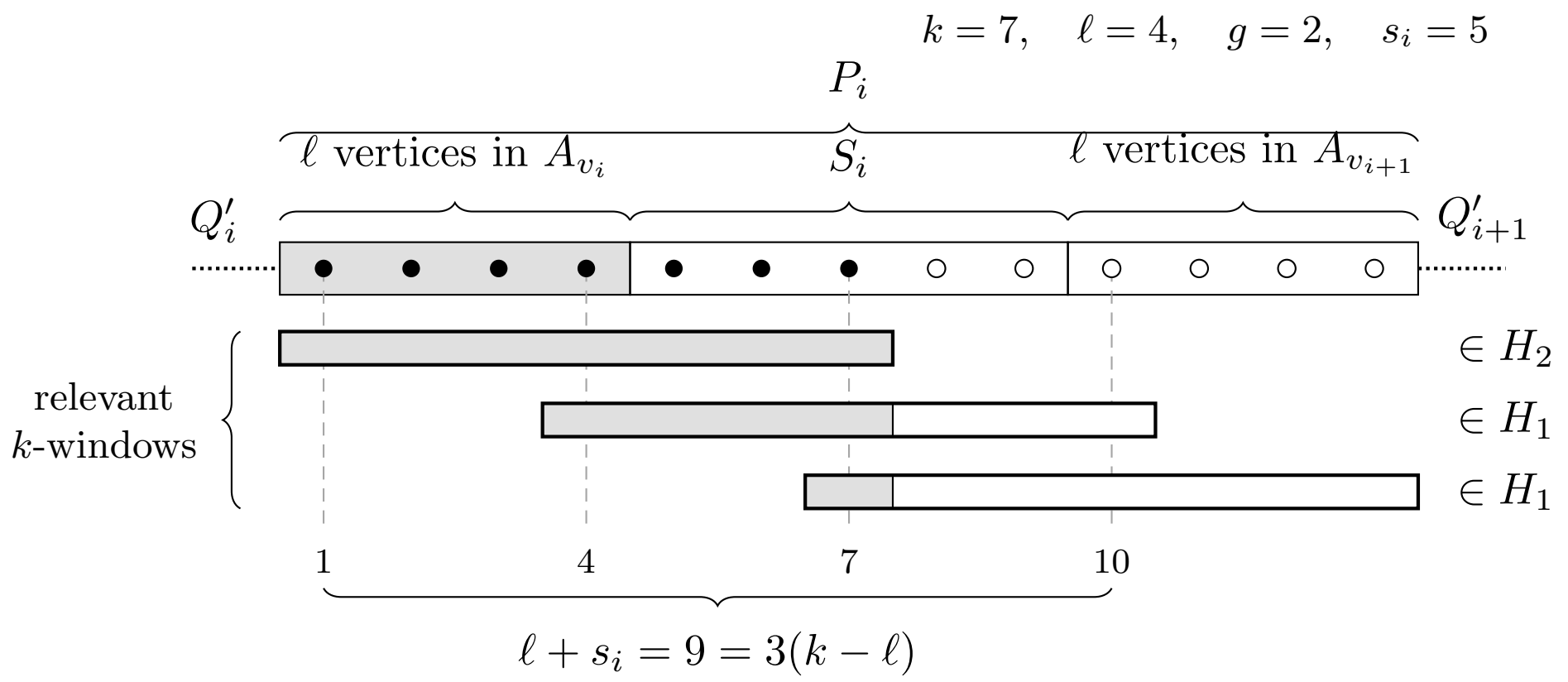
\begin{figure}[htbp]
\centering
\begin{lrbox}{\reviewfigurebox}\begin{tikzpicture}[
  x=0.70cm,y=0.70cm,
  every node/.style={font=\small},
  brace/.style={decorate,decoration={brace,amplitude=4pt}},
  guide/.style={densely dashed,black!35,thin}
]
  \node[anchor=east,font=\footnotesize] at (14.0,2.70)
    {$k=7,\quad\ell=4, \quad g = 2, \quad s_i=5$};

  \fill[black!12] (0,-0.30) rectangle (4,0.30);
  \draw (0,-0.30) rectangle (4,0.30);
  \draw (4,-0.30) rectangle (9,0.30);
  \draw (9,-0.30) rectangle (13,0.30);
  \foreach \p in {1,4,7,10}{
    \draw[guide] ({\p-0.5},0)--({\p-0.5},-3.00);
    \node[below,font=\scriptsize] at ({\p-0.5},-3.02) {$\p$};
  }
  \foreach \j in {0,...,6}
    \fill (\j+0.5,0) circle (2pt);
  \foreach \j in {7,...,12}
    \draw[fill=white] (\j+0.5,0) circle (2pt);
  \draw[densely dotted,thick] (-1.0,0)--(0,0);
  \draw[densely dotted,thick] (13,0)--(14.0,0);
  \node[above] at (-0.75,0.10) {$Q_i^{\prime}$};
  \node[above] at (13.75,0.10) {$Q_{i+1}^{\prime}$};

  \draw[brace] (0,0.55)--(4,0.55)
    node[midway,above=6pt] {$\ell$ vertices in $A_{v_i}$};
  \draw[brace] (4,0.55)--(9,0.55)
    node[midway,above=6pt] {$S_i$};
  \draw[brace] (9,0.55)--(13,0.55)
    node[midway,above=6pt] {$\ell$ vertices in $A_{v_{i+1}}$};
  \draw[brace] (0,1.45)--(13,1.45)
    node[midway,above=6pt] {$P_i$};

  \fill[black!12] (0,-1.10) rectangle (7,-0.72);
  \draw[thick] (0,-1.10) rectangle (7,-0.72);
  \node[right] at (13.25,-0.91) {$\in H_2$};

  \fill[black!12] (3,-1.92) rectangle (7,-1.54);
  \fill[white] (7,-1.92) rectangle (10,-1.54);
  \draw[thick] (3,-1.92) rectangle (10,-1.54);
  \draw (7,-1.92)--(7,-1.54);
  \node[right] at (13.25,-1.73) {$\in H_1$};

  \fill[black!12] (6,-2.74) rectangle (7,-2.36);
  \fill[white] (7,-2.74) rectangle (13,-2.36);
  \draw[thick] (6,-2.74) rectangle (13,-2.36);
  \draw (7,-2.74)--(7,-2.36);
  \node[right] at (13.25,-2.55) {$\in H_1$};
  \draw[brace] (-0.45,-2.74)--(-0.45,-0.72)
    node[midway,left=6pt,align=center,font=\footnotesize] {relevant\\$k$-windows};

  \draw[brace] (9.5,-3.65)--(0.5,-3.65)
    node[midway,below=6pt] {$\ell+s_i=9=3(k-\ell)$};
\end{tikzpicture}\end{lrbox}
\ifdim\wd\reviewfigurebox>\linewidth\resizebox{\linewidth}{!}{\usebox{\reviewfigurebox}}\else\usebox{\reviewfigurebox}\fi
\caption{Joining $Q'_i$ to $Q'_{i+1}$ through $P_i$.
Filled and open circles represent vertices of $A_{v_i}$ and $A_{v_{i+1}}$, respectively.
The numbers below are vertex positions in $P_i$; the three displayed relevant windows start at positions $1$, $4$, and $7$.}
\label{fig:joining-paths}
\end{figure}

To complete the construction of $C'$,
for $i=1,\dots,n-1$,  join \(Q'_i\) with \(Q'_{i+1}\) by a path $P_i$ containing the set \(S_i\) (in any order) and the two 
corresponding terminal \(\ell\)-sets. Every resulting relevant \(k\)-window is
either contained in one of \(A_{v_i}\) and \(A_{v_{i+1}}\), in which
case it belongs to \(H_2\), or meets both sets, in which case it
belongs to \(H_1\) (since \(v_iv_{i+1}\in E(G_1)\)). 
The connection works, because $\ell+s_i$ is divisible by $k-\ell$ (see Figure \ref{fig:joining-paths}).

This completes the construction of an auxiliary
\((\ell,k)\)-cycle $C'$ containing~\(e\).
Observe that, indeed, 
\begin{equation}\label{Cprim}
V(C')\cap U=\bigcup_{i=1}^nU_i\cup  V(\Pi).
\end{equation}

It remains to estimate its length \(N':=|V(C')|\).
By construction, recalling that $e\cap V(Q'_1\cup Q'_n)=\emptyset$,  we have
\[
    N'=|V(Q'_1\cup\Pi\cup Q'_n)|+\sum_{i=2}^{n-1}|V(Q'_i)|+\sum_{i=1}^{n-1}|S_i|
    \ge
    \sum_{i=1}^{n}|V(Q'_i)|+k
    +\sum_{i=1}^{n-1}s_i.
\]

By~\eqref{eq:Qi-length} and Lemma~\ref{nu1}(e) (with $k$ instead of $k-1$ for convenience), 
\[
        |V(Q'_i)|
        \ge\nu(\alpha_{v_i}-\sigma_i-r_i)\ge \nu(\alpha_{v_i}-r_i)-\sigma_ik,
\]
while, again, by Lemma~\ref{nu1}(e) (when $r_i>r_i')$) or by Lemma~\ref{nu1}(c) (when $r_i<r_i')$), 
$$\nu(\alpha_{v_i}-r_i)\ge\nu(\alpha_{v_i}-r'_i)-\max\{0,r_i-r_i'\}k\ge \nu(\alpha_{v_i}-r'_i)-|r_i-r_i'|k.$$
Summing up and using Lemma~\ref{nu1}(e) one more time, \eqref{eq:corrected-stability}, and the bound
\(\sum_i\sigma_i\le 5k\), we get
\[
        \sum_{i=1}^n |V(Q'_i)|
        \ge
        \sum_{i=1}^n\nu(\alpha-r'_i)-(5k+2k^3)k-(k-\ell)kh_N,
\]
where, recall, $h_N$ counts the number of indices $i$ for which $\alpha_i=\alpha-k+\ell$ (see~\eqref{alphai}).
Moreover, by the definition of admissible sequences (cf. the right-hand side of~\eqref{adm}) and, again, by
\eqref{eq:corrected-stability},
\[
        \sum_{i=1}^{n-1}s_i
        =
        \sum_{i=1}^n r_i
        \ge\sum_{i=1}^n\left( r'_i-|r_i-r_i'|\right)
        \ge\sum_{i=1}^n r'_i-2k^3.
\]
Consequently, by definitions of function $\Phi$ (cf.  Definition~\ref{defphi}) and constant $c_+$, and the right-hand side of~\eqref{nn},

\begin{equation}
\begin{aligned}\label{NprimN}
        N'
        &\ge
        \sum_{i=1}^n
        \bigl(\nu(\alpha-r'_i)+r'_i\bigr)-c_+-(k-\ell)kh_N +k\\
        &=\Phi(\alpha,n)-c_+-(k-\ell)kh_N+k>N.
\end{aligned}
\end{equation}

\subsection{Getting a Hamilton cycle}\label{getham}
In this section we will first truncate $C'$ to the desired length $N$ and then replace all reservoir vertices by the remaining vertices of $U$.
The following fact will take care of the truncation part.

\begin{fact}\label{trunc}
For every $(D,W^*)$-dominant path $Q'$ with at least $k-\ell$ reservoir vertices there exists a $(D,W^*)$-dominant path $Q''$ with $|V(Q'')|=|V(Q')|-(k-\ell)$ and $V(Q'')\subset V(Q')$ which preserves the order of vertices of $Q'$.
\end{fact}

\proof To obtain $Q''$ just delete the first $k-\ell$ reservoir vertices from $Q'$. Let the last of them be $v$ at position $p$ in $Q'$. We need to check that $Q''$ still satisfies condition (ii) in Definition~\ref{dom}.

Let \(p\) be the position in \(Q'\) of the last deleted reservoir
vertex. Then the first \(p-(k-\ell)\) vertices of \(Q''\) all belong to \(D\). Let $V_1$ be the set of those vertices and set $V_2=V(Q'')\setminus V_1$.
Consider a relevant \(k\)-window \(w''\) of \(Q''\).

If $w''\subset V_1$, then all its vertices
belong to \(D\) and (ii) holds. If \(w''\cap V_1=\emptyset\), then, as all deleted
vertices precede $w''$, we have $w''=w'$, where $w'$ is a relevant $k$-window of $Q'$ whose initial position is shifted by \(k-\ell\). 

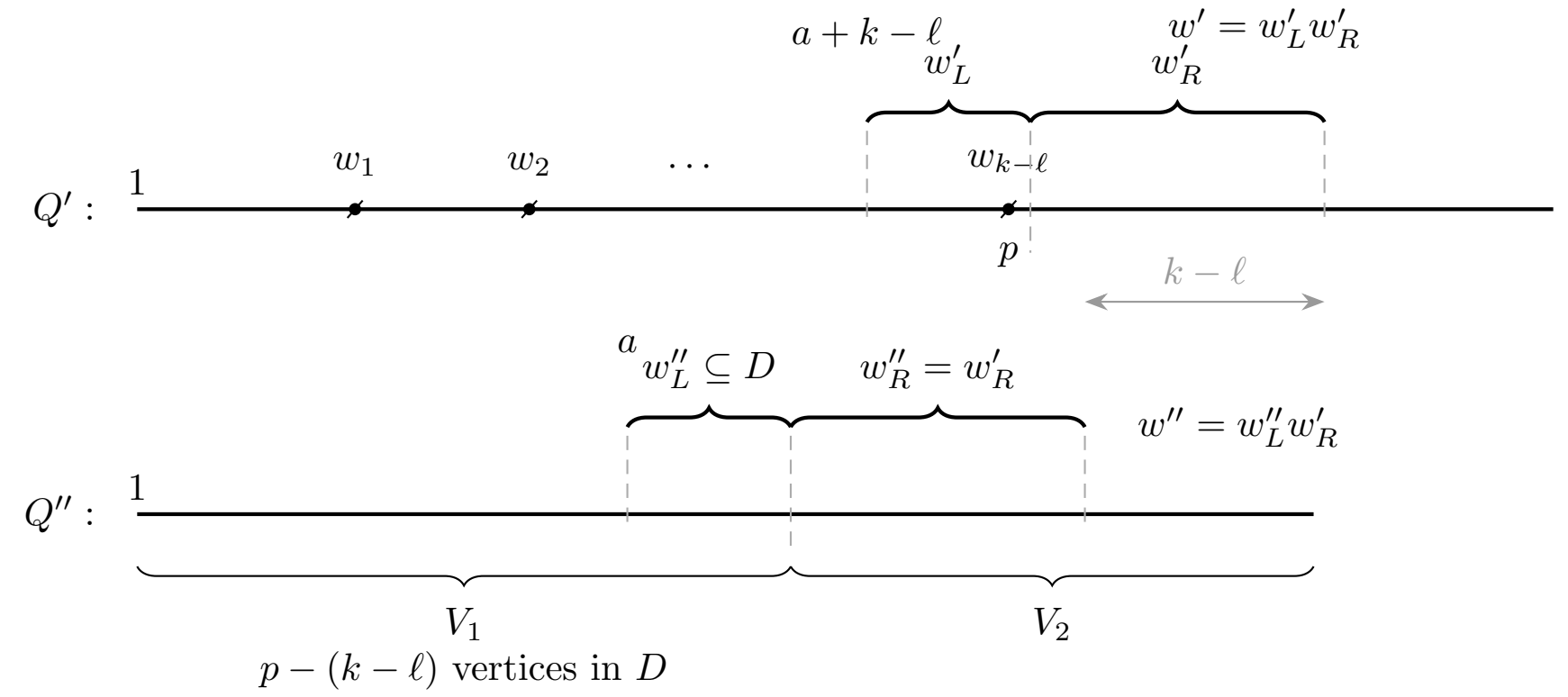
\begin{figure}[ht]
\centering
\begin{lrbox}{\reviewfigurebox}\begin{tikzpicture}[scale=0.8,transform shape, x=1cm,y=1cm,>=Stealth,
  every node/.style={font=\small},
  leftpart/.style={black,thick},
  rightpart/.style={black,thick},
  guide/.style={gray!65,densely dashed},
  brace/.style={decorate,decoration={brace,amplitude=4pt}}]

\def\shift{2.2}
\def\startnew{4.5}
\def\cutnew{6.0}
\def\endnew{8.7}
\def\rownew{-2.8}

\draw[thick] (0,0) -- (13,0);
\draw[thick] (0,\rownew) -- (10.8,\rownew);
\node[left] at (-.25,0) {$Q^{\prime}:$};
\node[left] at (-.25,\rownew) {$Q^{\prime\prime}:$};
\node[above] at (0,0) {$1$};
\node[above] at (0,\rownew) {$1$};

\foreach \xx/\lab in {2.0/{w_1},3.6/{w_2},8.0/{w_{k-\ell}}}{
  \fill[black] (\xx,0) circle (1.6pt);
  \draw[black] (\xx-.07,-.08) -- (\xx+.07,.08);
  \node[above=5pt,text=black] at (\xx,0) {$\lab$};
}
\node[above=5pt] at (5.1,0) {$\cdots$};
\node[below=5pt] at (8.0,0) {$p$};

\draw[leftpart,brace]
  ({\startnew+\shift},.80) -- ({\cutnew+\shift},.80)
  node[midway,above=6pt] {$w'_L$};
\draw[rightpart,brace]
  ({\cutnew+\shift},.80) -- ({\endnew+\shift},.80)
  node[midway,above=6pt] {$w'_R$};
\node[above] at ({\startnew+\shift},1.35) {$a+k-\ell$};
\node[above] at (10.35,1.35) {$w'=w'_Lw'_R$};
\draw[guide] ({\startnew+\shift},.72) -- ({\startnew+\shift},-.15);
\draw[guide] ({\cutnew+\shift},.72) -- ({\cutnew+\shift},-.40);
\draw[guide] ({\endnew+\shift},.72) -- ({\endnew+\shift},-.15);

\draw[leftpart,brace]
  (\startnew,{\rownew+.80}) -- (\cutnew,{\rownew+.80})
  node[midway,above=6pt] {$w''_L\subseteq D$};
\draw[rightpart,brace]
  (\cutnew,{\rownew+.80}) -- (\endnew,{\rownew+.80})
  node[midway,above=6pt] {$w''_R=w'_R$};
\node[above] at (\startnew,{\rownew+1.35}) {$a$};
\node[right] at (9.05,{\rownew+.80}) {$w''=w''_Lw'_R$};
\draw[guide] (\startnew,{\rownew+.72}) -- (\startnew,{\rownew-.12});
\draw[guide] (\cutnew,{\rownew+.72}) -- (\cutnew,{\rownew-.30});
\draw[guide] (\endnew,{\rownew+.72}) -- (\endnew,{\rownew-.12});

\draw[<->,gray!80] (\endnew,-.85) -- ({\endnew+\shift},-.85)
  node[midway,above] {$k-\ell$};

\draw[brace] (\cutnew,{\rownew-.48}) -- (0,{\rownew-.48})
  node[midway,below=7pt,align=center]
  {$V_1$\\$p-(k-\ell)$ vertices in $D$};
\draw[brace] (10.8,{\rownew-.48}) -- (\cutnew,{\rownew-.48})
  node[midway,below=7pt] {$V_2$};
\end{tikzpicture}\end{lrbox}
\ifdim\wd\reviewfigurebox>\linewidth\resizebox{\linewidth}{!}{\usebox{\reviewfigurebox}}\else\usebox{\reviewfigurebox}\fi
\caption{A window $w''$ meeting both $V_1$ and $V_2$ in the proof of Fact~\ref{trunc}. The right-hand part is shared with $w'$, which starts $k-\ell$ positions later in $Q'$.}
\label{2windows}
\end{figure}

Finally, suppose that \(w''\cap V_1\neq\emptyset\) and \(w''\cap V_2\neq\emptyset\). Compare $w''$ with the
\(k\)-window $w'$ of \(Q'\) beginning \(k-\ell\) positions later. Their right-hand
parts are identical (see Figure~\ref{2windows}). The left-hand part of \(w''\) consists
entirely of vertices of \(D\), whereas the corresponding part $w'$
 ends with the \(k-\ell\)-th deleted reservoir vertex. Hence
\(w''\) contains more mandatory  vertices than $w'$, therefore at least \(k-\ell+1\).

Thus, \(Q''\) satisfies condition (i) and (ii) of Definition~\ref{dom}. Moreover, the order of vertices in $Q''$ is preserved, so are, in particular its initial
and terminal \(\ell\)-sets (as they were not removed from $Q'$ when choosing $Q''$). Thus, $Q''$ is a \((D,W^*)\)-dominant path satisfying the conclusion of Fact~\ref{trunc}. \qed

\medskip

Now, one can easily truncate $C'$ down to the length of $N$ by removing just the vertices of $W^*$, not touching the terminal $\ell$-sets of any $Q'_i$, $i=1,\dots,n$, and not  changing the order of remaining vertices. One just applies Fact~\ref{trunc} iteratively, each time to one of the $Q'_i$'s. But for this, we need to know that whenever $N'>N$ there are always $k-\ell$ reservoir vertices lying entirely in one of the sets $V(Q'_i)$.

To see that this is the case, note that at any stage of the process, by~\eqref{Cprim}, there are still at least 
$$N- \sum_{i=1}^n\alpha_i-5k$$
reservoir vertices present. In turn, by~\eqref{total} and~\eqref{moc_bi}, this is at least
$$(k-1)n-5k+\sum_{i=n+1}^{2n}b_i\ge\sum_{i=n+1}^{2n}b_i\overset{\eqref{moc_bi}}{\ge}(k-1)n>(k-\ell-1)n, $$
and, thus,  by the Pigeon-Hole Principle, at least $k-\ell$ of them belong to the same path $Q'_i$, for some $i\in[n]$.

 Since $N$ and $N'$ are divisible by $k-\ell$,  this removal operation can be continued until the total number of
vertices becomes exactly \(N\). Note, that we
shorten the paths \(Q'_i\) by removing the reservoir vertices only and in multiples of $k-\ell$, without changing the orders of their distinguished
vertices, in particular, without changing the prescribed terminal $\ell$-sets.  Hence,  the resulting cycle, call it $C''$, is still an
\((\ell,k)\)-cycle.

\medskip

It remains to replace the  reservoir vertices still present in $C''$  by the remaining vertices of $U$. To this end, we map them
bijectively onto the vertices of
\[
    \bigcup_{j=n+1}^{2n}U_j\setminus V(\Pi),
\]
 requiring additionally that no reservoir
vertex of \(Q'_i\) be mapped into \(U_{n+v_i}\) (below, we explain why this does the job).

Hall's theorem gives such a bijection for sufficiently large $n$.  Define an auxiliary bipartite graph $B=(X,Y,E)$, where $X=V(C'')\cap W^*$, $Y=\bigcup_{i=n+1}^{2n}U_i\setminus V(\Pi)$, and for $x\in X\cap V(Q_i')$ and $y\in Y$, $\{x,y\}\in E$ if $y\not\in U_{n+v_i}$. Note that $|X|=|Y|$.

As estimated above, $|X|\ge (k-1)n$. On the other hand, by the definition of $Q_i'$ and by Lemma~\ref{nu1}(b), $|V(Q_i')\cap W^*|\le |V(Q_i')|\le \nu(\alpha+\ell)\le\ell(\alpha+\ell)$, while $|U_{n+v_i}|=k-1+b_{n+v_i}=O_\alpha(1)$ by~\eqref{biO}. Thus, for $N$ (and so $n$) sufficiently large, the minimum degree $\delta(B)> |X|/2$ (with a big margin) and, quite trivially, Hall's condition holds (see, e.g.,~\eqref{su1} below with $\beta\equiv 1$). Thus, $B$ has a perfect matching, which determines the desired bijection which, in turn,  pairs the vertices to be swapped. 

As a result, we have replaced each path $Q_i'$ with a path $Q_i$. It remains to verify that $Q_i\subset H_2$.
Note that, for every edge $f$ of \(Q_i\) we have \(\min(f)=v_i\)  and $f$
contains at least \(k-\ell+1\) vertices of \(U_{v_i}\). To verify that $f\in\mathcal K$, consider a set $L\in \binom{f\cap A}\ell\cap\mathcal L$. Clearly, $\min(L)=v_i$, since $L\cap U_{v_i}\neq\emptyset$. Moreover, since $\operatorname{tr}(L)$ must be a subset of an edge of $G$ and, by the above replacement scheme, $f\cap A_{n+v_i}=\emptyset$, we have $\operatorname{tr}(L)=\{v_i\}$, that is, $L\subset A_{\min(f)}$.
Hence, $f\in H_2$.

In summary, the final \((\ell,k)\)-cycle $C$ uses all \(N\) vertices of $U$, contains edge \(e\), while all
 other edges of $C$ belong to \(H_1\cup H_2\). This completes the proof of Lemma~\ref{main2}.

\section{Proof of Lemma~\ref{88}}\label{sec:correct-loads}

In this section we prove Lemma~\ref{88}, the last missing ingredient of the proof of Lemma~\ref{main2}. Itself it relies on two further lemmas the proofs of which are presented in two subsequent subsections. But first we state these two lemmas and use them to prove Lemma~\ref{88}.

\subsection{Outline of the proof of Lemma~\ref{88}}\label{sec:correct-loads:outline}

Let \(r'_1,\ldots,r'_n\) be an admissible sequence with a
witnessing sequence $s'_1,\ldots,s'_{n-1}$ with $s'_i\in\{g,g+k-\ell\}$ for every $i$, while $p_1,\dots,p_n$ be a nonnegative integer sequence satisfying the conclusions of Fact~\ref{sumpi}.
The first lemma replaces $s'_1,\ldots,s'_{n-1}$   with only a slightly lessened sequence $s_1,\ldots,s_{n-1}$ which still satisfies~\eqref{eq:phi-two-value-witness} but with the R-H-S of~\eqref{follows} diminished by $\sum_{i=a}^bp_i$.

\begin{lem}
\label{lem:finite-shortening}
For integers $n>0$, $k>0$,  $0<\ell<k$, and  \(0\le g<k-\ell\),  let $(s'_1,\dots,s'_{n-1})$ and $(p_1,\dots,p_n)$ be two sequences of nonnegative integers such that
\[
        s'_i\in\{g,g+k-\ell\},\qquad i=1,\ldots,n-1,
\]
and
\[
        0\le p_i\le k-\ell,\quad\mbox{with}\quad|\{i: p_i>0\}|\le k.
\]

Assume further that for all \(1\le a\le b\le n\),
\begin{equation}\label{fs1}
        \sum_{i=a}^b p_i\le k-\ell+b-a
\end{equation}
and 
\begin{equation}\label{fs2}
        \sum_{i=a}^{b-1}s'_i\le (b-a+1)(k-\ell).
\end{equation}
Then, there exists a set of indices \(J\subseteq[n-1]\) of size $|J|\le 2k^2$ such that \(s'_i=g+k-\ell\) for all $i\in J$ and, putting
\begin{equation}\label{ssJ}
        s_i=
        \begin{cases}
        g, & i\in J,\\
        s'_i, & i\notin J,
        \end{cases}
\end{equation}
 for all \(1\le a\le b\le n\),
\begin{equation}\label{fs3}
        \sum_{i=a}^{b-1} s_i
        \le
        \sum_{i=a}^b(k-\ell-p_i).
\end{equation}
\end{lem}

In the forthcoming proof of  Lemma~\ref{88} we are going to represent an admissible sequence by means of a bipartite graph.
To this end, below we define a star matching, a.k.a. polygamous matching, in a bipartite graph.
A \emph{star forest} is a graph whose every connected component is a star.

\begin{df}\rm
Let \(B=(X,Y;E)\) be a bipartite graph, and let
\(\beta:Y\to \mathbb Z_{\ge 0}\).  A subgraph \(M \subseteq B\) is a \emph{$\beta$-star-matching
saturating $X$}, or just a \emph{$\beta$-matching}, if $M$ is a star forest with all star centers in $Y$ and such that  $deg_M(x)=1$ for every $x\in X$ and $deg_M(y)\le \beta(y)$ for every $y\in Y$.
\end{df}
For $T \subseteq Y$, we write
\[
\beta(T):=\sum_{y \in T} \beta(y).
\]

As an immediate consequence of Hall's Theorem (see \cite[Exercise 9, Chapter 2]{Diestel}), $B$ has a $\beta$-matching if and only if

\begin{equation}\label{su1}
        |S|\le \beta(N(S))
        \qquad\text{for every }S\subseteq X.
\end{equation}

For a sequence $\sigma=(s_1,\ldots,s_{n-1})$, define a bipartite graph $B_\sigma=(X,Y;E)$ by
\[
X=\bigcup_{i=1}^{n-1} X_i,\;\mbox{where}\;\;X_i=\{x_{i,j}: 1\le j\le s_i\},
\qquad Y=\{y_1,\ldots,y_n\},
\]
and $$E=\{\{x_{i,j},y_i\}, \{x_{i,j},y_{i+1}\}:\;\; 1\le i\le n-1,\;1\le j\le s_i\}.$$


\begin{claim}\label{nowy}
Let $B:=B_\sigma$ be as above. For nonnegative integers $q_1,\ldots,q_n$ and $\beta(y_i)=q_i$,
Hall's condition~\eqref{su1} is equivalent to

\begin{equation}\label{follows1}
\sum_{i=a}^{b-1}s_i\le\sum_{i=a}^bq_i
\qquad\mbox{for all}\quad 1\le a\le b\le n.
\end{equation}

\end{claim}
\begin{proof}
Necessity follows by applying~\eqref{su1} to
$S=\bigcup_{i=a}^{b-1}X_i$: its neighborhood is contained in
$\{y_a,\ldots,y_b\}$.

Conversely, suppose~\eqref{follows1} holds, and fix $S\subseteq X$.
If $S=\varnothing$, there is nothing to prove. Otherwise, put

\[
I:=\{i:S\cap X_i\ne\varnothing\},\qquad
\widehat S:=\bigcup_{i\in I}X_i.
\]

Then $S\subseteq\widehat S$ and $N_B(S)=N_B(\widehat S)$,
because all vertices in each $X_i$ have the same neighbors.
Write $I$ as a disjoint union of maximal sets of consecutive integers,
$[a_1,b_1],\ldots,[a_t,b_t]$, ordered so that $b_j+2\le a_{j+1}$.
The neighborhood of $\bigcup_{i=a_j}^{b_j}X_i$ is exactly
$\{y_{a_j},\ldots,y_{b_j+1}\}$, and these neighborhoods are disjoint.
Applying~\eqref{follows1} with limits $a_j,b_j+1$ gives

\[
|S|\le|\widehat S|
=\sum_{j=1}^t\sum_{i=a_j}^{b_j}s_i
\le\sum_{j=1}^t\sum_{i=a_j}^{b_j+1}q_i
=\beta(N_B(S)),
\]
yielding~\eqref{su1}.
\end{proof}

The last lemma needed for the proof of Lemma~\ref{88} establishes some kind of continuity for  $\beta$-matchings.

For two star forests $M$ and $M'$ with the same set of centers $Y$, let 
$$dist(M,M'):=\sum_{y\in Y}|deg_M(y)-deg_{M'}(y)|$$
be the $L^1$-distance of their degree vectors.

\begin{lem}\label{lem:matching-stability}
Let \(B=(X,Y;E)\) be a bipartite graph and let functions \(\beta,\beta':Y\to\mathbb Z_{\ge0}\) satisfy \(\beta(y)\le \beta'(y)\) for every \(y\in Y\). If  \(M'\) is a \(\beta'\)-matching in $B$ and condition~\eqref{su1} holds for $\beta$, then
there exists a \(\beta\)-matching \(M\)  such that

\begin{equation}\label{su3}
        dist(M,M')
        \le
        2(\beta'(Y)-\beta(Y)).
\end{equation}
\end{lem}

\medskip

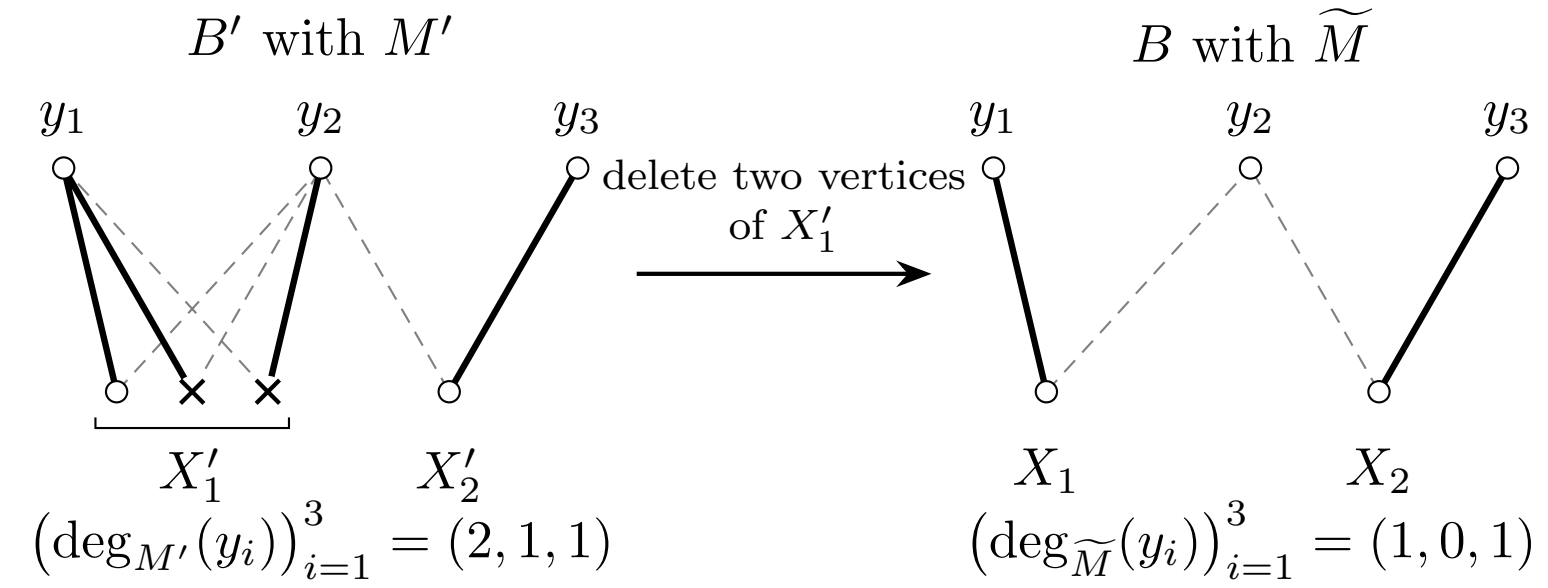
\begin{figure}[htbp]
\centering
\begin{tikzpicture}[
  x=1cm,y=1cm,>=Stealth,
  every node/.style={font=\small},
  vertex/.style={circle,draw=black,fill=white,inner sep=0pt,minimum size=4pt},
  matched edge/.style={black,very thick},
  other edge/.style={black!50,densely dashed},
  group/.style={black,thin}
]
\begin{scope}
  \node at (1.7,2.36) {$B'$ with $M'$};
  \coordinate (y1) at (0,1.48);
  \coordinate (y2) at (1.7,1.48);
  \coordinate (y3) at (3.4,1.48);
  \coordinate (x11) at (.35,0);
  \coordinate (x12) at (.85,0);
  \coordinate (x13) at (1.35,0);
  \coordinate (x21) at (2.55,0);
  \draw[other edge] (x11)--(y2) (x12)--(y2) (x13)--(y1) (x21)--(y2);
  \draw[matched edge] (x11)--(y1) (x12)--(y1) (x13)--(y2) (x21)--(y3);
  \foreach \i in {1,2,3} {
    \node[vertex,label=above:$y_{\i}$] at (y\i) {};
  }
  \node[vertex] at (x11) {};
  \node[vertex] at (x21) {};
  \foreach \v in {x12,x13} {
    \fill[white] (\v) circle (3pt);
    \draw[black,thick] ([shift={(-2.2pt,-2.2pt)}]\v)--([shift={(2.2pt,2.2pt)}]\v)
                     ([shift={(-2.2pt,2.2pt)}]\v)--([shift={(2.2pt,-2.2pt)}]\v);
  }
  \draw[group] (.21,-.17)--(.21,-.24)--(1.49,-.24)--(1.49,-.17);
  \node[below] at (.85,-.24) {$X'_1$};
  \node[below] at (2.55,-.24) {$X'_2$};
  \node at (1.7,-.98) {$\bigl(\deg_{M'}(y_i)\bigr)_{i=1}^3=(2,1,1)$};
\end{scope}
\draw[->,thick] (3.79,.78)--(5.74,.78)
  node[midway,above,align=center,font=\scriptsize] {delete two vertices\\of $X'_1$};
\begin{scope}[xshift=6.15cm]
  \node at (1.7,2.36) {$B$ with $\widetilde M$};
  \coordinate (y1) at (0,1.48);
  \coordinate (y2) at (1.7,1.48);
  \coordinate (y3) at (3.4,1.48);
  \coordinate (x11) at (.35,0);
  \coordinate (x21) at (2.55,0);
  \draw[other edge] (x11)--(y2) (x21)--(y2);
  \draw[matched edge] (x11)--(y1) (x21)--(y3);
  \foreach \i in {1,2,3} {
    \node[vertex,label=above:$y_{\i}$] at (y\i) {};
  }
  \node[vertex] at (x11) {};
  \node[vertex] at (x21) {};
  \node[below] at (.35,-.24) {$X_1$};
  \node[below] at (2.55,-.24) {$X_2$};
  \node at (1.7,-.98) {$\bigl(\deg_{\widetilde M}(y_i)\bigr)_{i=1}^3=(1,0,1)$};
\end{scope}
\end{tikzpicture}
\caption{Restricting $M'$ after deleting two vertices of $X'_1$ ($k-\ell=2$, $g=1$).
Crosses mark deleted vertices. Solid edges belong to the matchings $M'$ and $\widetilde M$ in the left and right panels, respectively. Dashed lines are the remaining edges.}
\label{fig:shorten-matching}
\end{figure}

\begin{proof}[Proof of Lemma~\ref{88}]
Given an admissible sequence $(r'_i)$ with witness $\sigma'=(s'_i)$, create a corresponding bipartite graph $B':=B_{\sigma'}=(X',Y',E')$ as described above. According to Claim~\ref{nowy} with $\beta(y_i)=r'_i$, in view of~\eqref{follows}, $B'$ satisfies Hall's condition~\eqref{su1}  and thus contains a matching $M'$ with $deg_{M'}(y_i)=r_i'\le k-\ell$, $i=1,\dots,n$. In particular, $M'$ is a $\beta'$-matching, where  $ \beta'\equiv k-\ell$.
 
 Note that, by~\eqref{follows}, sequence $(s_i')$ satisfies the assumptions of Lemma~\ref{lem:finite-shortening}. Further, let $(p_i)$ be a sequence satisfying the conclusions of Fact~\ref{sumpi} and, thus, the assumptions of Lemma~\ref{lem:finite-shortening}. 
 By applying that lemma
 we obtain a
    sequence $(s_i)$ with $s_i\le s_i'$ for all $i\in[n-1]$ satisfying~\eqref{ssJ} and~\eqref{fs3}. 
    
    We modify $B'$ by removing $s_i'-s_i$ vertices from each set $X_i'$, $i=1,\dots, n-1$, obtaining a new bipartite graph $B$ and a new $\beta'$-matching $\tilde M$ (see Figure~\ref{fig:shorten-matching}). Note that, by~\eqref{ssJ}, only at most $2k^2$ sets $X'_i$ have been diminished and if so, then by exactly $k-\ell$. Thus,
   $$dist(\tilde M,M')\le 2k^2(k-\ell).$$
    Note also that by~\eqref{fs3} and Claim~\ref{nowy} with $q_i=\beta(y_i):=k-\ell-p_i$, Hall's condition~\eqref{su1} holds.
   
   Finally, apply Lemma~\ref{lem:matching-stability} to $B$, $\beta'$,  $\tilde M$, and $\beta$, obtaining a $\beta$-matching $M$ such that $$dist(M,\tilde M)\le 2(\beta'(Y)-\beta(Y))=2\sum_ip_i\le 2k(k-\ell).$$
   Set $r_i=deg_M(y_i)$, $i=1,\dots,n$. Then
   $$\sum_i|r_i-r'_i|=dist(M,M')\le dist(\tilde M,M')+dist(M,\tilde M)\le2k(k-\ell)(k+1)\le 2k^3,$$
   so, condition~\eqref{eq:corrected-stability} holds. Condition~\eqref{eq:corrected-loads} holds by the definition of a $\beta$-matching.
  It remains to 
   establish the admissibility of $(r_i)$ with witness $(s_i)$. By~\eqref{eq:corrected-loads}, $r_i\le k-\ell$, while by~\eqref{ssJ}, $s_i\in\{g,g+k-\ell\}$.  
   As for~\eqref{adm}, 
   note that every vertex of $\bigcup_{i=1}^{t-1}X_i$ must be matched to one of $y_1,\dots, y_t$ only, while, in turn, these vertices of $Y$  can be adjacent only to the vertices in $\bigcup_{i=1}^{t}X_i$. These observations yield  the two inequalities in~\eqref{adm}. The equality of the total sums follows by
the saturation of $X$ by $M$.
    \end{proof}

\subsection{Proof of Lemma~\ref{lem:finite-shortening}}\label{sec:correct-loads:shortening}
Let $(s'_1,\dots,s'_{n-1})$ and $(p_1,\dots,p_n)$ satisfy the assumptions of Lemma~\ref{lem:finite-shortening}.
Set 
$$I:=\{i:p_i>0\}\quad\mbox{and}\quad\Lambda=\{j:s'_j=g+k-\ell\}.$$ For every
\(i\in I\), let $\lambda_i=|\{j\in \Lambda:\; j<i\}|$ and  $\rho_i=|\{j\in \Lambda:\; j\ge i\}|$, and let $J_i$ be the set of
$\min\{\lambda_i,k\}$  indices $j\in\Lambda,\;j<i$,  nearest to $i$ 
and $\min\{\rho_i,k\}$  indices $j\in\Lambda,\;j\ge i$,
 nearest to $i$.
Finally, let \(J=\bigcup_{i\in I}J_i\). Then
\[
        |J|\le 2k|I|\le 2k^2.
\]
\begin{example}\rm
Let $n=100$, $k=3$, $I=\{8,14,65\}$, and 
\[
\begin{aligned}
\Lambda=\{&3,7,8,11,13,14,17,22,30,31,33,41,\\
          &42,44,47,49,52,54,61,66,68,70,82,99\}.
\end{aligned}
\] Then, $\lambda_8=2$, while $\min\{\rho_8,\lambda_{14},\rho_{14},\lambda_{65},\rho_{65}\}\ge3$. Thus, $$J_8=\{3,7,8,11,13\},\qquad J_{14}=\{8,11,13,14,17,22\},$$ and 
$$J_{65}=\{52,54,61,66,68,70\},$$
so, 
$$J=\{3,7,8,11,13,14,17,22,52,54,61,66,68,70\}.$$

\end{example}

We are going to prove that this choice of \(J\) yields \eqref{fs3} (with $s_i$'s defined by~\eqref{ssJ}).
Fix $1\le a\le b\le n$, and put
\[
        L:=b-a+1,
        \quad\mbox{and}\quad
        P:=\sum_{i=a}^b p_i.
\]
Further, let
\[
        Y:=|\Lambda\cap[a,b-1]|,
        \quad\mbox{and}\quad
        Y_J:=|J\cap[a,b-1]|,
\]
and observe that, by the definition of the set $J$, if $I\cap[a,b]\neq\emptyset$, then
 \(Y_J\ge \min\{Y,k\}\).  \space
Using the above notation, by \eqref{fs2},
\begin{equation}\label{s'bound}
\sum_{i=a}^{b-1}s'_i=(L-1)g+(k-\ell)Y\le L(k-\ell),
\end{equation}
while
\[
        \sum_{i=a}^{b-1}s_i=(L-1)g+(Y-Y_J)(k-\ell).
\]
Thus, the desired inequality~\eqref{fs3} is equivalent to

       $$ (L-1)g+(Y-Y_J)(k-\ell)\le L(k-\ell)-P,$$
or
\begin{equation}\label{fs4}D:=(L-1)g+Y(k-\ell)-L(k-\ell)+P\le Y_J(k-\ell).\end{equation}

By the right-hand inequality in~\eqref{s'bound} and by~\eqref{fs4}, $D\le P$. If $D\le0$,
then~\eqref{fs4} holds. Assume thus that $D\ge 1$. 
Since $p_i\le k-\ell$ for all $i$,
 $$1\le h:=\left\lceil \frac{D}{k-\ell}\right\rceil\le\left\lceil\frac{P}{k-\ell}\right\rceil=\left\lceil\frac{\sum_{i=a}^bp_i}{k-\ell}\right\rceil\le |I\cap[a,b]|\le k.$$
 It follows that, indeed,  $I\cap[a,b]\neq\emptyset$, and by an earlier observation, \(Y_J\ge \min\{Y,k\}\).
Moreover, by \eqref{fs1} which can be now restated as $P\le k-\ell+L-1$,
\[
\begin{aligned}
        D
        &=(L-1)g+Y(k-\ell)-L(k-\ell)+P \\
        &\le (L-1)g+Y(k-\ell)-L(k-\ell)+k-\ell+L-1 \\
        &=Y(k-\ell)+(L-1)(g-k+\ell+1)\le Y(k-\ell),
\end{aligned}
\]
because \(g\le k-\ell-1\) (cf.~\eqref{gap}). Hence,  by the definition of $h$, we have $ h\le Y$, and
thus, because also $h\le k$, \(Y_J\ge\min\{Y,k\}\ge h\). Consequently, again by the definition of $h$,  \(D\le Y_J(k-\ell)\), which is  \eqref{fs4}.

\subsection{Proof of Lemma~\ref{lem:matching-stability}}\label{sec:correct-loads:stability}

Let $B,\beta,\beta'$ and $M'$ be as in the assumptions of Lemma~\ref{lem:matching-stability}.
If the $\beta'$-matching \(M'\) satisfies $deg_{M'}(y)\le \beta(y)$ for all $y\in Y$, then it is  also a $\beta$-matching and one may take $M=M'$. Otherwise, there is a vertex \(y_0\in Y\)
such that $deg_{M'}(y_0)>\beta(y_0)$.
We call such a vertex \emph{overfull}. Any vertex \(y\in Y\) with  $deg_{M'}(y)<\beta(y)$ is called \emph{underfull}. The \emph{overload} of $M'$ (with respect to $\beta$) is defined as
$$over(M')=\sum_{y\in Y}\max\{0,deg_{M'}(y)-\beta(y)\}$$
(so, $over(M')=0$ means that $M'$ is a $\beta$-matching).

We are going to show that for every overfull vertex $y_0$ there is an
alternating path from $y_0$ to an underfull vertex $y_t$ 
\[
     P=   y_0,x_1,y_1,x_2,y_2,\ldots,x_t,y_t,
\]
where for each $i=1,\dots,t$ we have \(\{y_{i-1},x_i\}\in M'\) and  \(\{x_i,y_{i}\}\in E(B)\).
Given such a path, one can alter the
$\beta'$-matching $M'$ along this path by taking the symmetric difference $E(M')\triangle E(P)$, that is,  each \(x_i\) is reassigned from \(y_{i-1}\) to~\(y_i\).

Let us denote the obtained subgraph by $M^*$.
Since this operation keeps all vertices of \(X\) matched, $M^*$ is still a $\beta'$-matching with
\[
\begin{gathered}
deg_{M^*}(y_0)=deg_{M'}(y_0)-1,\qquad deg_{M^*}(y_t)=deg_{M'}(y_t)+1,\\
deg_{M^*}(y)=deg_{M'}(y)\qquad\text{for all }y\notin\{y_0,y_t\}.
\end{gathered}
\]
Hence, $dist(M',M^*)=2$, while  $over(M^*)=over(M')-1$,  because $y_0$ is overfull and  $y_t$ is underfull.
After iterating this operation  $over (M')$ times one would arrive at a $\beta$-matching $M$ with
$$dist(M,M')\le2\cdot over(M')\le 2(\beta'(Y)-\beta(Y)),$$
where the last inequality follows from $deg_{M'}(y)\le \beta'(y)$ for all $y\in Y$.

It remains to prove that such a path $P$ always exists. 
Consider all alternating paths $P$ described above, but without the assumption that $y_t$ is underfull. Let \(Y_0\subseteq Y\) be the set of all their end-vertices, including $y_0$ as the path of length zero.
We claim that at least one vertex in $Y_0$ is underfull. Indeed, suppose otherwise
and  let
\[X_0=\{x\in X: \{x,y\}\in E(M')\quad\mbox{for some}\quad y\in Y_0\}.\]
From this definition it follows that
\[
        |X_0|=\sum_{y\in Y_0}deg_{M'}(y)
\]
and
\begin{equation}\label{su22}
        N(X_0)\subseteq Y_0.
\end{equation}
To see the latter, observe that whenever $\{x,y\}\in M'$ and $y\in Y_0$, then for every $y'\in Y$ such that $\{x,y'\}\in E(B)$, we have $y'\in Y_0$ too.

Since \(y_0\in Y_0\) is overfull and, by assumption, no vertex of \(Y_0\) is
underfull, we have
\[
        deg_{M'}(y)\ge \beta(y)\quad\text{for all }y\in Y_0,
\]
and the inequality is strict for \(y_0\). Hence
\begin{equation}\label{su21}
        |X_0|
        =
        \sum_{y\in Y_0}deg_{M'}(y)
        >
        \sum_{y\in Y_0}\beta(y)
        =
        \beta(Y_0).
\end{equation}

Combining \eqref{su21} and \eqref{su22}, we obtain
\[
        \beta(N(X_0))\le \beta(Y_0)<|X_0|,
\]
which violates Hall's condition~\eqref{su1} for $\beta$ (and $S=X_0$) and, thus, we arrive at a contradiction with the assumptions of the lemma.

\paragraph{Use of generative AI.}
The authors used OpenAI's GPT-5.6 and Astra-6.0 primarily in the later stages of
the work to help check proofs, prepare figures, and edit the exposition.
GPT-5.6 also proposed draft proofs of Lemmas~\ref{lem:finite-shortening} and~\ref{lem:matching-stability} and draft
arguments for the Appendix. The authors developed the overall proof
strategy and the ideas used to resolve the main technical difficulties.
They reviewed the proposed arguments, rewrote all AI-assisted material,
and take full responsibility for the final manuscript.

\section*{Appendix}

Recall from Section~\ref{eHam} that a \emph{$k$-window} is a block of
$k$ consecutive entries of a sequence. It is \emph{relevant} if its
initial position is congruent to $1$ modulo $k-\ell$.
We use this terminology for both vertex sequences and binary sequences used in Appendix A1.

\subsection*{A1. Proofs of Observation~\ref{obs1} and Lemma~\ref{nu1}}\label{obs_lem2_3}
Given a $U$-dominant $(\ell, k)$-path $P$ with $V(P)=(u_1,\ldots,u_t)$, we encode it by
a binary sequence $\hat{P}=\hat{u}_1\cdots \hat{u}_t$, where $\hat{u}_i=1$ if $u_i\in U$
and $\hat{u}_i=0$ otherwise. Clearly, the edges of $P$
correspond exactly to the relevant $k$-windows of $\hat{P}$.
We write $|\hat{P}|$ for the length of $\hat{P}$ and denote by $0^d$ and $1^d$ the sequences consisting
of $d$ zeros and $d$ ones, respectively.

With an infinite reservoir $W^*$, a binary sequence of length $t\ge k$ with $t\equiv\ell\pmod{k-\ell}$
represents a $U$-dominant path if it contains exactly
$|U|$ ones, begins and ends with $1^\ell$, and every relevant
$k$-window contains at least $k-\ell+1$ ones. Such a path is obtained
by assigning distinct vertices of $U$ to the ones and distinct
vertices of $W^*$ to the zeros.

\begin{proof}[Proof of Observation~\ref{obs1}]
Let $P$ be an $(\ell,k)$-path satisfying assumption (ii) in Definition~\ref{dom} for a given set $U$. Set $t=|V(P)|$ and $q=|U\cap V(P)|\le|U|$.
Let the binary sequence $\hat{P}=\hat{u}_1\cdots \hat{u}_t$ encode $P$, and let
$j_1<\cdots<j_q$ be the positions occupied by ones. 
Note that every block of $\ell$ consecutive entries of $\hat P$ is contained in a relevant
$k$-window and, on the other hand,  by  assumption (ii), each such window
contains at most $\ell-1$ zeros. Hence,  $\hat{P}$ contains no $0^\ell$ as a block. This implies that
\[
 j_{i+1}-j_i\le\ell\quad\mbox{for all} \quad1\le i<q\quad\mbox{and}
\quad t-j_q\le\ell-1.\]
Moreover, the first relevant $k$-window contains at least
$k-\ell+1$ ones which implies that $q\ge k-\ell+1$ and
$j_{k-\ell+1}\le k$. Using these inequalities we conclude that
\begin{align*}
 t &= j_{k-\ell+1}
      +\sum_{i=k-\ell+1}^{q-1}(j_{i+1}-j_i)
      +(t-j_q)\\
 &\le k+\bigl(q-(k-\ell+1)\bigr)\ell+\ell-1\\
 &=\ell q-(k-\ell-1)(\ell-1)
 \le\ell q\le\ell|U|.
\end{align*}
\end{proof}

\begin{proof}[Proof of Lemma~\ref{nu1}] For parts (a)-(e), let $\alpha\ge 2\ell$ be an integer.

\smallskip
\noindent\emph{(a)} Choose $q\in\{1,\ldots,k-\ell\}$ such that
$\alpha+q\equiv\ell\pmod{k-\ell}$, and consider the binary sequence
\[
 \beta=1^\ell 0^q1^{\alpha-\ell}.
\]
Since $\alpha\ge2\ell$, the first and last $\ell$ entries of $\beta$ are ones.
Moreover, by the assumption on $\ell$, every $k$-window of $\beta$ contains at least
$k-q\ge\ell\ge k-\ell+1$ ones. Hence, $\beta$ represents a $U$-dominant path
with  $\alpha+q$  vertices, and, by the definition of $\nu$, 
\[
 \nu(\alpha)\ge\alpha+q\ge\alpha+1.
\]

\smallskip
\noindent\emph{(b)} Apply Observation~\ref{obs1} to a path attaining
$\nu(\alpha)$.

\smallskip
\noindent\emph{(c)} Let $P$ attain $\nu(\alpha)$ with a mandatory set $U$. By (a), there is a vertex $u\in V(P)\setminus U$.
Then, $P$ is also $U\cup\{u\}$-dominant, and so
$\nu(\alpha+1)\ge|V(P)|=\nu(\alpha)$.

\smallskip
\noindent\emph{(d)} Let $\hat P$ encode a path $P$ attaining the maximum in $\nu(\alpha)$.
The sequence $1^{k-\ell}\hat{P}$ has $\alpha+k-\ell$ ones.
Its first relevant $k$-window is $1^k$, while all subsequent
relevant $k$-windows are precisely those of $\hat{P}$. Moreover, its first
 $\ell$ entries, as well as its last $\ell$ entries are all ones. Consequently, $1^{k-\ell}\hat P$ represents a $U$-dominant path with $|U|=\alpha+k-\ell$. Hence,
\[
 \nu(\alpha+k-\ell)\ge|1^{k-\ell}\hat{P}|
 =\nu(\alpha)+k-\ell,
\]
which was to be proved.

\smallskip
\noindent\emph{(e)} We will show a slightly stronger bound with $k-1$ replaced by $\ell+g$, where $g$ is the gap (cf.~\eqref{gap}).
Let $\hat{P}=\hat{u}_1\cdots \hat{u}_t$ encode a path $P$ attaining the maximum in $\nu(\alpha+1)=:t$. In particular, $\hat P$ contains exactly $\alpha+1$ ones.
If $t-(\ell+g)\le\alpha$, the bound follows from (a).
Assume therefore that $t-(\ell+g)>\alpha$ and recall that $\hat P$ begins with $w_0:=1^\ell$.

Since $\ell+g$ is divisible by $k-\ell$, the $k$-window $w$  ending
at position $2\ell+g$ is relevant and $|w\cap w_0|=k-\ell-g\le k-\ell$. Consequently,
the block $w_1:= \hat{u}_{\ell+1}\cdots \hat{u}_{2\ell+g}$ contains at least one digit one.
Thus, the  subsequence
\[
 \hat{ P'}=\hat P-w_1=1^\ell \hat{u}_{2\ell+g+1}\cdots \hat{u}_t
\]
obtained from $\hat P$ by  removing  block $w_1$, has length $|\hat{P}'|=t-(\ell+g)>\alpha$ and contains at most
$\alpha$ ones.

Observe that $\hat{P}'$ may also be obtained from the block
$w_2:=\hat{u}_{\ell+g+1}\cdots \hat{u}_t$ by replacing its first $\ell$ entries
by ones (see Figure~\ref{fig:nu-delete-block}).

\begin{figure}[htbp]
\centering
\begin{tikzpicture}[
  x=0.97cm,y=1cm,
  every node/.style={font=\small,inner sep=2pt},
  brace/.style={decorate,decoration={brace,amplitude=3.5pt}},
  edge/.style={black!75,thin}
]
  \foreach \y/\a in {0/0,-1.75/4.25}{
    \fill[black!7] (7.05,{\y-.26}) rectangle (10.40,{\y+.26});
    \draw[edge] (\a,{\y-.26}) rectangle (10.40,{\y+.26});
    \draw[edge] (7.05,{\y-.26}) -- (7.05,{\y+.26});
    \node at (8.725,\y) {$\hat u_{2\ell+g+1}\cdots\hat u_t$};
  }

  \draw[edge] (2.80,-.26) -- (2.80,.26);
  \draw[edge] (4.25,-.26) -- (4.25,.26);
  \node[left=4pt] at (0,0) {$\widehat P$};
  \node at (1.40,0) {$1^\ell$};
  \node[font=\footnotesize,inner sep=0pt] at (3.525,0) {$g$ entries};
  \node at (5.65,0) {$\ell$ entries};

  \draw[brace] (0,.37) -- (2.80,.37)
    node[midway,above=5pt] {$w_0$};
  \draw[brace] (2.80,.37) -- (7.05,.37)
    node[midway,above=5pt] {$w_1$ (contains a $1$)};
  \draw[brace] (7.05,-.36) -- (2.04,-.36)
    node[midway,below=5pt] {$w$ \,($|w\cap w_0|=k-\ell-g$)};

  \node[left=4pt] at (4.25,-1.75) {$\widehat{P}'$};
  \node at (5.65,-1.75) {$1^\ell$};
\end{tikzpicture}
\caption{$\widehat{P}'=\hat{P}-w_1$ in the proof of Lemma~\ref{nu1}(e).}
\label{fig:nu-delete-block}
\end{figure}
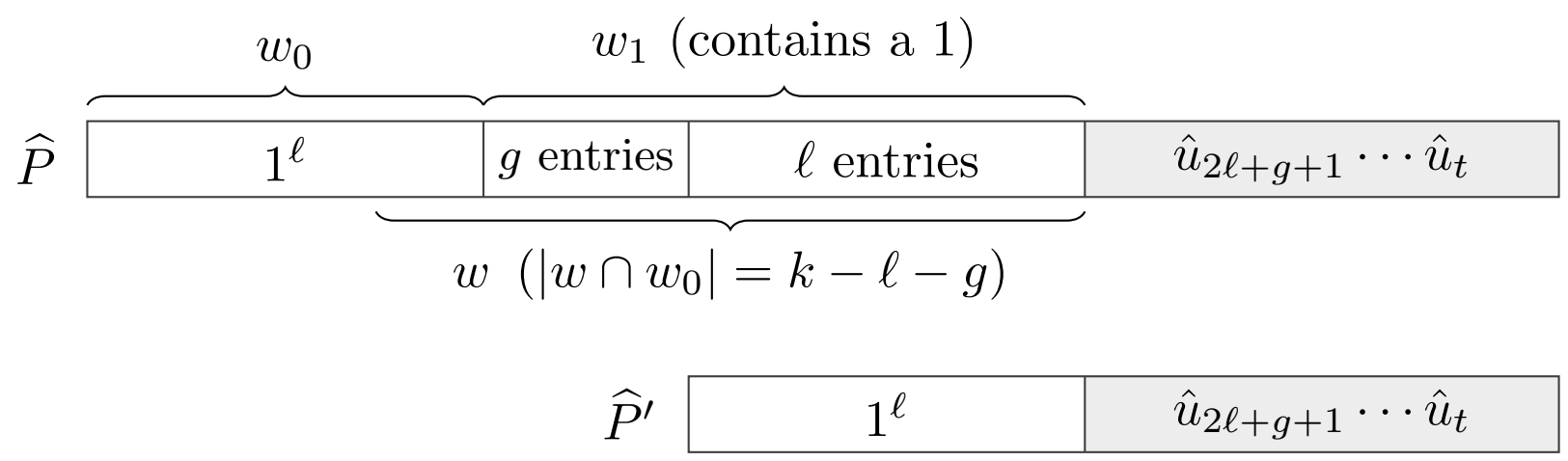

Since $\ell+g$ is divisible by $k-\ell$, every relevant $k$-window of $w_2$ is a relevant window in $\hat P$, and thus, 
 contains at least $k-\ell+1$ ones. This property cannot be destroyed by replacing zeros by ones,  so every relevant window of $\hat P'$ contains at least $k-\ell+1$ ones too. If $\hat P'$ contains fewer than $\alpha$ ones, we replace some  zeros by ones to get a binary  sequence $\hat P''$   with exactly $\alpha$ ones which 
represents a $U$-dominant path with $|U|=\alpha$. Hence
\[
 \nu(\alpha)\ge |\hat{P}''|
 =\nu(\alpha+1)-(\ell+g),
\]
which was to be proved.

\smallskip
\noindent\emph{(f)} Recall from~\eqref{gap} that $k+g$ is divisible
by $k-\ell$. For each integer $t\ge1$, consider  a binary sequence
$$\beta_t=1^\ell\bigl(0\,1^{k+g-1}\bigr)^t.$$
It has exactly $ \alpha_t=\ell+t(k+g-1)$ one and $t$ zeros, and, thus,  the length
$\ell+t(k+g)\equiv\ell\pmod{k-\ell}$. Moreover, its first
 $\ell$ entries, as well as its last $\ell$ entries are all ones. And most importantly, any consecutive zeros are $k+g\ge k$ positions apart,
so every $k$-window contains at most one zero and hence at least
$k-1\ge k-\ell+1$ ones. 
Consequently, $\beta_t$ represents a $U$-dominant path with $|U|=\alpha_t$
and so
\[
 \nu(\alpha_t)\ge|\beta_t|=\alpha_t+t.
\]
For every sufficiently large $\alpha$, let $t:=t_\alpha$ be the unique integer $t$ such that
$\alpha_t\le\alpha<\alpha_{t+1}$. Then, by (c),
\[
 \nu(\alpha)-\alpha
 \ge\nu(\alpha_t)-\alpha
 \ge t-(\alpha-\alpha_t)
 \ge t-(k+g-1).
\]
Since $t\to\infty$ as $\alpha\to\infty$, this proves (f).
\end{proof}

\subsection*{A2. Proofs of Claim~\ref{just2} and Lemmas~\ref{phi-step} and~\ref{phi-glue}}

\begin{proof}[Proof of Claim~\ref{just2}]
By~\eqref{follows}, for all $j$, we have $s_j\le r_j+r_{j+1}\le2(k-\ell)$,
so the additional value $g+2(k-\ell)$ is impossible if $g>0$.
If $g=0$ and $s_j=2(k-\ell)$, then $r_j=r_{j+1}=k-\ell$ and
\[
\sum_{i=1}^{j-1}r_i=\sum_{i=1}^{j-1}s_i.
\]
Indeed, by~\eqref{adm} and $r_j,r_{j+1}\le k-\ell$,
\[
\begin{aligned}
\sum_{i=1}^{j-1}s_i+2(k-\ell)
&=\sum_{i=1}^{j}s_i
 \le\sum_{i=1}^{j+1}r_i\\
&=\sum_{i=1}^{j-1}r_i+r_j+r_{j+1}
 \le\sum_{i=1}^{j-1}s_i+2(k-\ell).
\end{aligned}
\]
Thus,  equality holds throughout, proving both assertions.
Replace $s_j$ by $k-\ell$ and $r_j$ by $0$.
Conditions in~\eqref{adm} remain valid: for $t<j$ their sums
are unchanged while for $t=j$ the new middle sum equals the left-hand
sum. Furthermore, for $t>j$ all three sums decrease by $k-\ell$.
Both total sums also decrease by $k-\ell$.
Since  $\alpha-k+\ell\ge 2\ell$, by Lemma~\ref{nu1}(d) with $\alpha:=\alpha-k+\ell$, 
\[
\nu(\alpha-k+\ell)+k-\ell\le\nu(\alpha),
\]
so this replacement does not decrease the expression being maximized.
Repeating this operation if necessary, proves the assertion. \end{proof}

\begin{proof}[Proof of Lemma~\ref{phi-step}]
Let $\alpha\ge k+\ell$ and $m\ge1$.
Since~\eqref{phi-bounds} follows 
by summing
\eqref{phi-dif-bounds} sidewise (for $m'=1,\dots,m$), it suffices to prove~\eqref{phi-dif-bounds} only.
The case $m=1$ holds, because $\Phi(\alpha,0)=0$ and
$\Phi(\alpha,1)=\nu(\alpha)$. Assume $m\ge2$.

For the lower bound, choose a sequence
$(r'_1,\ldots,r'_{m-1})\in\mathcal A_{m-1}$ attaining the maximum in
$\Phi(\alpha,m-1)$, with witness $(s'_1,\ldots,s'_{m-2})$
(when $m=2$, the latter sequence  empty). Put
\[
 (r_1,\ldots,r_m)=(g,r'_1,\ldots,r'_{m-1}),
 \qquad
 (s_1,\ldots,s_{m-1})=(g,s'_1,\ldots,s'_{m-2}).
\]
These sequences satisfy~\eqref{adm}. Indeed, the first pair of inequalities 
becomes $0\le g\le g$, while the remaining inequalities and the equality follow from the
original ones by adding $g$ sidewise. Therefore,
 $(r_1,\ldots,r_m)\in\mathcal A_m$ and thus, by the definition of function $\Phi$ and by Lemma~\ref{nu1}(e),
\begin{align*}
 \Phi(\alpha,m)
 &\ge\sum_{i=1}^m(\nu(\alpha-r_i)+r_i)\\
 &=\Phi(\alpha,m-1)+\nu(\alpha-g)+g\\
 &\ge\Phi(\alpha,m-1)+\nu(\alpha)-g(k-2).
\end{align*}

For the upper bound, choose $(r_1,\ldots,r_m)$ attaining the maximum in
$\Phi(\alpha,m)$, with a witness $(s_1,\allowbreak\ldots,s_{m-1})$
satisfying~\eqref{eq:phi-two-value-witness}. Put
$a=s_{m-1}-r_m$. Then 
\begin{align*}
a &= \left(\sum_{i=1}^{m-1}s_i -\sum_{i=1}^{m-2}s_i\right) 
  - \left(\sum_{i=1}^{m}r_i-\sum_{i=1}^{m-1}r_i\right)\\
 &=\sum_{i=1}^{m-1}r_i-\sum_{i=1}^{m-2}s_i +\left(\sum_{i=1}^{m}r_i-\sum_{i=1}^{m-1}s_i\right)
\end{align*}
and so, by~\eqref{adm}, 
\[
0\leq a =\sum_{i=1}^{m-1}r_i-\sum_{i=1}^{m-2}s_i\leq r_{m-1}.
\]

Set $\rho=(r_1,\ldots,r_{m-2},r_{m-1}-a)$. We are going to show that $\rho\in\mathcal A_{m-1}$ with the truncated witness $\sigma=(s_1,\dots,s_{m-2})$.
Since $0\le a\le r_{m-1}$, we have
\[
0\le r_{m-1}-a\le k-\ell,
\]
and it remains to prove that~\eqref{adm} holds for $\rho$ and $\sigma$.

Since for $m=2$, $\rho=(0)\in\mathcal A_1$,
we assume that $m\ge3$.
The inequalities in~\eqref{adm} hold for $\rho$ and $\sigma$, because they are exactly the same inequalities which hold, by assumption, for $(r_i)_1^m$ and $(s_i)_1^{m-1}$ for $t=1,\dots,m-2$.

 Moreover, since $a+r_m=s_{m-1}$,
\[
\sum_{i=1}^{m-2}r_i+(r_{m-1}-a)
=\sum_{i=1}^{m}r_i-s_{m-1}
=\sum_{i=1}^{m-2}s_i.
\]

It follows that
$(r_1,\ldots,r_{m-2},r_{m-1}-a)\in\mathcal A_{m-1}$,
with witness $(s_1,\ldots,s_{m-2})$, and, by the definition of $\Phi$, 
\[
\Phi(\alpha,m-1)
\ge
\sum_{i=1}^{m-2}\bigl(\nu(\alpha-r_i)+r_i\bigr)
+\nu(\alpha-r_{m-1}+a)+r_{m-1}-a.
\]
Consequently, using the extremality of $(r_1,\ldots,r_m)$ together with $a\ge0$ and monotonicity of $\nu$ (cf. Lemma~\ref{nu1}(c)),
we obtain
\begin{align*}
\Phi(\alpha,m-1)
&\ge
\sum_{i=1}^{m-1}\bigl(\nu(\alpha-r_i)+r_i\bigr)-a\\
&=\Phi(\alpha,m)-\nu(\alpha-r_m)-r_m-a.
\end{align*}
Rearranging and recalling that $a+r_m=s_{m-1}$, we conclude
that
\begin{align*}
\Phi(\alpha,m)
&\le \Phi(\alpha,m-1)+\nu(\alpha-r_m)+r_m+a\\
&=\Phi(\alpha,m-1)+\nu(\alpha-r_m)+s_{m-1}\\
&\le \Phi(\alpha,m-1)+\nu(\alpha)+g+k-\ell,
\end{align*}
where the last inequality follows from $r_m\ge0$,
monotonicity of $\nu$, and
$s_{m-1}\in\{g,g+k-\ell\}$.

Thus,~\eqref{phi-dif-bounds} holds with
\[
 c=\max\{g(k-2),\,g+k-\ell\}.
\]
\end{proof}

\begin{proof}[Proof of Lemma~\ref{phi-glue}]

We only show the first inequality, since the second one then follows by induction on $q$. Let $\alpha\ge k+\ell$ and $a,b\ge1$ be all integers.
Choose sequences
\[
 (r'_1,\ldots,r'_a)\in\mathcal A_a,
 \qquad
 (r'_{a+1},\ldots,r'_{a+b})\in\mathcal A_b
\]
attaining the maxima in, respectively, $\Phi(\alpha,a)$ and $\Phi(\alpha,b)$ with   witnesses
\[
 (s'_1,\ldots,s'_{a-1}),
 \qquad
 (s'_{a+1},\ldots,s'_{a+b-1})
\]
whose entries belong to $\{g,g+k-\ell\}$ (cf. Claim~\ref{just2}). Recall that for $a=1$ or $b=1$ the witnesses are empty.

Based on these sequences, our goal is to construct  a sequence (and its witness) belonging to $\mathcal A_{a+b}$. To this end, apply Lemma~\ref{88} to sequence $(r_1',\dots,r_a')$ with
\[
 p_1=p_2=\cdots =p_{a-1}=0\quad\mbox{and}\quad p_a=g.
\]
Its assumptions, that is, inequalities~\eqref{eq:p-interval-corrected} and~\eqref{npk}, hold, since at most one $p_i$ is positive and, by \eqref{gap}, for every
$I\subseteq[a]$,
\[
 \sum_{i\in I}p_i\le g\le k-\ell-1+|I|.
\]
As an output, we obtain an admissible sequence
$(\widetilde r_1,\ldots,\widetilde r_a)$ with witness
$(\widetilde s_1,\ldots,\widetilde s_{a-1})$ such that
\begin{align*}
& \widetilde r_a\le k-\ell-g,\\
& \widetilde r_i\le k-\ell, \qquad\qquad\;\; i=1,\dots,a-1,\\
&   \widetilde s_i\in\{g,g+k-\ell\},\quad i=1,\dots,a-1,\\
& \sum_{i=1}^a|\widetilde r_i-r'_i|\le 2k^3.
\end{align*}

Define sequences $(r_i)_{1}^{a+b}$ and $(s_i)_{1}^{a+b-1}$  by
\[
 r_i=
 \begin{cases}
  \widetilde r_i,&1\le i<a,\\
  \widetilde r_a+g,&i=a,\\
  r'_i,&a<i\le a+b,
 \end{cases}
 \qquad
 s_i=
 \begin{cases}
  \widetilde s_i,&1\le i<a,\\
  g,&i=a,\\
  s'_i,&a<i<a+b.
 \end{cases}
\]
We now verify that $(r_i)\in\mathcal A_{a+b}$ with $(s_i)$ as a witness.
The inequalities $0\le r_i\le k-\ell$ follow from the construction.
For $t<a$, conditions~\eqref{adm} are exactly those for
$(\widetilde r_i)$. 
At $t=a$, the required inequalities follow from
\[
\sum_{i=1}^{a-1}s_i
=\sum_{i=1}^{a-1}\widetilde s_i
=\sum_{i=1}^{a}\widetilde r_i
\le \sum_{i=1}^{a}\widetilde r_i+g
=\sum_{i=1}^{a}r_i
=\sum_{i=1}^{a}s_i.
\]
For $t>a$, subtracting $\sum_{i=1}^{a}r_i=\sum_{i=1}^{a}s_i$
from all three sums in~\eqref{adm} gives
\[
\sum_{i=a+1}^{t-1}s'_i
\le \sum_{i=a+1}^{t}r'_i
\le \sum_{i=a+1}^{t}s'_i,
\]
which are precisely the admissibility inequalities
satisfied by  $(r'_{a+1},\ldots,r'_{a+b})$.
Thus $(r_i)_{i=1}^{a+b}\in\mathcal A_{a+b}$. 

Moreover,
\begin{equation}\label{rrprim}
 \sum_{i=1}^{a+b}|r_i-r'_i|
 \le\sum_{i=1}^a|\widetilde r_i-r'_i|+g
 \le 2k^3+g.
\end{equation}

By Lemma~\ref{nu1}(c) and (e), for all integers $\beta,\gamma\ge2\ell$,
\begin{equation}\label{betagama}
|\nu(\beta)-\nu(\gamma)|\le (k-1)|\beta-\gamma|.
\end{equation}
Since both, $\alpha-r_i$ and $\alpha-r'_i$, are at least $k+\ell-(k-\ell)=2\ell$, the triangle
inequality together with~\eqref{betagama} gives, for every $i=1,\ldots,a+b$,
\begin{equation}\label{nurrprim}
\left|
\bigl(\nu(\alpha-r_i)+r_i\bigr)
-\bigl(\nu(\alpha-r'_i)+r'_i\bigr)
\right|
\le k|r_i-r'_i|.
\end{equation}
Consequently, by the definitions of $\Phi(\alpha,a+b)$, $\Phi(\alpha,a)$, and $\Phi(\alpha,b)$,
\begin{align*}
\Phi(\alpha,a+b)
&\ge \sum_{i=1}^{a+b}\bigl(\nu(\alpha-r_i)+r_i\bigr)\\
&\overset{\eqref{rrprim},\eqref{nurrprim}}{\ge} \sum_{i=1}^{a+b}\bigl(\nu(\alpha-r'_i)+r'_i\bigr)
     -k(2k^3+g)\\
&=\Phi(\alpha,a)+\Phi(\alpha,b)-k(2k^3+g).
\end{align*}

This proves the lemma with $c=k(2k^3+g)$.
\end{proof}

\end{document}